\documentclass[11pt]{article}
\usepackage{amsmath,amsfonts,amssymb,amsthm,amscd}

\title{Reconstruction of Manifolds in Noncommutative Geometry}

\author{Adam Rennie\dag\ddag
\thanks{email: \texttt{adam.rennie@maths.anu.edu.au},
\texttt{varilly@cariari.ucr.ac.cr}}
\word{and}%
Joseph C. V\'arilly\P$^*$ \\[6pt]
\dag Institute for Mathematical Sciences,
University of Copenhagen\\
Universitetsparken 5, DK-2100 Copenhagen, Denmark\\[6pt]
\ddag Mathematical Sciences Institute,
Australian National University,\\
Canberra, ACT 0200, Australia\\[6pt]
\P Departamento de Matem\'aticas,
Universidad de Costa Rica, \\
2060 San Jos\'e, Costa Rica}

\date{31 January 2008}

\advance\topmargin by -\headheight
\advance\topmargin by -\headsep
\evensidemargin=\oddsidemargin
\makeatletter
\def\section{\@startsection{section}{1}{\z@}{-3.5ex plus -1ex minus
  -.2ex}{2.3ex plus .2ex}{\large\bf}}
\def\subsection{\@startsection{subsection}{2}{\z@}{-3.25ex plus -1ex
  minus -.2ex}{1.5ex plus .2ex}{\normalsize\bf}}
\makeatother

\numberwithin{equation}{section} 

\theoremstyle{plain} 
\newtheorem{thm}{Theorem}[section]
\newtheorem{lemma}[thm]{Lemma}
\newtheorem{prop}[thm]{Proposition}
\newtheorem{corl}[thm]{Corollary}

\theoremstyle{definition} 
\newtheorem{defn}{Definition}[section]
\newtheorem{cond}{Condition} 

\newtheorem{rmk}[thm]{Remark} 

\DeclareMathOperator{\Ad}{Ad}     
\DeclareMathOperator{\ad}{ad}     
\DeclareMathOperator{\Cliff}{{\C\ell}} 
\DeclareMathOperator{\Dom}{Dom}   
\DeclareMathOperator{\End}{End}   
\DeclareMathOperator{\Hom}{Hom}   
\DeclareMathOperator{\Id}{Id}     
\DeclareMathOperator{\intr}{Int}  
\DeclareMathOperator{\linspan}{span} 
\DeclareMathOperator{\rank}{rank} 
\DeclareMathOperator{\spec}{sp}   
\DeclareMathOperator{\supp}{supp} 
\DeclareMathOperator{\Tr}{Tr}     
\DeclareMathOperator{\tr}{tr}     
\DeclareMathOperator{\tsum}{{\textstyle\sum}} 
\DeclareMathOperator{\Vol}{Vol}   
\DeclareMathOperator{\vol}{vol}   

\newcommand{\al}{\alpha}      
\newcommand{\bt}{\beta}       
\newcommand{\dl}{\delta}      
\newcommand{\eps}{\varepsilon} 
\newcommand{\Ga}{\Gamma}      
\newcommand{\ga}{\gamma}      
\newcommand{\La}{\Lambda}     
\newcommand{\la}{\lambda}     
\newcommand{\sg}{\sigma}      
\newcommand{\vf}{\varphi}     

\newcommand{\A}{\mathcal{A}}  
\newcommand{\B}{\mathcal{B}}  
\newcommand{\C}{\mathbb{C}}   
\newcommand{\cc}{\mathbf{c}}  
\newcommand{\D}{\mathcal{D}}  
\newcommand{\E}{\mathcal{E}}  
\renewcommand{\H}{\mathcal{H}}  
\renewcommand{\L}{\mathcal{L}} 
\newcommand{\N}{\mathbb{N}}   
\newcommand{\R}{\mathbb{R}}   
\newcommand{\Sf}{\mathbb{S}}  
\newcommand{\Z}{\mathbb{Z}}   

\newcommand{\ac}{\mathrm{ac}} 
\newcommand{\del}{\partial}   
\newcommand{\hookto}{\hookrightarrow} 
\newcommand{\less}{\setminus} 
\newcommand{\ol}{\overline}   
\newcommand{\ox}{\otimes}     
\newcommand{\w}{\wedge}       
\newcommand{\x}{\times}       
\newcommand{\8}{\bullet}      
\renewcommand{\.}{\cdot}      
\renewcommand{\:}{\colon}     

\newcommand{\as}{\quad\mbox{as}\enspace} 
\newcommand{\CDA}{\mathcal{C_D(A)}} 
\newcommand{\Coo}{C^\infty}   
\renewcommand{\d}{\underline{\mathrm{d}}} 
\newcommand{\Dhat}{\widehat{\mathcal{D}}} 
\newcommand{\Dreg}{\langle\D\rangle} 
\newcommand{\Dslash}{{D\mkern-11.5mu/\,}} 
\newcommand{\Gaoo}{\Gamma_\infty} 
\newcommand{\half}{\tfrac{1}{2}} 
\newcommand{\hatox}{\mathrel{\widehat\otimes}} 
\newcommand{\omlim}{\mathop{\omega\mbox{-lim}}\limits} 
\newcommand{\op}{\circ}       
\newcommand{\otto}{\leftrightarrow} 
\newcommand{\oxyox}{\otimes\cdots\otimes} 
\newcommand{\Shat}{\widehat{S}} 
\newcommand{\Trw}{\Tr_\omega} 
\newcommand{\wyw}{\wedge\cdots\wedge} 

\newcommand{\bbraket}[2]{\langle\!\langle#1\mathbin{|}
                          #2\rangle\!\rangle} 
\newcommand{\bigpairing}[2]{\bigl(#1\mathbin{\big|}#2\bigr)} 
\newcommand{\braket}[2]{\langle#1\mathbin{|}#2\rangle} 
\newcommand{\hideqed}{\renewcommand{\qed}{}} 
\newcommand{\pairing}[2]{(#1\mathbin{|}#2)} 
\newcommand{\piso}[1]{\lfloor#1\rfloor} 
\newcommand{\row}[3]{{#1}_{#2},\dots,{#1}_{#3}} 
\newcommand{\set}[1]{\{\,#1\,\}}  
\newcommand{\twobytwo}[4]{\begin{pmatrix}#1 & #2 \\
                           #3 & #4\end{pmatrix}} 
\newcommand{\word}[1]{\quad\mbox{#1}\quad} 

\newbox\ncintdbox \newbox\ncinttbox 
	\setbox0=\hbox{$-$}
	\setbox2=\hbox{$\displaystyle\int$}
	\setbox\ncintdbox=\hbox{\rlap{\hbox
		to \wd2{\hskip-.125em \box2\relax\hfil}}\box0\kern.1em}
	\setbox0=\hbox{$\vcenter{\hrule width 4pt}$}
	\setbox2=\hbox{$\textstyle\int$}
	\setbox\ncinttbox=\hbox{\rlap{\hbox
		to \wd2{\hskip-.175em \box2\relax\hfil}}\box0\kern.1em}
\newcommand{\ncint}{\mathop{\mathchoice{\copy\ncintdbox}%
					{\copy\ncinttbox}{\copy\ncinttbox}%
					{\copy\ncinttbox}}\nolimits}

\begin{document}

\maketitle

\vspace{-8pt}

\begin{abstract}
We show that the algebra $\A$ of a commutative unital spectral triple
$(\A,\H,\D)$ satisfying several additional conditions, slightly
stronger than those proposed by Connes, is the algebra of smooth
functions on a compact spin manifold.
\end{abstract}

\tableofcontents

\section{Introduction}
\label{sec:intro}

Noncommutative Geometry, as developed over the past several years by
Connes and coworkers, has produced a profusion of examples of
``noncommutative spaces'' \cite{ConnesMaPardis}, many of which partake
of the characteristics of smooth Riemannian manifolds, whose metric
and differential structure is determined by a generalized Dirac
operator. To find a common framework for those examples, Connes
proposed in \cite{ConnesGrav} an axiomatic framework for
``noncommutative spin manifolds''.

The geometry is carried by the notion of spectral triple $(\A,\H,\D)$;
the familiar Riemannian spin geometry is recovered when $\A$ is a
coordinate algebra of smooth functions on a manifold, $\H$ a Hilbert
space of spinors, and $\D$ the Dirac operator determined by the spin
structure and Riemannian metric. The question of reconstruction is
whether the operator-theoretic framework proposed by Connes, or some
variation of it, suffices to determine this spin manifold structure
whenever the algebra $\A$ is commutative.

In \cite{ConnesGrav}, Connes held out the hope that it could be so;
but the extraction of a manifold from these postulates has proved
elusive. In \cite{RennieComm}, a first attempt at doing so was
presented, but was subsequently shown to fall short of the goal
\cite{Gorokhovsky}. A detailed description of the reconstruction of
many geometric features of a spin manifold was presented in
\cite{Polaris}, but there the starting algebra $\A$ was assumed
\emph{a priori} to be the smooth functions on a compact manifold.

In this paper, using a slightly stronger set of conditions on a
spectral triple, we show that from the further assumption of a
commutative coordinate algebra $\A$ one can indeed recover a compact
boundaryless manifold whose smooth functions coincide with~$\A$.

One of the key themes of the axioms proposed by Connes was Poincar\'e
duality in $K$-theory. Earlier results of Sullivan \cite{Sullivan}
indicated that in high dimensions, in the absence of $2$-torsion and
in the simply connected setting, Poincar\'e duality in $KO$-theory
characterizes the homotopy type of a compact manifold. While as a
guiding principal such an idea is very attractive, we have not found a
way to implement this approach to reconstruct a manifold.

Instead, we utilize an earlier formulation of Poincar\'e duality
in noncommutative geometry which is phrased at the level of Hochschild
chains, and thus is more useful for elaborating a proof. This concrete
version of Poincar\'e duality is described by a ``closedness
condition'' \cite[VI.4.$\ga$]{Book}, which historically arose from
attempts to fine-tune the Lagrangian of the Standard Model of
elementary particles, and conceptually is an analogue of Stokes'
theorem.

Poincar\'e duality in $K$-theory plays no role in our reconstruction
of a manifold as a compact space $X$ with charts and smooth transition
functions. However, once that has been achieved, it is needed to show
that $X$ carries a spin$^c$ structure and to identify the class of
$(\A,\H,\D)$ as the fundamental class of the spin$^c$ manifold. The
key to this is Plymen's characterization of spin$^c$ structures
\cite{Plymen} as Morita equivalence bimodules for the Clifford action
induced by the metric. Indeed, it would be economical to replace
Poincar\'e duality by postulating instead the existence of such
bimodules; we touch on this in our final section.

Compactness of the manifold, or equivalently the condition that the
coordinate algebra have a unit, is an essential technical feature of
our proof. However, the reconstruction of noncompact manifolds should
also be possible, under some alternative conditions along the lines
suggested in~\cite{Himalia,RennieSmooth}. Indeed, many of the crucial
arguments used in reconstructing the coordinate charts are completely
local.

The proof that the Gelfand spectrum $X = \spec(\A)$ is a differential
manifold is quite long, but may be conceptually broken into two steps.
The first is to construct a vector bundle over~$X$ which plays the
role of the cotangent bundle. Already at this stage we need to deploy
all the conditions on our spectral triple (except Poincar\'e duality
in $K$-theory and a metric condition). In particular, we identify
local trivializations and bases of this bundle in terms of the
`$1$-forms' given by the orientability condition. These $1$-forms
$[\D,a^j_\al]$, for $j = 1,\dots,p$, $\al = 1,\dots,n$, generate the
sections of this bundle, and the aim now is to show that the maps
$a_\al = (a^1_\al,\dots,a^p_\al) : X \to \R^p$ provide coordinates on
suitable open subsets of~$X$.

This is accomplished by proving that $a_\al$ is locally one-to-one and
open. The tools used here are a Lipschitz functional calculus, some
measure theoretic results of Voiculescu~\cite{Voiculescu}, some basic
point set topology and properties of the map~$a_\al$, and finally
the unique continuation properties for Dirac-type operators
\cite{BoossMW,Kim}.

The main tools in the proof are a multivariate $\Coo$ functional
calculus for regular spectral triples \cite{RennieSmooth}, which we
present here; as well as a Lipschitz functional calculus. The first 
of these enables us to construct partitions of unity and local
inverses within the algebra~$\A$.

\vspace{6pt}

The plan of the paper is as follows. In Section \ref{sec:defns} we
give some standard definitions and background results, including the
$\Coo$ functional calculus and its immediate consequences. In Section
\ref{sec:conds}, we introduce the conditions on a spectral triple
needed to establish our main result.

Section \ref{sec:cotg-bdl} details the construction of the cotangent
bundle, while Section \ref{sec:lip-open} develops a Lipschitz
functional calculus needed to deal with the topology of our coordinate
charts. Sections \ref{sec:point-sets} and~\ref{sec:spec-mfld} contain
the detailed proof that we do indeed recover a manifold. We develop
the necessary point set topology to establish that the spectrum of our
algebra is a manifold: the main issue is the absence of branch points
in the chart domains. We show that the algebra generated by $\A$ and
$[\D,\A]$ is locally a direct sum of Clifford actions arising from one
or several Riemannian metrics, for which $\D$ is (again, locally) a
direct sum of Dirac-type operators. Then we use the unique
continuation properties of Dirac-type operators and the local
description of $\D$ to banish any branch points and thereby get a
manifold. That done, we assemble the Clifford actions globally, and
so produce the Clifford action of a single Riemannian metric.

In Section \ref{sec:PD} we explain in some detail how the (unique)
spin$^c$ structure arises from Poincar\'e duality in $K$-theory.
The Dirac operator is shown to differ from $\D$ by at most an
endomorphism of the corresponding spinor bundle.

Section \ref{sec:more-conds} collects some further remarks on our
postulates and their possible variants.

Appendix~A establishes some basic results about Hermitian pairings on
finite projective modules. Appendix~B examines additional results
about our conditions, in particular the redundancy of the metric
condition.

\subsubsection*{Acknowledgments}
This work has profited from discussions with Alan Carey, Alain Connes,
Nigel Higson, Steven Lord, Ryszard Nest, and Iain Raeburn. JCV is
grateful to Iain Raeburn for warm hospitality at the University of
Newcastle and to Ryszard Nest for a timely visit to Copenhagen. AR
thanks JCV for unique hospitality whilst visiting Costa Rica. This
work was supported by an ARC grant, DP0211367, by the SNF, Denmark, by
a University of Newcastle Visitor Grant, and by the European
Commission grant MKTD--CT--2004--509794 at the University of Warsaw.
Support from the Universidad de Costa Rica is also acknowledged. Both
authors would especially like to thank Piotr Hajac for generous
hospitality at IMPAN during the course of this work.

\section{Spectral triples and smooth functional calculus}
\label{sec:defns}

The central notion of this paper is that of a spectral triple
\cite{ConnesSpec} over a commutative algebra. We begin by recalling
several basic definitions, in order to establish a suitable functional
calculus for them.

\begin{defn}
\label{df:spec-tri}
A \textit{spectral triple} $(\A,\H,\D)$ is given by:
\begin{enumerate}
\item
A faithful representation $\pi\: \A \to \B(\H)$ of a unital
$*$-algebra $\A$ by bounded operators on a Hilbert space~$\H$; and
\item
A selfadjoint operator $\D$ on~$\H$, with dense domain $\Dom\D$, such
that for each $a \in \A$, $[\D,\pi(a)]$ extends to a bounded operator
on $\H$ and $\pi(a)(1 + \D^2)^{-1/2}$ is a compact operator.
\end{enumerate}
The spectral triple is said to be \textit{even} if there is an
operator $\Ga = \Ga^* \in \B(\H)$ such that $\Ga^2 = 1$
(this determines a $\Z_2$-grading on~$\H$), for which
$[\Ga, \pi(a)] = 0$ for all $a \in \A$ and
$\Ga\D + \D\Ga = 0$ (i.e., $\pi(\A)$ is even and $\D$ is odd
with respect to the grading). If no such grading is available, the
spectral triple is called \textit{odd}.
\end{defn}

\begin{rmk}
Since $\A$ is faithfully represented on~$\H$, we may and shall
omit~$\pi$, regarding $\A$ as a subalgebra of $\B(\H)$. As such, its
norm closure $\ol{\A} = A$ is a $C^*$-algebra.
\end{rmk}

\begin{rmk}
In this paper, we shall always assume that $\A$ is \emph{unital}.
Nonunital spectral triples have been studied in
\cite{RennieSmooth,RennieSumm} under the assumption that $\A$ has a
dense ideal with local units. Another class of nonunital spectral
triples are those arising from Moyal products, analyzed in detail in
\cite{Himalia} (and anticipated in~\cite{Selene}). The Moyal example
shows that it is important to treat a certain unitization of $\A$ as
part of the data of a (nonunital) spectral triple, so that it is
proper to focus first on the unital case.
\end{rmk}

\begin{defn}
\label{df:qc-infty}
The operator $\D$ gives rise to two (commuting) derivations of
operators on~$\H$; we shall denote them by
$$
\d x := [\D,x],  \qquad  \dl x := [|\D|,x],  \word{for} x \in \B(\H).
$$
According to Definition~\ref{df:spec-tri}, $\A$ lies within
$\Dom\d := \set{x \in \B(\H) : [\D,x] \in \B(\H)}$.

A \textit{spectral triple} $(\A,\H,\D)$ is called $QC^\infty$ if
$$
\A \cup \d\A  \subseteq  \bigcap_{m=1}^\infty \Dom \dl^m.
$$
\end{defn}

\begin{rmk}
The terminology $QC^\infty$ was introduced in \cite{CareyPRSone}, to
distinguish ``quantum'' differentiability of operators from
``classical'' differentiability of smooth functions. One can also
define $QC^k$, for $k \in \N$, by requiring only that
$\A \cup \d\A \subseteq \Dom \dl^m$ for $m = 1,\dots,k$. Such spectral
triples are more often referred to as \textit{regular}
\cite{ConnesGrav,Polaris}, and have been called \textit{smooth}
in~\cite{RennieSmooth}.
\end{rmk}

\begin{defn}
\label{df:delta-top}
If $(\A,\H,\D)$ is a $QC^\infty$ spectral triple, the family of
seminorms
\begin{equation}
q_m(a) := \|\dl^m a\|  \word{and}
q'_m(a) := \|\dl^m([\D,a])\|,  \quad  m = 0,1,2,\dots
\label{eq:dl-snorms}
\end{equation}
determine a locally convex topology on~$\A$ which is finer than the
norm topology of~$A$ (that is given by $q_0$ alone) and in which the
involution $a \mapsto a^*$ is continuous. Let $\A_\dl$ denote the
completion of~$\A$ in the topology of~\eqref{eq:dl-snorms}.
\end{defn}

We quote Lemma~16 of \cite{RennieSmooth}.

\begin{lemma}
\label{lm:smo}
Let $(\A,\H,\D)$ be a $QC^\infty$ spectral triple. The Fr\'echet
algebra $\A_\dl$ is a pre-$C^*$-algebra, and $(\A_\dl,\H,\D)$ is also
a $QC^\infty$ spectral triple.
\qed
\end{lemma}

Recall that a pre-$C^*$-algebra is a dense subalgebra of a
$C^*$-algebra which is stable under the holomorphic functional
calculus of that $C^*$-algebra. There is little loss of generality in
assuming that $\A$ is complete in the topology given by
\eqref{eq:dl-snorms}, thus is a \textit{Fr\'echet pre-$C^*$-algebra},
and we shall do so. This condition guarantees that the spectrum of an
element $a \in \A$ coincides with its spectrum in the
$C^*$-completion~$A$, and that any character of the pre-$C^*$-algebra
$\A$ extends to a character of $A$ as well. We shall denote the
character space by $X := \spec(\A) = \spec(A)$.

Moreover, when $\A$ is a Fr\'echet pre-$C^*$-algebra, so also is the
algebra $M_n(\A)$ of $n \x n$ matrices with entries in~$A$, whose
$C^*$-completion is $M_n(A)$; for a proof, see \cite{Schweitzer}. By a
theorem of Bost \cite{Bost,Polaris}, the (topological) $K$-theories of
$\A$ and~$A$ coincide: $K_i(\A) = K_i(A)$ for $i = 0,1$.

By replacing any seminorm $q$ by $a \mapsto q(a) + q(a^*)$ if
necessary, we may suppose that $q(a) = q(a^*)$ for all $a \in \A$.
We note in passing that the multiplication in the Fr\'echet algebra
$\A$ is jointly continuous \cite{Mallios}.

\begin{lemma}
\label{lm:dense-proj}
Let $\A$ be a Fr\'echet pre-$C^*$-algebra and let $A$ be its
$C^*$-completion. If $\tilde q \in A$ is a projector (i.e., a
selfadjoint idempotent), and if $0 < \eps < 1$, then there is a
projector $q \in \A$ such that $\|q - \tilde q\| < \eps$.
\end{lemma}

\begin{proof}
Choose $\dl$ with $0 < \dl < \eps/32$, and let $b = b^* \in \A$
be such that $\|b - \tilde q\| < \dl$. Observe that
$$
\|b^2 - b\| = \|b^2 - \tilde q^2 + \tilde q - b\|
\leq (\|b + \tilde q\| + 1)\,\|b - \tilde q\|
\leq (3 + \dl)\|b - \tilde q\| < \dl(3 + \dl) < 4\dl.
$$
Since $\A$ is a Fr\'echet pre-$C^*$-algebra, one may, provided $\dl$
is sufficiently small, use holomorphic functional calculus to
construct a homotopy within $\A$ from $b$ to $e \in \A$ such that
$e^2 = e$ and $\|e - b\| < 2\|b^2 - b\|$; see
\cite[Lemma~3.43]{Polaris}, for instance.

Let $q := ee^* (ee^* + (1 - e^*)(1 - e))^{-1}$. Then $q$ is a
projector in $A$, and it lies in $\A$ since $ee^* + (1 - e^*)(1 - e)$
is invertible in the pre-$C^*$-algebra $\A$. By taking $A$ to be
faithfully represented on a Hilbert space $\H$, we can write $e$, $q$
and $b$ as operators on $\H = e\H \oplus (1-e)\H$, as follows:
$$
e = \twobytwo{1}{T}{0}{0},    \qquad
q = \twobytwo{1}{0}{0}{0},  \qquad
b = \twobytwo{R}{V}{V^*}{S},
$$
with $R$, $S$ selfadjoint and $V,T \: (1 - e)\H \to e\H$ bounded. Then
$\|e - b\| < 8\dl < \eps/4$ means $\|(e - b)^*(e - b)\| < \eps^2/16$,
which entails
$$
\|(R - 1)^2 + VV^*\| < \frac{\eps^2}{16},  \qquad
\|(V - T)^*(V - T) + S^2\| < \frac{\eps^2}{16},
$$
so that $\|V\| < \eps/4$, $\|V - T\| < \eps/4$, and therefore
$\|q - e\| = \|T\| < \eps/2$. Hence
$$
\|q - \tilde q\| \leq \|q - e\| + \|e - b\| + \|b - \tilde q\|
< \frac{\eps}{2} + \frac{\eps}{4} + \dl < \eps.
\eqno \qed
$$
\hideqed
\end{proof}

A $QC^\infty$ spectral triple $(\A,\H,\D)$ for which $\A$ is
complete has not only a holomorphic functional calculus for
$\A$, but also a \textit{$\Coo$ functional calculus} for
selfadjoint elements: we quote \cite[Prop.~22]{RennieSmooth}.

\begin{prop}[$\Coo$ Functional Calculus]
\label{pr:cinfty}
Let $(\A,\H,\D)$ be a $QC^\infty$ spectral triple, and suppose $\A$ is
complete. Let $f \: \R \to \C$ be a $\Coo$ function in a
neighbourhood of the spectrum of $a = a^* \in \A$. If we define
$f(a) \in A$ using the continuous functional calculus, then in fact
$f(a)$ lies in~$\A$.
\qed
\end{prop}

\begin{rmk}
For each $a = a^* \in \A$, the $\Coo$-functional calculus defines
a continuous homomorphism $\Psi \: \Coo(U) \to \A$, where
$U \subset \R$ is any open set containing the spectrum of $a$, and the
topology on $\Coo(U)$ is that of uniform convergence of all
derivatives on compact subsets.
\end{rmk}

The following proposition extends this result to the case of smooth
functions of several variables, yielding a \textit{multivariate
$\Coo$ functional calculus}. Before stating it, we recall the
continuous functional calculus for a finite set $\row{a}{1}{n}$ of
commuting selfadjoint elements of a unital $C^*$-algebra~$A$. These
generate a unital $*$-algebra whose closure in~$A$ is a
$C^*$-subalgebra $C^*\langle 1,\row{a}{1}{n}\rangle$; let $\Delta$ be
its (compact) space of characters. Evaluation of polynomials
$p \mapsto p(\row{a}{1}{n})$ yields a surjective morphism from
$C\bigl(\prod_{j=1}^n \spec a_j \bigr)$ onto
$C^*\langle 1,\row{a}{1}{n}\rangle \simeq C(\Delta)$ which
corresponds, via the Gelfand functor, to a continuous injection
$\Delta \hookto \prod_{j=1}^n \spec a_j$; this joint spectrum $\Delta$
may thus be regarded as a compact subset of~$\R^n$. If
$h \in C(\Delta)$, we may define $h(\row{a}{1}{n})$ as the image of
$h|_\Delta$ in $C^*\langle 1,\row{a}{1}{n}\rangle$ under the Gelfand
isomorphism.

\begin{prop}
\label{pr:mult-cinfty}
Let $(\A,\H,\D)$ be a $QC^\infty$ spectral triple. Let $\row{a}{1}{n}$
be mutually commuting selfadjoint elements of~$\A$, and let
$\Delta \subset \R^n$ be their joint spectrum. Let $f\: \R^n \to \C$
be a $\Coo$ function supported in a bounded open neighbourhood $U$
of~$\Delta$. Then $f(\row{a}{1}{n})$ lies in~$\A_\dl$.
\end{prop}

\begin{proof}
We first define the operator $f(\row{a}{1}{n})$ lying in~$A$, the
$C^*$-completion of~$\A$, using the continuous functional calculus.

Since $f$ is a compactly supported smooth function on~$\R^n$, we may
alternatively define $f(\row{a}{1}{n}) \in A$ by a Fourier integral:
\begin{equation}
f(\row{a}{1}{n}) = (2\pi)^{-n/2} \int_{\R^n} \hat f(\row{s}{1}{n})
\exp(i\,s\.a) \,d^n s,
\label{eq:Fourier-calc}
\end{equation}
where $s\.a = s_1 a_1 +\cdots+ s_n a_n$. Since $\dl$ (and likewise
$\d = \ad\D$) is a norm-closed derivation from $\A$ to $\B(\H)$, we
may conclude that $f(\row{a}{1}{n}) \in \Dom\dl$ with
\begin{equation}
\dl(f(\row{a}{1}{n})) = (2\pi)^{-n/2} \int_{\R^n}
\hat f(\row{s}{1}{n}) \dl(\exp(i\,s\.a)) \,d^n s,
\label{eq:delta-f}
\end{equation}
provided we can establish dominated convergence for the integral on
the right hand side \cite{BratteliRoI}. Just as in the one-variable
case \cite{RennieSmooth}, since each $a_j \in \Dom\dl$, we find that
$\exp(i\,s\.a) = \prod_j \exp(i s_j a_j)$ lies in $\Dom\dl$ also: its
factors are given by the expansion
\begin{equation}
\dl(\exp(i s_j a_j)) = i s_j \int_0^1
\exp(its_j a_j)\, \dl(a_j) \,\exp(i(1-t)s_j a_j) \,dt,
\label{eq:iter-delta}
\end{equation}
and in particular,
$$
\|\dl(\exp(i\,s\.a))\| \leq C \tsum_j |s_j|,  \qquad
C = \max_j \biggl( \|\dl(a_j)\| \prod_{i\neq j} \|a_i\| \biggr).
$$
A norm bound which dominates the right hand side of~\eqref{eq:delta-f}
is thus given by
$$
\int_{\R^n} |\hat f(\row{s}{1}{n})|\, \|\dl(\exp(i\,s\.a))\| \,d^n s
\leq C \sum_{j=1}^n (2\pi)^{-n/2} \int_{\R^n}
|\hat f(\row{s}{1}{n})|\,|s_j| \,d^n s.
$$

Let $A_0$ be the completion of~$\A$ for the norm
$\|a\|_\D := \|a\| + \|\d a\|$; notice that $A_0 \subseteq A$.
Replacing $\dl$ by~$\d$ in the previous argument, we find that
$$
\|f(\row{a}{1}{n})\|_\D \leq \|\hat f\|_1 + \|\d a\|
\sum_{j=1}^n (2\pi)^{-n/2} \int_{\R^n}
|\hat f(\row{s}{1}{n})|\,|s_j| \,d^n s.
$$
Therefore, $f(\row{a}{1}{n})$ can be approximated, in the $\|\.\|_\D$
norm, by Riemann sums for~\eqref{eq:Fourier-calc} belonging to~$\A$,
and thus $f(\row{a}{1}{n}) \in A_0$.

Since $\dl$ and $\d$ are commuting derivations, we obtain that
$\dl(f(\row{a}{1}{n})) \in \Dom\d$ and
$\d(f(\row{a}{1}{n})) \in \Dom\dl$ for $a \in \A$, and
$\|\dl(\d(f(\row{a}{1}{n})))\|$ is bounded by a linear combination of
expressions $\|\dl(\d a_j)\| \int|\hat f(\row{s}{1}{n})|\,|s_j|\,d^ns$
and
$\|\dl a_j\|\,\|\d a_k\| \int|\hat f(\row{s}{1}{n})|\,|s_js_k|\,d^ns$.
In particular, $\|\dl(f(\row{a}{1}{n}))\|_\D$ also has a bound of this
type.

For each $m = 1,2,3,\dots$, let $A_m$ be the completion of~$\A$ for
the norm $\sum_{k\leq m} \|\dl^k(a)\|_\D$. Then $\dl$ extends to
a norm-closed derivation from $A_m$ to $\B(\H)$, and an ugly but
straightforward induction on~$m$ shows that each
$\dl^k(f(\row{a}{1}{n}))$ and $\dl^k(\d(f(\row{a}{1}{n})))$ lies in
its domain, using the convergence of
$\int |\hat f(\row{s}{1}{n})|\,|p(\row{s}{1}{n})|\,d^n s$ for $p$ a
polynomial of degree $\leq m + 1$. Thus, $f(\row{a}{1}{n}) \in A_m$.
Since $\A_\dl = \bigcap_{m\geq 0} A_m$, we conclude that
$f(\row{a}{1}{n}) \in \A_\dl$.
\end{proof}

\begin{rmk}
For most of this paper, the spectral triples we consider will be
commutative and will satisfy the first order property
\cite[VI.4.$\ga$]{Book}, meaning that $[[\D,a],b] = 0$ for all
$a,b \in \A$. In such cases, to define $[\D,f(\row{a}{1}{n})]$, we may
note that for any polynomial~$p$, the first-order property allows us
to write
\begin{equation}
[\D, p(\row{a}{1}{n})]
= \sum_{j=1}^n \del_j p(\row{a}{1}{n}) \, [\D,a_j],
\label{eq:Dcomm-poly}
\end{equation}
with $\del_j p$ being the $j$-th partial derivative of~$p$. By a
$C^1$-approximation argument ---see Proposition~\ref{pr:Nachbin}
below--- we obtain
$[\D,f(\row{a}{1}{n})] = \sum_j \del_j f(\row{a}{1}{n}) \,[\D,a_j]$
for any $f$ satisfying the hypotheses of
Proposition~\ref{pr:mult-cinfty}. Since $(\A,\H,\D)$ is $QC^\infty$,
one sees immediately that the right hand side of~\eqref{eq:Dcomm-poly}
belongs to the smooth domain of~$\dl$.
\end{rmk}

We now use the $\Coo$ functional calculus to prove the existence
of certain elements of $\A$, where $(\A,\H,\D)$ is a (unital)
commutative $QC^\infty$ spectral triple. The algebra elements we are
looking for are smooth partitions of unity and local inverses.

\begin{lemma}
\label{lm:partn-unity}
Let $(\A,\H,\D)$ be a $QC^\infty$ spectral triple where $\A$ is
commutative and complete. Let $\set{U_\al : \al = 1,\dots,n}$ be any
finite open cover of the compact Hausdorff space $X = \spec(\A)$. Then
there exist $\phi_\al \in \A$, for $\al = 1,\dots,n$, such that
\begin{equation}
\supp \phi_\al \subseteq U_\al,  \quad  0 \leq \phi_\al \leq 1,
\word{and} \sum_{\al=1}^n \phi_\al = 1.
\label{eq:partn-unity}
\end{equation}
\end{lemma}

\begin{proof}
Since $\A = \A_\dl$ is a pre-$C^*$-algebra, its character space is
the same as that of its $C^*$-completion, $A$; thus $X = \spec(A)$ is
a compact Hausdorff space. Now, $X$ always admits a continuous
partition of unity~\cite{Dydak} subordinate to the cover $\{U_\al\}$,
i.e., we can find $\row{\tilde\phi}{1}{n} \in A$
satisfying~\eqref{eq:partn-unity}. Let $p \in M_n(A)$ be the matrix
whose $(\al,\bt)$-entry is $(\tilde\phi_\al \tilde\phi_\bt)^{1/2}$;
then $p$ is a projector, that is, a selfadjoint idempotent:
$p^2 = p = p^*$.

Since $M_n(\A)$ is a Fr\'echet pre-$C^*$-algebra, it is known
\cite{Bost,Polaris} that the inclusion $M_n(\A) \hookto M_n(A)$
induces a homotopy equivalence between the respective sets of
idempotents in these algebras. Thus, there is a norm-continuous path
of idempotents $t \mapsto e_t = e_t^2 \in M_n(A)$ linking $p = e_0$ to
an idempotent $e_1 \in M_n(\A)$. Moreover, such a path may be chosen
so that each $\|p - e_t\| < \eps$ for a preassigned $\eps > 0$
\cite[Lemma~3.43]{Polaris}. We choose $\eps < 1/3n$. Replacing $e_t$ by
$q_t := e_te_t^* (e_te_t^* + (1 - e_t^*)(1 - e_t))^{-1}$, we may link
$p$ to $q = q_1$ by a path of projectors in $M_n(A)$; by the proof of
Lemma~\ref{lm:dense-proj}, we obtain $\|p - q_t\| < 3\eps < 1/n$ for
$0 \leq t \leq 1$. Since the positive element
$e_1e_1^* + (1-e_1^*)(1-e_1)$ is invertible in the pre-$C^*$-algebra
$M_n(\A)$, $q$ lies in $M_n(\A)$.

By the Serre--Swan theorem~\cite{SwanVect}, the projectors $p$ and~$q$
define vector bundles over~$X$ of the same rank: the rank is given by
the matrix trace $\tr p = \tr q$, a locally constant integer-valued
function in $A = C(X)$. Write $\psi_\al := q_{\al\al} \in \A$, and
notice that $\psi_\al \geq 0$ in~$A$; then
\begin{equation}
\sum_{\al=1}^n \psi_\al = \tr q = \tr p
= \sum_{\al=1}^n \tilde\phi_\al = 1.
\label{eq:part}
\end{equation}

We now modify the elements $\psi_\al$ to obtain a partition of unity
subordinate to the cover $\{U_\al\}$. By construction,
$\|\psi_\al - \tilde\phi_\al\| \leq \|q - p\| < 1/n$ for each~$\al$.
Choose a smooth function $g \: \R \to [0,1]$ with support in
$[\eps, 1+\eps]$ such that $0 < g(t) \leq 1$ for $\eps < t \leq 1$.
Then define
$V_\al := \set{x \in X : \psi_\al(x) > \eps} \subset U_\al$. Setting
$\chi_\al(x) := g(\psi_\al(x))$ gives $\chi_\al > 0$ on $V_\al$ and
$\supp \chi_\al \subset U_\al$. For all $x \in X$, there is some
$\chi_\bt$ with $\chi_\bt(x) > 0$: for if not, then
$\psi_\bt(x) \leq \eps$ for each~$\bt$, and
$\sum_\bt \psi_\bt(x) \leq n\eps < 1$, contradicting~\eqref{eq:part}.
We now define $\phi_\al := \chi_\al \big/ \sum_\bt \chi_\bt$, which
clearly satisfies \eqref{eq:partn-unity}. Since
$\chi_\al = g(\psi_\al)$, Proposition~\ref{pr:cinfty} shows that
$\chi_\al \in \A$. Moreover, $\sum_\bt \chi_\bt$ is invertible
in~$\A$, and hence $\phi_\al \in \A$ for each~$\al$, as required.
\end{proof}

\begin{corl}
\label{cr:partn-unity}
Given a $QC^\infty$ spectral triple $(\A,\H,\D)$ where $\A$ is
commutative and complete, let $K \subset U \subset \spec(\A)$ where
$K$ is compact and $U$ is open. Then there is some $\psi \in \A$ such
$0 \leq \psi \leq 1$, $\psi \equiv 1$ on~$K$, and $\psi \equiv 0$
outside~$U$.
\end{corl}

\begin{proof}
There is a partition of unity $\{\psi, 1 - \psi\}$ subordinate to the
open cover $\{U, \spec(\A) \less K\}$, with $\psi \in \A$.
\end{proof}

\begin{lemma}
\label{lm:local-unit}
Let $(\A,\H,\D)$ again be a $QC^\infty$ spectral triple with $\A$
commutative and complete. Let $a \in \A$ have compact support
contained in an open subset $U \subset \spec(\A)$. Then there exists
$\phi \in \A$ such that $\phi a = a$ and $\supp\phi \subset U$.
\end{lemma}

\begin{proof}
Choose $b \in A = C(X)$, by Urysohn's lemma, such that
$0 \leq b \leq 1$, $b(x) = 1$ for $x \in \supp a$ and $b(x) = 0$ for
$x \notin U$. Pick $\dl \in (0,\half)$ and choose $\psi \in \A$
satisfying $\|b - \psi\| < \half\dl$. Let $f \: \R \to [0,1]$ be a
compactly supported smooth function such that $f(t) = 0$ for
$t \leq \dl$ and $f(t) = 1$ for $1-\dl \leq t \leq 2$. Then
$\phi := f(\psi)$ lies in~$\A$ by the $\Coo$-functional calculus.
Also, for all $x \in \supp(a)$, the estimates
$1 - \half\dl \leq \psi(x) \leq 1 + \half\dl$ hold, and so
$\phi(x) = f(\psi(x)) = 1$; this shows that $\phi a = a$.

The continuity of~$b$ shows that
$$
\supp(f(\psi)) \subseteq \overline{\set{x : \psi(x) > \dl}}
\subseteq \overline{\set{x : b(x) > \half\dl}}
\subseteq \set{x : b(x) \geq \half\dl} \subset U.
$$
Thus, $\supp\phi \subset U$, as required.
\end{proof}

\begin{prop}
\label{pr:smooth-quot}
Let $(\A,\H,\D)$ be a $QC^\infty$ spectral triple with $\A$
commutative and complete. Let $U \subset \spec(\A)$ be an open subset,
and let $h \in \A$ satisfy $h(x) \neq 0$ for all $x \in U$. Then
whenever $a \in \A$ with $\supp a \subset U$, $\A$ contains the
element $ah^{-1} \in C(\spec(\A))$ defined by
\begin{equation}
(ah^{-1})(x) := \begin{cases}
a(x)/h(x) & \text{if } h(x) \neq 0, \\
0 & \text{otherwise}. \end{cases}
\label{eq:ah-inv}
\end{equation}
\end{prop}

\begin{proof}
The formula \eqref{eq:ah-inv} clearly defines a continuous function on
$\spec(\A)$, so that $ah^{-1} \in A = C(\spec(\A))$.
To check that it lies in $\A$, it is enough to replace $h$ by an
invertible element $\tilde h \in \A$ for which
$\tilde h(x) = h(x)$ whenever $a(x) \neq 0$: since $\A = \A_\dl$ is a
pre-$C^*$-algebra, $\tilde h^{-1}$ will lie in $\A$, and thus
$ah^{-1} = a\tilde h^{-1} \in \A$.

Let $\eps := \inf\set{|h(x)| : x \in \supp a}$; since $\A$ is unital,
$\supp a$ is compact and therefore $\eps > 0$. Let
$V := U \cap \set{x : |h(x)| > \eps/2}$. By
Corollary~\ref{cr:partn-unity}, we can find $\phi \in \A$ so that
$0 \leq \phi \leq 1$, $\phi \equiv 0$ on~$\supp a$ and $\phi \equiv 1$
outside~$V$. The element $h + \half\eps \phi \in \A$ coincides with
$h$ on $\supp a$ and vanishes only on the compact set
$\set{x \notin V : h(x) = -\half\eps}$. If this set is nonvoid, we
can likewise find $\psi \in \A$, which is nonzero on this set and
vanishes on~$V$, so that $\tilde h := h + \half\eps \phi + \psi$
vanishes nowhere on $\spec(\A)$. This gives the required
$\tilde h \in \A$ such that $\tilde h \equiv h$ on $\supp a$.
\end{proof}

\begin{corl}
\label{cr:mat-inv}
Let $(\A,\H,\D)$ be a $QC^\infty$ spectral triple with $\A$
commutative and complete. Let $U \subset \spec(\A)$ be an open subset,
and let $h \in M_k(\A)$ be such that $h(x) \in M_k(\C)$ is invertible
for all $x \in U$. Let $a \in \A$ with $\supp a \subset U$; then the
element $ah^{-1} \in M_k(A)$ defined by
$$
(ah^{-1})(x) := \begin{cases}
(a(x) \ox 1_k)\,h(x)^{-1} & \text{if } h(x) \neq 0, \\
0 & \text{otherwise}, \end{cases}
$$
actually lies in the subalgebra $M_k(\A)$.
\end{corl}

\begin{proof}
The proof of Proposition~\ref{pr:smooth-quot} goes through with minor
modifications; for instance, one may take
$V := U \cap \set{x : |\det(h(x))| > \eps/2}$. By adding to~$h$
suitable scalar matrices which vanish on $\supp a$, one
constructs an invertible element $\tilde h \in M_k(\A)$ such that
$ah^{-1} = (a \ox 1_k)\tilde h^{-1}$, where $\tilde h^{-1} \in M_k(\A)$
since $M_k(\A)$ is a pre-$C^*$-algebra.
\end{proof}

\section{Geometric properties of noncommutative manifolds}
\label{sec:conds}

The conditions on a spectral triple that we introduce below will
control several interdependent features. Before introducing these
conditions, we first discuss several such features: summability,
metrics and differential structures.

We recall the symmetric operator ideals $\L^{p,\infty}(\H)$, for
$1 \leq p < \infty$; these are discussed in detail in
\cite[IV.2.$\al$]{Book} and \cite[7.C]{Polaris}; in
\cite{ConnesAction} and \cite{Polaris} they are called $\L^{p+}(\H)$.
The Dixmier ideal $\L^{1,\infty}(\H)$ is the common domain of the
Dixmier traces
$$
\Trw : \L^{1,\infty}(\H) \to \C.
$$
These are labelled by an uncountable index set of generalized limits
($\omega$-limits \cite{CareyPS}), but they are effectively computable
only on the subspace of ``measurable'' operators $T$ for which all
values $\Trw T$ coincide; for instance, trace-class operators satisfy
$\Trw T = 0$. Once we have established that a certain operator $T$ is
indeed measurable, we shall write $\ncint T$ instead of $\Trw T$ to
denote the common value of its Dixmier traces. If the limit
$$
\lim_{n\to\infty} \frac{1}{\log n} \sum_{k=0}^n \mu_k(T)
$$
exists, where the $\mu_k(T)$ are the singular values of $T$ in
nonincreasing order, then $T$ is measurable and this limit equals
$\ncint T$. A partial converse has been established by Lord, Sedaev
and Sukochev \cite{LordSS}: for \textit{positive} $T$, measurability
is \textit{equivalent} to the existence of this limit.

\begin{defn}
\label{df:p+-summ}
A spectral triple $(\A,\H,\D)$ is \textit{$p^+$-summable}, with
$1 \leq p < \infty$, if $(1 + \D^2)^{-1/2} \in \L^{p,\infty}(\H)$.
\end{defn}

For convenience, we abbreviate
$$
\Dreg := (1 + D^2)^{1/2},
$$
recalling that $\Dreg - |\D|$ is bounded, by functional calculus.

\begin{rmk}
If $\Dreg^{-1} \in \L^{p,\infty}(\H)$, it follows that
$A := \Dreg^{-p} \in \L^{1,\infty}(\H)$, and hence that
$\Trw \Dreg^{-p}$ is finite for any Dixmier trace~$\Trw$. Now
if $q > p$, then $\Dreg^{-q} = A^{q/p}$ lies in the ideal
$\L^1(\H)$ of trace-class operators, so $\Trw \Dreg^{-q} = 0$. It
follows that there is at most one value of~$p$ (independent
of~$\omega$) for which $\Trw \Dreg^{-p}$ can be both finite and
positive. Since
$$
\Trw \Dreg^{-p} = \omlim_{n\to\infty}
\frac{1}{\log n} \sum_{k=0}^n \mu_k(\Dreg^{-p})
$$
and \cite{CareyPS} for any bounded positive sequence $\{t_n\}$, one 
can estimate:
$$
\liminf_{n\to\infty} t_n \leq \omlim_{n\to\infty} t_n
\leq \limsup_{n\to\infty} t_n,
$$
then if $\liminf t_n > 0$, every $\omlim t_n$ is also positive. If
$0 < \Trw \Dreg^{-p} < \infty$ for all $\omega$, we shall call $p$ the
\textit{metric dimension} of the spectral triple $(\A,\H,\D)$.
\end{rmk}

If $\eta\: A \to A/\C\,1$ is the linear quotient map, then
$\|[\D,a]\|$ depends only on the image $\eta(a)$ of $a \in \A$.
Suppose that the set
\begin{equation}
\set{\eta(a) \in \A/\C\,1 : \|[\D,a]\| \leq 1}
\label{eq:bddset}
\end{equation}
is norm bounded in the Banach space $A/\C\,1$. Then the following
formula defines a bounded metric distance on the state space of~$A$,
as follows from \cite{ConnesMetric}:
\begin{equation}
d(\phi,\psi) := \sup\set{|\phi(a) - \psi(a)| : \|[\D,a]\| \leq 1}.
\label{eq:metric}
\end{equation}
(In Appendix~B, we show that any irreducible unital spectral triple
$(\A,\H,\D)$ determines a possibly unbounded distance function, which
actually suffices for the purposes of our proof.)

When $\A$ is commutative, the distance function \eqref{eq:metric} is
determined by its restriction to the subspace of pure states, which
may be identified with $X = \spec(A)$. In the case of $\A = \Coo(M)$
where $M$ is a compact spin$^c$ manifold, and $\D$ is a Dirac operator
arising from a Riemannian metric $g$ on~$M$, this $d$ coincides
\cite{ConnesMetric} with the Riemannian distance function determined
by~$g$.

\begin{defn}
If $(\A,\H,\D)$ is a commutative spectral triple for which the set
\eqref{eq:bddset} is bounded, we define the \textit{metric topology}
on the pure state space of $A$ to be the topology defined by the
distance function~\eqref{eq:metric}.
\end{defn}

\begin{rmk}
The equation \eqref{eq:metric} entails the inequality
\begin{equation}
|\phi(a) - \psi(a)| \leq \|[\D,a]\| \,d(\phi,\psi),
\label{eq:Lip-d}
\end{equation}
so that all $a\in\A$ are Lipschitz for the metric topology on~$X$. When
$[\D,a] \neq 0$ for $a \notin \C\,1$, the two topologies coincide
\cite{Pavlovic,RieffelMetr} if and only if the set \eqref{eq:bddset}
is precompact in $A/\C\,1$.
\end{rmk}

\begin{rmk}
\textit{A priori}, the metric topology may be finer than the original
weak$^*$ topology on~$X$. In particular, the metric topology need not
be compact unless the two topologies coincide. We shall henceforth
adopt the convention, when discussing continuous functions on~$X$ and
so forth, that \textit{the topology of~$X$ is by default its weak$^*$
topology}, unless the metric topology is explicitly invoked.
\end{rmk}

To exhibit the differential structure of a spectral triple, we first
recall the universal graded differential algebra $\Omega^\8\A$ of any
associative algebra $\A$ \cite{ConnesNCDiffG}. It is generated as an
algebra by symbols $a$, $da$ for $a \in \A$ subject to the preexisting
algebra relations of $\A$, the derivation rule
$d(ab) = a\,db + da\,b$, and the relations
$$
d(a_0 \,da_1 \dots da_k) = da_0\,da_1 \dots da_k, \qquad
d(da_1\,da_2 \dots da_k) = 0.
$$
We may then identify $\Omega^\8\A$ with the (normalized) Hochschild
complex of $\A$, that is,
$$
\Omega^k\A \simeq C^k(\A) := \A \ox (\A/\C\,1)^{\ox k}, \qquad
a_0\,da_1 \dots da_k \otto a_0 \ox \eta(a_1) \oxyox \eta(a_k).
$$
When $\A$ is a Fr\'echet algebra, one generally uses the projective
topological tensor product to topologize $\Omega^\8\A$.

\begin{defn}
If $(\A,\H,\D)$ is a spectral triple, we shall use the notation $\CDA$
for the subalgebra of $\B(\H)$ generated by $\A$ and $\d\A = [\D,\A]$.
We can define an (algebra) representation
$\pi_\D \: \Omega^\8\A \to \CDA$ by setting
$$
\pi_\D(a_0\,da_1 \dots da_k) := a_0\,[\D,a_1]\dots[\D,a_k].
$$
We may regard $\Omega^\8\A$ as an involutive algebra by setting
$(da)^* := - d(a^*)$; then $\pi_\D $ is a $*$-repre\-sentation of the
Hochschild chains of $\A$ as operators on~$\H$.
\end{defn}

However, this $\CDA$ is not a graded algebra (although the count of
$[\D,a]$ factors does give a filtration), and $\pi_\D$ is
not a representation of graded differential algebras. As is well known
from physical examples \cite{Book,Cordelia,SchueckerZ}, there may be
nontrivial ``junk forms'' $\omega \in \Omega^\8\A$ such that
$\pi_\D(\omega) = 0$ but $\pi_\D(d\omega) \neq 0$. On quotienting out
the junk, we obtain a graded differential algebra \cite[VI.1]{Book}:
$$
\La^\8_\D \A := \CDA/J_\pi,  \word{where}
J_\pi = \pi_\D(d(\ker \pi_\D)).
$$
In particular, $J_\pi$ is a differential ideal. The subspaces
$\La^k_\D \A =
\pi_\D(\Omega^k\A)/\pi_\D(d(\Omega^{k-1}\A \cap \ker \pi_\D))$ give
the grading; the differential $\d$ on $\La^\8_\D \A$ is defined on
equivalence classes of operators, modulo junk terms, by
$$
\d\bigl( a_0\,[\D,a_1]\dots[\D,a_k] + \text{junk}\bigr)
:= [\D,a_0]\,[\D,a_1]\dots[\D,a_k] + \text{junk}.
$$
Here $\La^0_\D \A = \A$, and $\La^1_\D \A$ may be identified with the
$\A$-bimodule of operators of the form $\sum_j a_j\,[\D,b_j]$ (finite
sum), with $a_j,b_j \in \A$. Nontrivial junk terms appear in higher
degrees.

\subsection{Axiomatic conditions on commutative spectral triples}
\label{ssc:geom-cond}

{}From now on, let $(\A,\H,\D)$ be a spectral triple whose algebra
$\A$ is \textit{commutative} (and unital). We shall also assume that
the $C^*$-algebra $A = \ol{\A}$ is \textit{separable}.

In \cite{ConnesReal}, Connes introduces several conditions on such an
$(\A,\H,\D)$ in order to specify axiomatically what a noncommutative
spin geometry should be. We now list these conditions (for the
commutative case), as well as a few supplementary requirements which
we need to establish our main results.

\begin{cond}[\textit{Dimension}] 
\label{cn:metr-dim}
The spectral triple $(\A,\H,\D)$ is $p^+$-summable for a fixed
\textit{positive integer} $p$, for which $\Trw \Dreg^{-p} > 0$
for all~$\omega$.

By the remark after Definition~\ref{df:p+-summ}, this condition
determines $p$ uniquely; we then say that the critical summability
parameter $p$ is the ``metric dimension'' of $(\A,\H,\D)$.
\end{cond}

\begin{rmk}
A priori, there is no reason why the growth of the eigenvalues of~$\D$
should be such that $p$ is an integer. However, the orientability
condition below introduces another dimensionality parameter~$p$ as the
degree of a certain Hochschild cycle, which is necessarily an integer,
and we require that these two quantities coincide.

This formulation excludes certain interesting ``$0$-dimensional''
cases, such as occur when $\A$ has finite (linear) dimension. For
noncommutative $0$-dimensional spectral triples built over matrix
algebras, we refer to \cite{IochumKM,Krajewski,PaschkeS}.
\end{rmk}

\begin{cond}[\textit{Metric}] 
\label{cn:metric}
The set $\set{\eta(a) \in \A/\C\,1 : \|[\D,a]\| \leq 1}$ is
norm-bounded in the Banach space $A/\C\,1$. This ensures that the
character space $X = \spec(A)$ is a metric space \cite{ConnesMetric}
with the metric distance \eqref{eq:metric}. (See Appendix~B).
\end{cond}

\begin{rmk}
\label{rk:metrizable}
Since we have assumed that $A$ is separable, the space $X$ with its
weak$^*$ topology is metrizable. However, the metric on $X$ defined by
the equation \eqref{eq:metric} does not necessarily give the weak$^*$
topology.
\end{rmk}

\begin{cond}[\textit{Regularity}] 
\label{cn:qc-infty}
The spectral triple $(\A,\H,\D)$ is $QC^\infty$, as set forth in
Definition~\ref{df:qc-infty}. Without loss of generality, we assume
that $\A$ is complete in the topology given by~\eqref{eq:dl-snorms}
and so is a \textit{Fr\'echet pre-$C^*$-algebra}.
\end{cond}

\begin{cond}[\textit{Finiteness}] 
\label{cn:finite}
The dense subspace of $\H$ which is the smooth domain of~$\D$,
$$
\H_\infty := \bigcap_{m\geq 1} \Dom \D^m
$$
is a \textit{finitely generated projective $\A$-module}. Moreover, 
there exists a Hermitian pairing
$$
\pairing{\cdot}{\cdot} \: \H_\infty \x \H_\infty \to \A,
$$
making $\bigl( \H_\infty, \pairing{\cdot}{\cdot} \bigr)$ a full
pre-$C^*$ right $\A$-module; and such that, for some particular
Dixmier trace $\Tr_\Omega$, the following relation holds:
\begin{equation}
\Tr_\Omega\bigl( \pairing{\xi}{\eta}\,\Dreg^{-p} \bigr)
= \braket{\xi}{\eta},
\word{for all}  \xi,\eta \in \H_\infty.
\label{eq:herm-pairing}
\end{equation}
Here $\braket{\.}{\.}$ denotes the scalar product on~$\H$.
\end{cond}

\begin{rmk}
Since $\A$ is commutative, we are free to regard $\H_\infty$ as
either a right or left $\A$-module. As a dense subspace of~$\H$, it
is naturally a left module via the representation~$\pi$, but it is
algebraically more convenient to treat it as a right module; thus
$\H_\infty \simeq q\,\A^m$ where $q \in M_m(\A)$ is 
(a selfadjoint) idempotent.
\end{rmk}

\begin{rmk}
It is proved in Appendix~A that (up to positive scalar multiples) the
only Hermitian form which can satisfy the conditions listed here is
the standard one, expressible as
$\pairing{\xi}{\eta} = \sum_{j,k} \xi_j^* q_{jk} \eta_k$ on
identifying $\H_\infty$ with~$q\A^m$ with $q$~selfadjoint.
\end{rmk}

\begin{cond}[\textit{Absolute Continuity}] 
\label{cn:abs-cont}
For all nonzero $a \in \A$ with $a \geq 0$, and for any
$\omega$-limit, the following Dixmier trace is positive:
$$
\Trw(a\Dreg^{-p}) > 0.
$$
\end{cond}

\begin{rmk}
The absolute continuity condition actually subsumes the dimension
condition, and they should properly be regarded as one and the same.
This formulation eases the adaptation of the conditions for nonunital
algebras~\cite{Himalia,RennieSumm}.
\end{rmk}

\begin{cond}[\textit{First Order}] 
\label{cn:first-ord}
The bounded operators in $[\D,\A]$ commute with $\A$; in other words,
\begin{equation}
[[\D, a], b] = 0  \word{for all}  a,b \in \A.
\label{eq:first-ord}
\end{equation}
This condition says that operators in $[\D,\A]$ can be regarded as
\textit{endomorphisms} of the $\A$-module $\H_\infty$; and more
generally, that $\CDA \subseteq \End_\A(\H_\infty)$.
\end{cond}

\begin{rmk}
For spectral triples over noncommutative algebras, the first-order
condition is more elaborate: as well as the representation
$\pi\: \A \to \B(\H)$ we require a commuting representation
$\pi^\op \: \A^\op \to \B(\H)$ of the opposite algebra $\A^\op$ (or
equivalently, an antirepresentation of~$\A$ that commutes with~$\pi$):
writing $a$ for $\pi(a)$ as usual, and $b^\op$ for $\pi^\op(b)$, we
ask that $[a, b^\op] = 0$. Now $\H_\infty$ can be regarded as a right
$\A$-module under the action $\xi \. b := \pi^\op(b)\,\xi$. The
first-order condition is then expressed as: $[[\D,a], b^\op] = 0$ for
$a,b \in \A$; and once again it entails that
$\CDA \subseteq \End_\A(\H_\infty)$.

Coming back to the commutative case, we take $\pi^\op = \pi$ from now
on. (But see \cite{Krajewski}, for instance, for examples of
commutative algebras with different left and right actions on~$\H$.)
\end{rmk}

\begin{cond}[\textit{Orientability}] 
\label{cn:orient}
Let $p$ be the metric dimension of $(\A,\H,\D)$. We require that the
spectral triple be \textit{even}, with $\Z_2$-grading $\Ga$, if and
only if $p$ is even. For convenience, we take $\Ga = 1$ when $p$ is
odd. We say the spectral triple $(\A,\H,\D)$ is \textit{orientable} if
there exists a Hochschild $p$-cycle
\begin{subequations}
\label{eq:Hoch-cycle}
\begin{equation}
\cc = \sum_{\al=1}^n a_\al^0 \ox a_\al^1 \oxyox a_\al^p \in Z_p(\A,\A)
\label{eq:Hoch-cycle-defn}
\end{equation}
whose Hochschild class $[\cc] \in HH_p(\A)$ may be called the
``orientation'' of $(\A,\H,\D)$, such that
\begin{equation}
\label{eq:Hoch-cycle-repn}
\pi_D(\cc) \equiv \sum_\al a_\al^0 \,[\D,a_\al^1] \dots [\D,a_\al^p]
= \Ga.
\end{equation}
\end{subequations}
\end{cond}

\begin{cond}[\textit{Poincar\'e duality}] 
\label{cn:pdual}
The spectral triple $(\A,\H,\D)$ determines a $K$-homology class for
$A \ox A$. Let
$\mu = [(\A \ox \A,\H,\D)] \in K^\8(\A \ox \A) = K^\8(A \ox A)$ denote
this $K$-homology class. We require that $\mu$ be a fundamental class,
i.e., that the Kasparov product~\cite{HigsonR}
$$
\8 \ox_A \mu : K_\8(A) \to K^\8(A)
$$
be an isomorphism. (More on this in Section~\ref{sec:PD}).
\end{cond}

\begin{cond}[\textit{Reality}] 
\label{cn:real}
There is an antiunitary operator $J \: \H \to \H$ such that
$J a^* J^{-1} = a$ for all $a \in \A$; and moreover, $J^2 = \pm 1$,
$J \D J^{-1} = \pm\D$ and also $J \Ga J^{-1} = \pm\Ga$ in the even
case, according to the following table of signs depending only on
$p \bmod 8$:
$$
\begin{array}[t]{|c|cccc|}
\hline
p \bmod 8            & 0 & 2 & 4 & 6 \rule[-5pt]{0pt}{17pt} \\
\hline
J^2 = \pm 1          & + & - & - & + \rule[-5pt]{0pt}{17pt} \\
J\D J^{-1} = \pm\D   & + & + & + & + \rule[-5pt]{0pt}{17pt} \\
J\Ga J^{-1} = \pm\Ga & + & - & + & - \rule[-5pt]{0pt}{17pt} \\
\hline
\end{array}
\qquad\qquad
\begin{array}[t]{|c|cccc|}
\hline
p \bmod 8          & 1 & 3 & 5 & 7 \rule[-5pt]{0pt}{17pt} \\
\hline
J^2 = \pm 1        & + & - & - & + \rule[-5pt]{0pt}{17pt} \\
J\D J^{-1} = \pm\D & - & + & - & + \rule[-5pt]{0pt}{17pt} \\
\hline
\end{array}
\belowdisplayskip=1pc
$$
For the origin of this sign table in $KR$-homology, we refer to
\cite[Sec.~9.5]{Polaris}.
\end{cond}

\begin{rmk}
For a noncommutative algebra $\A$, we would require
$J a^* J^{-1} = a^\op$ or, more precisely,
$J \pi(a)^* J^{-1} = \pi^\op(a)$. Thus, $J$ implements on~$\H$ the
involution $\tau : a \ox b^\op \mapsto b^* \ox {a^*}^\op$ of
$\A \ox \A^\op$. 
\end{rmk}

\begin{cond}[\textit{Irreducibility}] 
\label{cn:irred}
The spectral triple $(\A,\H,\D)$ is \textit{irreducible}: that is, the
only operators in $\B(\H)$ (strongly) commuting with $\D$ and with all
$a \in \A$ are the scalars in $\C\,1$.
\end{cond}

The foregoing conditions are an elaboration of those set forth in
\cite{ConnesGrav} for the reconstruction of a spin manifold. We add a
final condition, which provides a cohomological version of Poincar\'e
duality: see the discussion in \cite[VI.4.$\ga$]{Book}.

\begin{cond}[\textit{Closedness}] 
\label{cn:closed}
The $p^+$-summable spectral triple $(\A,\H,\D)$ satisfies the
following closedness condition: for any $\row{a}{1}{p} \in \A$, the
operator $\Ga\,[\D,a_1]\dots [\D,a_p]\Dreg^{-p}$ has vanishing
Dixmier trace; thus, for any $\Trw$,
\begin{equation}
\Trw\bigl(\Ga\,[\D,a_1]\dots[\D,a_p]\,\Dreg^{-p}\bigr) = 0.
\label{eq:closed-chain}
\end{equation}
\end{cond}

\begin{rmk}
By setting $\Phi(\row{a}{0}{p}) :=
\Trw\bigl(\Ga\,a_0\,[\D,a_1]\dots[\D,a_p]\,\Dreg^{-p}\bigr)$,
the equation \eqref{eq:closed-chain} may be rewritten \cite[VI.2]{Book}
as $B_0\Phi = 0$, where $B_0$ is defined on $(k+1)$-linear functionals
by $(B_0\phi)(\row{a}{1}{k}) :=
\phi(1,\row{a}{1}{k}) + (-1)^k \phi(\row{a}{1}{k},1)$.
\end{rmk}

We quote Lemma 3 of \cite[VI.4.$\ga$]{Book}, adapted to the present
case where $\A$ is commutative and $(\A,\H,\D)$ is $p^+$-summable.

\begin{lemma}[Connes]
\label{lm:PDin-cohom}
Let $(\A,\H,\D)$ be $p^+$-summable and satisfy
Condition~\ref{cn:first-ord} (first order). Then for each
$k = 0,1,\dots,p$ and $\eta \in \Omega^k\A$, a Hochschild cocycle
$C_\eta \in Z^{p-k}(\A,\A^*)$ is defined by
$$
C_\eta(a^0,\dots,a^{p-k}) := \Trw\bigl( \Ga\,\pi_\D(\eta)\,
a^0\,[\D,a^1]\dots [\D,a^{p-k}]\,\Dreg^{-p} \bigr).
$$
Moreover, if Condition~\ref{cn:closed} (closedness) also holds, then
$C_\eta$ depends only on the class of $\pi_\D(\eta)$ in $\La_\D^k \A$,
and $B_0 C_\eta = (-1)^k \,C_{d\eta}$.
\qed
\end{lemma}

\subsection{First consequences of the geometric conditions}
\label{ssc:first-lemmas}

We now describe some immediate consequences of these conditions, which
already give a reasonable picture of the spaces and bundles we shall
employ. In this subsection, $(\A,\H,\D)$ will always denote a
$QC^\infty$ spectral triple whose algebra $\A$ is \textit{commutative}
and \textit{complete} (and unital, too). In other words,
Condition~\ref{cn:qc-infty} (regularity) is taken for granted. We
shall write, as before, $X = \spec(\A) = \spec(A)$ where $A$ is the
separable $C^*$-completion of~$\A$; it is a metrizable compact
Hausdorff space under its weak$^*$ topology.

\begin{lemma}
\label{lm:no-proj}
Under Conditions \ref{cn:first-ord} and \ref{cn:irred} (first order,
irreducibility), the algebra $\A$ contains no nontrivial projector.
\end{lemma}

\begin{proof}
Let $q \in \A$ be a projector. Then
$$
[\D,q] = [\D,q^2] = q\,[\D,q] + [\D,q]\,q = 2q\,[\D,q],
$$
where the first order condition gives the last equality. Hence
$(2q - 1)[\D,q] = 0$, implying $[\D,q] = 0$ since $2q - 1$ is
invertible. Thus, $q$ commutes with $\D$, and with all $a \in \A$
since $\A$ is commutative: by irreducibility, $q$ must be a scalar,
either $q = 0$ or $q = 1$.
\end{proof}

\begin{corl}
\label{cr:spec-conn}
Under the same Conditions \ref{cn:first-ord} and \ref{cn:irred}, the
space $X$ is connected.
\end{corl}

\begin{proof}
It is enough to show that there are no nontrivial projectors in the
$C^*$-algebra $A$. Let $\tilde q \in A$ be a projector; by
Lemma~\ref{lm:dense-proj}, we can find a projector $q \in \A$ such
that $\|q - \tilde q\| < \half$. Since $q$ is either $0$ or~$1$ by
Lemma~\ref{lm:no-proj}, the same must be true of $\tilde q$.
Therefore, $C(X)$ contains no nontrivial projectors, and so $X$ is
connected.
\end{proof}

\begin{lemma}
\label{lm:gotta-bundle}
Under Condition~\ref{cn:finite} (finiteness), the dense subspace
$\H_\infty \subset \H$ consists of continuous sections of a complex
vector bundle $S \to X$. If $(\A,\H,\D)$ is also irreducible, then $S$
has constant rank.
\end{lemma}

\begin{proof}
Condition~\ref{cn:finite} says that there is an integer $m > 0$ and a
projector $q \in M_m(\A)$ such that $\H_\infty \simeq q \A^m$ as a
(right) $\A$-module. We may regard $q$ as an element of $M_m(A)$; the
$A$-module $qA^m$ is well-defined as a finitely generated projective
module over $A = C(X)$. By the Serre--Swan theorem \cite{SwanVect},
$qA^m \simeq \Ga(X,S)$, the space of continuous sections of a complex
vector bundle $S \to X$. From the finiteness axiom, the Hermitian 
pairing on~$q\A^m$ gives $S \to X$ the structure of a Hermitian 
vector bundle.

If $(\A,\H,\D)$ is irreducible, then $X$ is connected by
Corollary~\ref{cr:spec-conn}, and the rank of $S$ must be constant.
We shall denote this rank by~$N$.
\end{proof}

\begin{lemma}
\label{lm:oper-endo}
Under Conditions \ref{cn:finite} to \ref{cn:orient} (finiteness,
absolute continuity, first order, orientability), the algebra $\CDA$
is a unital selfadjoint subalgebra of $\Ga(X,\End S)$. The operator
norm of each $T \in \CDA$ coincides with its norm as an endomorphism
of~$S$.
\end{lemma}

\begin{proof}
Since $(\A,\H,\D)$ is $QC^\infty$, each operator $T \in \CDA$ maps
$\H_\infty$ into itself; and the first order condition ensures that
$T$ is an $\A$-linear map on $\H_\infty$. If
$T = \sum_j a_j\,[\D,b_j]\in \La^1_\D\A = \pi_\D(\Omega^1\A)$,
the adjoint operator
$$
T^* = -\tsum_j [\D,b_j^*]\,a_j^*
= \tsum_j b_j^*\,[\D,a_j^*] - [\D,b_j^*a_j^*]
$$
lies in $\La^1_\D\A$ also, so $\La^1_\D\A$ is a selfadjoint linear
subspace of $\B(\H)$. Thus, the algebra $\CDA$ generated by $\A$ and
$\La^1_\D\A$ is a $*$-subalgebra of $\B(\H)$. Moreover, since the
pairing on $\H_\infty$ is determined by the scalar product on~$\H$ via
\eqref{eq:herm-pairing}, we conclude that
$\pairing{\xi}{T\eta} = \pairing{T^*\xi}{\eta}$ for each $T \in \CDA$.
Consequently, $T$ yields an adjointable $A$-module map of the
$C^*$-module $\Ga(X,S)$; that is,
$\CDA \subset \End_A(\Ga(X,S)) = \Ga(X, \End S)$. The algebra $\CDA$ 
contains $\Ga^2 = 1$.

We use the inequality \cite[Cor.~2.22]{RaeburnW} between positive
elements of the $C^*$-algebra $A$:
$$
\pairing{T\xi}{T\xi} \leq \|T\|^2_{\End S} \pairing{\xi}{\xi},
$$
where $\|T\|_{\End S}$ denotes the norm of $T$ in the $C^*$-algebra
$\Ga(X, \End S)$. Therefore, when $\xi \in \H_\infty$,
\begin{align*}
\braket{T\xi}{T\xi}
&= \Tr_\Omega\bigl(\pairing{T\xi}{T\xi}\,\Dreg^{-p}\bigr)
\\
&\leq \|T\|^2_{\End S}
\Tr_\Omega\bigl(\pairing{\xi}{\xi}\,\Dreg^{-p}\bigr)
= \|T\|^2_{\End S} \, \braket{\xi}{\xi}.
\end{align*}
Majorizing this inequality over
$\set{\xi \in \H_\infty : \braket{\xi}{\xi} \leq 1}$, we obtain
$\|T\| \leq \|T\|_{\End S}$.

To see that these norms are indeed equal, suppose that
$0 \leq M < \|T\|^2_{\End S}$, so that $M - T^*T$ is selfadjoint
and not positive in $\Ga(X, \End S)$. Then we can find a nonzero
$\xi \in \Ga(X,S)$ such that
$\pairing{T\xi}{T\xi} - M\,\pairing{\xi}{\xi}$ is positive and
nonzero. In view of Condition~\ref{cn:abs-cont}, this implies that
\begin{align*}
M\,\braket{\xi}{\xi}
&= \Tr_\Omega\bigl(M\,\pairing{\xi}{\xi}\,\Dreg^{-p}\bigr)
\\
&< \Tr_\Omega\bigl(\pairing{T\xi}{T\xi}\,\Dreg^{-p}\bigr)
\\
&= \braket{T\xi}{T\xi} = \|T\xi\|^2 \leq \|T\|^2 \,\braket{\xi}{\xi},
\end{align*}
so that $M < \|T\|^2$ since $\xi \neq 0$. This is true for all
$M < \|T\|^2_{\End S}$, thus $\|T\| = \|T\|_{\End S}$.
\end{proof}

\begin{corl}
\label{cr:full-rank}
Under the same Conditions \ref{cn:finite} to~\ref{cn:orient}, the
algebra of sections $\CDA$ is pointwise a direct sum of matrix
algebras:
$$
(\CDA)_x \simeq M_{k_1}(\C) \oplus\ldots\oplus M_{k_r}(\C)
\word{for}  x \in X,
$$
where $k_1 +\cdots+ k_r = N$.
\end{corl}

\begin{proof}
Lemma~\ref{lm:oper-endo} shows that $(\CDA)_x$ is a selfadjoint
subalgebra of the finite-dimensional algebra $\End S_x$; hence it is a
direct sum of full matrix algebras. This subalgebra has full rank,
since $\Ga = \pi_\D(\cc)$ lies in $\CDA$, so that $\Ga_x^2$ is the
identity element in $\End S_x$.
\end{proof}

The following result allows us to sidestep several questions of 
domains.

\begin{prop}
\label{pr:lin-bdd}
Let $T\: \H_\infty \to \H_\infty$ be $\A$-linear. Then $T$ extends to
a bounded operator on~$\H$.
\end{prop}

\begin{proof}
Let $\xi_1,\dots,\xi_m \in \H_\infty \simeq q\A^m$ be defined by
$\xi_j := q \eps_j = \sum_k q_{kj} \eps_k$ where $\eps_j \in \A^m$ is
the column-vector with $1$ in the $j$th slot and zeroes elsewhere.
Every $\xi \in \H_\infty$ can be written in the form
$\xi = \sum_{j=1}^m \xi_j a_j$, for some $a_j \in \A$. By the
properties of our chosen Hermitian pairing, we get
$$
\pairing{\xi_j}{\xi_k}_{\H_\infty}
= \bigpairing{\tsum_r q_{rj} \eps_r}{\tsum_s q_{sk} \eps_s}_{\A^m}
= \sum_{r,s} q_{jr} q_{sk} \dl_{rs} = q_{jk}.
$$
Thus, as already noted,
$\pairing{\xi}{\xi} = \sum_{j,k} a^*_j q_{jk} a_k$.

The $\A$-linearity of $T$ gives $T\xi = \sum_{j=1}^m (T\xi_j) a_j$,
and by hypothesis, $T\xi_j \in \H_\infty$. Therefore,
\begin{align*}
\pairing{T\xi_j}{T\xi_k}
&= \bigpairing{T \tsum_n q_{nj} \eps_n}{T \tsum_l q_{lk} \eps_l}
= \bigpairing{T \tsum_{n,r} q_{nr} q_{rj} \eps_n}
			 {T \tsum_{l,s} q_{ls} q_{sk} \eps_l}
\\
&= \bigpairing{T \tsum_{n,r} q_{nr} \eps_n q_{rj}}
			  {T \tsum_{l,s} q_{ls} \eps_l q_{sk}}
= \sum_{r,s} q_{jr} \pairing{T\xi_r}{T\xi_s} q_{sk}.
\end{align*}
Denote by $\Theta_{\xi,\eta}$, for $\xi,\eta \in \H_\infty$, the
``ketbra'' operator $\rho \mapsto \xi \pairing{\eta}{\rho}$. The
pairing $\pairing{T\xi}{T\xi}$ may be expanded as follows:
\begin{align*}
\pairing{T\xi}{T\xi}
&= \sum_{j,k} a_j^* \pairing{T\xi_j}{T\xi_k} a_k
= \sum_{j,k,r,s} a_j^* \pairing{\xi_j}{\xi_r} \pairing{T\xi_r}{T\xi_s}
\pairing{\xi_s}{\xi_k} a_k
\\
&= \sum_{r,s} \pairing{\xi}{\xi_r} \pairing{T\xi_r}{T\xi_s}
\pairing{\xi_s}{\xi}
= \sum_{r,s}
\bigpairing{T\xi_r \pairing{\xi_r}{\xi}}{T\xi_s \pairing{\xi_s}{\xi}}
\\
&= \sum_{r,s}
\bigpairing{\Theta_{T\xi_r,\xi_r} \xi}{\Theta_{T\xi_s,\xi_s} \xi}
= \sum_{r,s}
\bigpairing{\Theta_{\xi_s \pairing{T\xi_s}{T\xi_r},\xi_r} \xi}{\xi}
\\
&\leq \sum_{r,s}
\bigl\| \Theta_{\xi_s \pairing{T\xi_s}{T\xi_r},\xi_r} \bigr\|
\pairing{\xi}{\xi}.
\end{align*}
In the last line here the norm is both the operator norm and the
endomorphism norm, which coincide by Lemma~\ref{lm:oper-endo}. The
norm of each ketbra is finite, since each $T\xi_r \in \H_\infty$ by
hypothesis. Now we can estimate the operator norm of $T$; for
$\xi \in \H_\infty$ we get the bound
\begin{align*}
\braket{T\xi}{T\xi}
&= \Tr_\Omega\bigl( \pairing{T\xi}{T\xi}\,\Dreg^{-p} \bigr)
\\
&\leq \sum_{r,s}
\bigl\| \Theta_{\xi_s \pairing{T\xi_s}{T\xi_r},\xi_r} \bigr\| \,
\Tr_\Omega\bigl( \pairing{\xi}{\xi} \,\Dreg^{-p} \bigr)
\\
&= \sum_{r,s}
\bigl\| \Theta_{\xi_s \pairing{T\xi_s}{T\xi_r},\xi_r} \bigr\| \,
\braket{\xi}{\xi}.
\end{align*}
Since the $\xi_r$ are a fixed finite set of vectors, the operator norm
of $T$ is finite, with
$$
\|T\|^2 \leq \sum_{r,s} \bigl\| \pairing{T\xi_s}{T\xi_r} \bigr\|
\, \bigl\| \pairing{\xi_r}{\xi_r} \bigr\|^{1/2}
\, \bigl\| \pairing{\xi_s}{\xi_s} \bigr\|^{1/2}.
$$
This last expression for the norm follows from
\cite[Lemma~2.30]{RaeburnW} or \cite[Lemma~4.21]{Polaris}.
\end{proof}

\begin{lemma}
\label{lm:bundle-meas}
Under Conditions \ref{cn:metr-dim}, \ref{cn:finite}, \ref{cn:abs-cont}
and~\ref{cn:orient} ($p^+$-summability, finiteness, absolute
continuity, orientability), the $\A$-valued Hermitian pairing on
$\H_\infty$ given by \eqref{eq:herm-pairing} is independent of the
choice of Dixmier trace.
\end{lemma}

\begin{proof}
Connes' character theorem \cite[Thm.~IV.2.8]{Book} ---we refer to
\cite{Polaris} and \cite{CareyPRSone} for its detailed proof---
shows that any operator of the form
\begin{equation}
T = \Ga \tsum_\al a^0_\al\,[\D,a^1_\al]\dots [\D,a^p_\al]\,
\Dreg^{-p},
\label{eq:mble-oper}
\end{equation}
where $\cc = \sum_\al a^0_\al \ox a^1_\al \oxyox a^p_\al$ is a
Hochschild \textit{cycle}, is a measurable operator.
Condition~\ref{cn:orient} provides us with such a Hochschild $p$-cycle
$\cc$ for which $\pi(\cc) = \Ga$. Using $\Ga^2 = 1$, we can rewrite
\eqref{eq:herm-pairing} as
\begin{equation}
\braket{\xi}{\eta}
= \Tr_\Omega\bigl(\pairing{\xi}{\eta}\,\Dreg^{-p}\bigr)
= \Tr_\Omega\bigl(\Ga\pairing{\xi}{\eta}\Ga\,\Dreg^{-p}\bigr).
\label{eq:Hoch-mble}
\end{equation}
If $a = \pairing{\xi}{\eta}$, then $a\cc$ is also a Hochschild cycle
for $\A$ ---as an easy consequence of the cycle property of~$\cc$ and
the commutativity of~$\A$--- so the right hand side
of~\eqref{eq:Hoch-mble} is $\Tr_\Omega(T)$, where
$T = \Ga\,\pi_D(a\cc)\,\Dreg^{-p}$ is indeed of the form
\eqref{eq:mble-oper}. Thus $\Tr_\Omega$ may be replaced, in the
formula~\eqref{eq:herm-pairing}, by any other Dixmier trace~$\Trw$.
\end{proof}

It was noted in \cite{ConnesSpec} that the orientability condition
yields the following expression for $\D$ in terms of commutators
$\d a = [\D,a]$ and $[\D^2,a]$.

\begin{lemma}
\label{lm:dee-formula}
Under Condition~\ref{cn:orient} (orientability), the operator $\D$
verifies the following formula (as an operator on $\H_\infty$):
\begin{equation}
\D = \half (-1)^{p-1} \Ga \sum_{\al=1}^n \sum_{j=1}^p (-1)^{j-1}
a^0_\al \,\d a^1_\al \dots \d a^{j-1}_\al \,[\D^2,a^j_\al]
\,\d a^{j+1}_\al \dots \d a^p_\al + \half (-1)^{p-1} \Ga \,\d\Ga,
\label{eq:D-formula}
\end{equation}
where $\Ga = \sum_{\al=1}^n a^0_\al \,\d a^1_\al \dots \d a^p_\al$
and we write
$\d\Ga := \sum_{\al=1}^n \d a^0_\al\,\d a^1_\al \dots \d a^n_\al$.
\end{lemma}

\begin{proof}
First note that on the domain $\H_\infty$, the derivation $\ad D^2$
may be written as
\begin{equation}
[\D^2,a] = \D\,\d a + \d a\,\D  \word{for all}  a \in \A.
\label{eq:der}
\end{equation}
Thus, the summation over~$j$ in \eqref{eq:D-formula} telescopes, to
give
\begin{align*}
\sum_{\al=1}^n
& \sum_{j=1}^p (-1)^{j-1} a^0_\al \,\d a^1_\al \dots
\d a^{j-1}_\al \,[\D^2,a^j_\al] \,\d a^{j+1}_\al \dots \d a^p_\al
\\
&= \sum_{\al=1}^n (a^0_\al\,\D\,\d a^1_\al \dots \d a^p_\al
+ (-1)^{p-1} a^0_\al \,\d a^1_\al \dots \d a^p_\al \,\D)
\\
&= - \sum_{\al=1}^n \d a^0_\al \,\d a^1_\al \dots \d a^p_\al
+ \D\Ga + (-1)^{p-1} \Ga \D
\\
&= - \d\Ga + 2\,(-1)^{p-1} \Ga \D,
\end{align*}
and \eqref{eq:D-formula} follows on multiplying both sides by
$\half (-1)^{p-1} \Ga$.
\end{proof}

\begin{corl}
\label{cr:comm-formula}
Under Conditions \ref{cn:first-ord} and~\ref{cn:orient} (first order,
orientability), the commutator $[\D,a]$, for $a \in \A$, has the
expansion
$$
[\D,a] = \half(-1)^{p-1} \Ga \sum_{\al=1}^n \sum_{j=1}^p (-1)^{j-1}
a^0_\al \,\d a^1_\al \dots (\d a^j_\al\,\d a + \d a\,\d a^j_\al)
\dots \d a^p_\al.
$$
\end{corl}

\begin{proof}
The first order condition entails that $a$ commutes with all operator
factors in the expansion~\eqref{eq:D-formula}, except the
$[\D^2,a^j_\al]$ factors. For those, \eqref{eq:der} and
$[[\D,a^j_\al], a] = 0$ imply
$$
[[\D^2,a^j_\al], a] = [\D[\D,a^j_\al], a] + [[\D,a^j_\al]\D, a]
= [\D,a]\,[\D,a^j_\al] + [\D,a^j_\al]\,[\D,a].
\eqno \qed
$$
\hideqed
\end{proof}

\section{The cotangent bundle}
\label{sec:cotg-bdl}

Throughout this section, $(\A,\H,\D)$ will be a spectral triple whose
algebra $\A$ is (unital and) commutative and complete; and
$X = \spec(\A)$ will be its metrizable compact Hausdorff character
space. Moreover, we shall assume that Conditions~\ref{cn:metr-dim},
\ref{cn:qc-infty}--\ref{cn:orient} and~\ref{cn:closed} hold, namely
that $(\A,\H,\D)$ is $p^+$-summable, $QC^\infty$ and has the
properties of finiteness, absolute continuity, first order,
orientability and closedness.

\begin{lemma}
\label{lm:metr-junk}
The operator $[\D,a]\,[\D,b] + [\D,b]\,[\D,a]$ is a junk term, for any
$a,b \in \A$.
\end{lemma}

\begin{proof}
We must show that $da\,db + db\,da$ belongs to $d(\ker\pi_\D)$ in
the universal graded differential algebra $\Omega^\8\A$. Since
$da\,db + db\,da = d(a\,db - d(ba) + b\,da) = d(a\,db - db\,a)$,
it is enough to notice that the first-order condition gives
$$
\pi_\D(a\,db - db\,a) = a\,[\D,b] - [\D,b]\,a = 0.
\eqno \qed
$$
\hideqed
\end{proof}

\begin{lemma}
\label{lm:quot-nzero}
The image of $\Ga = \pi_\D(\cc)$ in $\La_\D^p \A$ is nonzero.
\end{lemma}

\begin{proof}
The Hochschild cycle $\cc \in Z_p(\A,\A)$ defines a Hochschild
$0$-cocycle (a trace) $C_\cc$ on $\A$ by Lemma~\ref{lm:PDin-cohom}.
Taking into account Lemma~\ref{lm:bundle-meas}, it is given by
$$
C_\cc(a) = \ncint \Ga\pi_\D(\cc)\,a\Dreg^{-p} = \ncint a\Dreg^{-p}.
$$
Since
\begin{equation}
C_\cc(1) = \ncint \Dreg^{-p} > 0,
\label{eq:dim-isp}
\end{equation}
this cocycle does not vanish. Moreover, Condition~\ref{cn:closed}
entails that $C_\cc$ depends only on the class of $\pi_\D(\cc)$ in
$\La_\D^p(\A)$. If this class were zero, so that
$\pi_\D(\cc) \in \pi_\D(d(\ker\pi_\D))$, then we could write it as a
finite sum of the form
$\pi_\D(\cc) = \tsum_\bt \d b_\bt^1 \dots \d b_\bt^p$. But the
closedness condition would then apply to show that
$$
C_\cc(1) = \ncint \Ga \pi_\D(\cc) \Dreg^{-p}
= \tsum_\bt \ncint \Ga\,\d b_\bt^1\,\d b_\bt^2
\dots \d b_\bt^p \Dreg^{-p} = 0,
$$
contradicting \eqref{eq:dim-isp}. Hence, the class of $\Ga$ has a
nonzero image in $\La_\D^p(\A)$.
\end{proof}

\begin{corl}
\label{cr:skew-symm}
Let $\Ga' \in \CDA$ be defined by
\begin{equation}
\Ga' := \frac{1}{p!} \sum_{\sg\in S_p} (-1)^\sg \sum_\al a^0_\al
\,\d a_\al^{\sg(1)} \,\d a_\al^{\sg(2)} \dots \d a_\al^{\sg(p)},
\label{eq:skew-Gamma}
\end{equation}
on skewsymmetrizing the expression for $\Ga$ obtained
from~\eqref{eq:Hoch-cycle}. If $a \in \A$ is positive and nonzero,
then $a\Ga' \neq 0$.
\end{corl}

\begin{proof}
Let $a \in \A$ be positive, $a \neq 0$. Since $\A$ is commutative and
$\cc \in Z_p(\A,\A)$, the product $a\cc$ is also a Hochschild
$p$-cycle. Now the absolute continuity condition implies that
$$
C_{a\cc}(1) = C_\cc(a) = \ncint a\Dreg^{-p} > 0.
$$
Since $\pi_\D(a\cc) = a\Ga$, the proof of Lemma~\ref{lm:quot-nzero}
shows that the class $[a\Ga]$ in $\La_\D^p \A$ is nonzero. Now,
Lemma~\ref{lm:metr-junk} shows that $[a\Ga'] = [a\Ga]$. In particular,
$a\Ga' \neq 0$ as an element of $\B(\H)$.
\end{proof}

Thus, the skewsymmetrization $\Ga'$ of $\Ga$ given by
\eqref{eq:skew-Gamma} is nonzero as an operator on~$\H$, and \textit{a
fortiori} as a section in $\Ga(X, \End S)$. In fact, this section
vanishes nowhere on~$X$, as the proof of the following Proposition
shows.

In what follows, whenever $T$ is a continuous (local) section of
$\End S \to X$, we write either $T(x)$ or~$T|_x$ to denote its value
in the fibre $\End S_x$. The \textit{support} of $T$ will mean its
support as a section, namely, the complement of the largest open
subset $V \subseteq X$ such that $T(x) = 0$ in $\End S_x$ for all
$x \in V$.

\begin{prop}
\label{pr:lin-indp}
There is an open cover $\{\row{U}{1}{n}\}$ of $X$ such that, for each
$\al = 1,\dots,n$, the operators $[\D,a_\al^1], \dots, [\D,a_\al^p]$
are pointwise linearly independent sections of $\Ga(U_\al, \End S)$.
\end{prop}

\begin{proof}
Let $Z := \set{x \in X : \Ga'(x) = 0}$ be the zero set of~$\Ga'$.
If $V$ were a nonvoid open subset of~$Z$, then, using
Lemma~\ref{lm:partn-unity}, we could find a nonzero positive
$b \in \A$ such that $\supp b \subset V$; but this would imply
$b\Ga' = 0$, contradicting Corollary~\ref{cr:skew-symm}. Therefore,
$Z$ has empty interior.

The pairing on $\H_\infty$, or rather, on the completed $A$-module
$\Ga(X,S)$, induces a $C^*$-norm on $\Ga(X,\End S) = \End_A(\Ga(X,S))$
which in turn determines a norm on each fibre $\End S_x$ so that
$\|T\|_{\End S} = \sup_{x\in X} \|T(x)\|_{\End S_x}$ for
$T \in \Ga(X,\End S)$. Choose $\eps > 0$; unless $Z = \emptyset$,
there is an open set $W \supset Z$ such that
$\sup_{y\in W} \|\Ga'(y)\|_{\End S_y} < \eps$. Next choose $a \in \A$,
positive and nonzero, with $\supp a \subset W$. By
Lemma~\ref{lm:local-unit}, there exists $\psi \in \A$ such that
$0 \leq \psi \leq 1$, $\psi a = a$ and $\supp\psi \subset W$. Hence,
by Lemma~\ref{lm:oper-endo},
$\|\psi\Ga'\| = \|\psi\Ga'\|_{\End S} < \eps$.

By Corollary \ref{cr:skew-symm}, the Hochschild $0$-cocycles
$C_{a\Ga}$ and $C_{a\Ga'}$ are equal, and they define positive
functionals on $\A$. For any Dixmier trace $\Trw$,
we know that
\begin{align}
C_{a\Ga}(1) &= C_\Ga(a) = \Trw(a\Dreg^{-p}),
\label{eq:AG-cocyc}
\\
C_{a\Ga'}(1) &= \Trw(\Ga\Ga'\,a\Dreg^{-p})
= \Trw(\Ga\Ga'\,\psi a\Dreg^{-p}).
\notag
\end{align}
This yields the following estimate:
\begin{align*}
|C_{a\Ga'}(1)| = C_{a\Ga'}(1)
&= \Trw(\Ga\Ga'\psi a\Dreg^{-p})
\\
&\leq \|\Ga\Ga'\,\psi\| \Trw( a\Dreg^{-p} )
\\
&< \eps \Trw(a\Dreg^{-p}).
\end{align*}
Since $\Trw(a\Dreg^{-p}) > 0$, this is inconsistent with
\eqref{eq:AG-cocyc} when $0 < \eps < 1$. We conclude that the set $Z$
must necessarily be empty, so that $\Ga'(x)$ is nonvanishing on~$X$.

Thus, at each $x \in X$, there is some $\al$ such that
$a_\al^0(x) \neq 0$ and $\d a^1_\al(x), \dots, \d a^p_\al(x)$ in
$\End S_x$ have a nonzero skewsymmetrized product, and therefore are
linearly independent. We may now define $U_\al$ to be the open set of
all~$x$ for which this linear independence holds; and
$U_1 \cup\cdots\cup U_n = X$ from the nonvanishing of~$\Ga'$.
\end{proof}

\begin{lemma}
\label{lm:wedge}
Fix $\al \in \{1,\dots,n\}$ and let $a \in \A$, writing
$a =: a^{p+1}_\al$ for notational convenience. Then
\begin{equation}
\frac{(-1)^{p-1}}{2(p+1)!} \sum_{\sg\in S_{p+1}} (-1)^\sg a^0_\al
[\D,a^{\sg(1)}_\al] \dots [\D,a^{\sg(p)}_\al]\,[\D,a^{\sg(p+1)}_\al]
= 0.
\label{eq:Gamma-wedge}
\end{equation}
\end{lemma}

\begin{proof}
By Corollary~\ref{cr:comm-formula}, we may write
$$
\Ga\,[\D,a] = \half (-1)^{p-1} \sum_{\al=1}^n \sum_{j=1}^p
(-1)^{j-1} a^0_\al \,\d a^1_\al
\dots (\d a^j_\al\,\d a^{p+1}_\al + \d a^{p+1}_\al\,\d a^j_\al)
\dots \d a^p_\al.
$$
For each $\al$, every term in the sum over $j$ contains a symmetric
product of one-forms, so its skewsymmetrization vanishes.
\end{proof}

\begin{rmk}
For brevity, we shall denote by $\Ga'_\al$ the $\al$th summand of
$\Ga'$ in~\eqref{eq:skew-Gamma}, and by $\Ga'_\al \w \d a$ the
operator on the left hand side of~\eqref{eq:Gamma-wedge}. Now, $\d a$
and each $\d a_\al^j$, and therefore each $\Ga'_\al$, is an
endomorphism of~$\H_\infty$; Lemma~\ref{lm:wedge} shows that
$\Ga'_\al(x) \w \d a(x) = 0$ in~$\End S_x$, for each $x \in U_\al$,
where the notation $\w$ now denotes skewsymmetrization with respect to
the several $\d a_\al^j(x)$.
\end{rmk}

\begin{rmk}
{}From now on we shall assume, without any loss of generality, that
each $a^j_\al$, for $j = 1,\dots,p$, $\al = 1,\dots,n$, is
\textit{selfadjoint}. (Otherwise, we just take selfadjoint and
skewadjoint parts, allowing some repetition of the sets $U_\al$.)
Consequently, \textit{each $[\D,a^j_\al]$ is skewadjoint}.
\end{rmk}

\begin{prop}
\label{pr:covec-span}
The operators $[\D,a^j_\al]$, for $\al = 1,\dots,n$ and
$j = 1,\dots,p$, generate $\La^1_\D\A$ as a finitely generated
projective $\A$-module.
\end{prop}

\begin{proof}
Let $a \in \A$; choose (and fix) $\al$ such that $a\Ga'_\al \neq 0$.
Then by Lemma~\ref{lm:wedge}, $\Ga'_\al \w \d a = 0$ in $\CDA$, and
thus $\Ga'_\al(x) \w \d a(x) = 0$ in $\End S_x$, for each
$x \in U_\al$. Let $E_x$ be the (complex) vector subspace of
$\End S_x$ spanned by the endomorphisms
$\d a^1_\al(x),\dots,\d a^p_\al(x)$. The exterior algebra $\La^\8 E_x$
is represented on $S_x$ by
$$
(v^1 \wyw v^k) \. \xi := \frac{1}{k!} \sum_{\sg\in S_k} (-1)^\sg
v^{\sg(1)} \dots v^{\sg(k)} \,\xi,
$$
and this representation is faithful, on account of
Proposition~\ref{pr:lin-indp}. Similarly, we can represent on~$S_x$
the exterior algebra of the vector subspace
$$
E'_x := \linspan \set{\d a^1_\al(x), \dots, \d a^p_\al(x), \d a(x)}
\supseteq E_x.
$$
Now Lemma~\ref{lm:wedge} implies that $\La^\8 E'_x = \La^\8 E_x$, and
thus $E'_x = E_x$; therefore, $\d a(x)$ lies in~$E_x$, for all
$x \in U_\al$.

Choose a partition of unity $\{\phi_\al\}_{\al=1}^n \subset \A$
subordinate to the open cover $\{U_\al\}$, as in
Lemma~\ref{lm:partn-unity}. Then for any $a \in \A$ we may write
$[\D,a] = \sum_\al \phi_\al\,[\D,a]$. Then for each $x \in U_\al$, the
linear independence of the $\d a^j_\al(x)$ yields unique
constants $c_{j\al}(x)$ such that
\begin{equation}
\phi_\al(x)\,\d a(x) = \sum_{j=1}^p c_{j\al}(x) \,\d a^j_\al(x).
\label{eq:single-sum}
\end{equation}
Since $\supp \phi_\al \subset U_\al$, \eqref{eq:single-sum} defines a
continuous local section in $\Ga(U_\al,\End S)$. By uniqueness,
$c_{j\al}(x) = 0$ outside $\supp\phi_\al$, so we may extend $c_{j\al}$
by zero to a function on all of~$X$, and thus we may regard these
local sections as elements of $\Ga(X,\End S)$. We claim that each
$c_{j\al}$ lies in~$\A$ (and in particular is continuous).

Define an $\A$-valued Hermitian pairing on $\La^1_\D\A$ by setting
\begin{equation}
\pairing{\d a}{\d b} := C_p \tr((\d a)^*\,\d b),
\label{eq:prng-defn}
\end{equation}
where $C_p$ is a suitable positive normalization constant, and $\tr$
denotes the matrix trace in $\End_\A \H_\infty = q M_m(\A) q$, where
$\Ga(X,S) = q A^m$. To see that this pairing takes values in $\A$, we
use the following localization argument.

Choose a finite open cover $\{V_\rho\}$ of $X$ such that $S$ is
trivial over each $V_\rho$. For each $\rho$, choose $N = \rank S$
elements $\xi^\rho_1,\dots,\xi_N^\rho \in \H_\infty$ which, when
regarded as sections of $S$, are linearly independent over~$V_\rho$.
Moreover, these sections can be chosen so that
$\{\xi^\rho_1(x),\dots,\xi_N^\rho(x)\}$ is an orthonormal basis in
each fibre $S_x$, for $x \in V_\rho$; this means that for all
$b \in \A$ with $\supp b \subset V_\rho$ the orthogonality relations
$b\pairing{\xi^\rho_i}{\xi^\rho_j} = b\delta_{ij}$ hold. That may be
achieved by Gram--Schmidt orthogonalization, on invoking
Proposition~\ref{pr:smooth-quot} to see that local inverses of
elements in $\A$ also lie in $\A$. Next, choose a partition of unity
$\{\psi_\rho\}$ subordinate to the cover $\{V_\rho\}$, with
$\psi_\rho \in \A$, as in Lemma~\ref{lm:partn-unity}. Then
$$
\tr((\d a)^*\,\d b) = \sum_\rho \psi_\rho\,\tr((\d a)^*\,\d b)
= \sum_\rho \sum_{j=1}^N
\psi_\rho\,\pairing{\d b\,\xi_j^\rho}{\d a\,\xi_j^\rho}.
$$
The right hand side is a finite sum of elements of $\A$, and so 
belongs to~$\A$.

If $b = \pairing{\d a}{\d a}$, then $b(x)$ is the trace of a positive
element of $\End S_x$, so $b(x) = 0$ if and only if $\d a(x) = 0$;
thus the pairing is positive definite. Consider the matrix
$g_\al = [g_\al^{jk}] \in M_p(\A)$ given by
\begin{equation}
g_\al^{jk} := \pairing{\d a_\al^j}{\d a_\al^k}
= - C_p \tr(\d a_\al^j \,\d a_\al^k).
\label{eq:metr-defn}
\end{equation}
The matrix $g_\al(x) \in M_p(\C)$ has the form
$C_p \tr(m_\al(x)^* m_\al(x))$ where $m_\al(x) \in (\End S_x)^p$ is
the $p$-column with linearly independent entries $\d a_\al^j(x)$.
Thus, for $x \in U_\al$, each $g_\al(x)$ is a positive definite Gram
matrix, hence invertible, when $x \in U_\al$. Let
$g_\al^{-1}(x) := [g_{\al,ij}](x)$ denote the inverse matrix.

We may now invoke Corollary~\ref{cr:mat-inv} ---recall that $\A$ is
complete--- to conclude that $\phi_\al g_{\al,ij}$ is an element
of~$\A$ for $i,j = 1,\dots,p$. Now if
$\phi_\al\,\d a = \sum_i c_{i\al}\,\d a_\al^i$ with
$\supp c_{i\al} \subset U_\al$ as in~\eqref{eq:single-sum}, we find
that
\begin{align*}
c_{j\al} &= \sum_{i,k} c_{i\al} \,g_{\al,jk} \,g_\al^{ki}
= - C_p \sum_{i,k} g_{\al,jk} \tr(\d a_\al^k\, c_{i\al}\,\d a_\al^i)
\\
&= - C_p \sum_k g_{\al,jk} \tr(\d a_\al^k\,\phi_\al\,\d a)
= \sum_k \phi_\al\,g_{\al,jk} \,
\pairing{\d a_\al^k}{\d a},
\end{align*}
where each $\phi_\al\,g_{\al,jk} \in \A$ and
$\pairing{\d a_\al^k}{\d a} \in \A$ by previous arguments; we conclude
that $c_{j\al} \in \A$.

Finally, for any $a \in \A$, we may now write
\begin{equation}
\d a = \sum_{\al=1}^n \phi_\al\,\d a
= \sum_{\al=1}^n \sum_{j=1}^p c_{j\al} \,\d a^j_\al \in \La^1_\D\A.
\label{eq:ctg-gen-global}
\end{equation}
Since the coefficients $c_{j\al}$ in this finite sum lie in~$\A$, the
$\d a^j_\al = [\D,a^j_\al]$ generate the $\A$-module $\La^1_\D\A$.

To see that $\La^1_\D\A$ is a \textit{projective} $\A$-module, we
rewrite the coefficients in \eqref{eq:ctg-gen-global} as
$c_{j\al} = \sum_{k=1}^p \psi_{jk\al} \,\pairing{\d a_\al^k}{\d a}$,
and we get, for $b \in \A$,
$$
b\,\d a = \sum_{\al=1}^n \sum_{j,k=1}^p
\psi_{jk\al} \,\pairing{\d a_\al^k}{b\,\d a} \,\d a_\al^j,
$$
and therefore $\La^1_\D\A \simeq Q\,\A^{np}$ via standard isomorphisms
\cite[Prop.~3.9]{Polaris}, where $Q \in M_{np}(\A)$ is the projector
with entries $Q_{\al k,\bt l}
:= \sum_{m=1}^p \psi_{lm\bt} \,\pairing{\d a_\al^k}{\d a_\bt^l}$.
\end{proof}

We shall frequently need to replace the expansion
\eqref{eq:ctg-gen-global} by a ``localized'' version for a
single~$\al$, as follows.

\begin{corl}
\label{cr:local-sum}
If $a \in \A$ is such that $\supp \d a \subset U_\al$, then there
exist $c_{1\al},\dots,c_{p\al} \in \A$, compactly supported
in~$U_\al$, such that
\begin{equation}
\d a = \sum_{j=1}^p c_{j\al} \,\d a^j_\al.
\label{eq:ctg-gen-local}
\end{equation}
More generally, if $b \in \A$, then for any open $V \subset U_\al$
there are continuous functions $b_{j\al} \: V \to \C$ for
$j = 1,\dots,p$, such that
\begin{equation}
\d b(x) = \sum_{j=1}^p b_{j\al}(x) \,\d a^j_\al(x)
\quad \text{for all } x \in V,
\label{eq:db-local}
\end{equation}
and such that each $c\,b_{j\al} \in \A$ whenever $c \in \A$ with
$\supp c \subset V$.
\end{corl}

\begin{proof}
If $\supp \d a \subset U_\al$, then we may choose the partition of
unity of the previous proof such that $\phi_\al(x) = 1$ on
$\supp \d a$, by Corollary~\ref{cr:partn-unity}. Thus
$\phi_\al\,\d a = \d a$ and $\phi_\bt\,\d a = 0$ for
$\bt \neq \al$. Thus both \eqref{eq:single-sum} and
\eqref{eq:ctg-gen-global} reduce to \eqref{eq:ctg-gen-local}. By
construction, $\supp c_{j\al} \subseteq \supp \phi_\al$.

In the same way, if $c \in \A$ with $\supp c \subset V$, we may
expand $c\,\d b =: \sum_{j=1}^p c'_{j\al} \,\d a^j_\al$ with
$c'_{j\al} \in \A$ and $\supp c'_{j\al} \subseteq \supp c$.
Uniqueness of the coefficients at each $x \in V$ shows that
$c'_{j\al}(x) = c(x) b_{j\al}(x)$, where each function $b_{j\al}$
does not depend on~$c$; also, $b_{j\al}$ is continuous because its
restriction to each compact subset of~$V$ is continuous.
\end{proof}

With the local linear independence and spanning provided by
Propositions \ref{pr:lin-indp} and~\ref{pr:covec-span}, we now obtain
a (complex) vector subbundle $E$ of $\End S$, such that
$\La^1_\D\A \subseteq \Ga(X,E)$. This vector bundle will eventually
play the role of the complexified cotangent bundle $T^*_\C(X)$,
although at this stage we have not yet identified a suitable
differential structure on~$X$.

\begin{prop}
\label{pr:gotta-bundle}
For each $x \in U_\al$, define a $p$-dimensional complex vector space
by
\begin{equation}
E_x := \linspan \set{\d a^1_\al(x),\dots, \d a^p_\al(x)}
\subseteq \End S_x.
\label{eq:cotg-fibre}
\end{equation}
Then these spaces form the fibres of a complex vector bundle $E \to X$.
\end{prop}

\begin{proof}
We prove that $E$ is a vector bundle by providing transition
functions satisfying the usual \v{C}ech cocycle condition.

For each pair of indices $\al,\bt$, Corollary \ref{cr:local-sum}
provides continuous functions
$c^k_{j\al\bt} \: U_\al \cap U_\bt \to \C$ such that
\begin{equation}
\d a^k_\al(x) = \sum_{j=1}^p  c^k_{j\al\bt}(x) \,\d a^j_\bt(x)
\word{for all}  x \in U_\al \cap U_\bt.
\label{eq:dak-combo}
\end{equation}
Whenever $x \in U_\al \cap U_\bt \cap U_\ga$, this entails the
additional relation
$$
\d a^k_\al(x) = \sum_{j=1}^p c^k_{j\al\bt}(x) \,\d a^j_\bt(x)
= \sum_{j,l=1}^p c^k_{j\al\bt}(x) \,c^j_{l\bt\ga}(x) \,\d a^l_\ga(x),
$$
and the linear independence of the $\d  a^l_\ga(x)$ shows that
$$
c^k_{l\al\ga}(x) = \sum_{j=1}^p c^k_{j\al\bt}(x) \, c^j_{l\bt\ga}(x)
\word{for all}  x \in U_\al \cap U_\bt \cap U_\ga.
$$
In particular, the matrix $c_{\al\bt}(x) = [c^k_{j\al\bt}(x)]$ is
invertible with $c_{\al\bt}^{-1}(x) = c_{\bt\al}(x)$ for
$x \in U_\al \cap U_\bt$. The relation \eqref{eq:dak-combo} and its
analogue with $\al$ and $\bt$ exchanged show that the vector space
$E_x$ of~\eqref{eq:cotg-fibre} is well defined, independently
of~$\al$.

Moreover, the cocycle conditions
$c_{\al\bt} c_{\bt\ga} = c_{\al\ga}$ hold over every nonvoid
$U_\al \cap U_\bt \cap U_\ga$, so these are continuous transition
matrices for a vector bundle $E \to X$, whose total space is the
disjoint union $E := \biguplus_{x\in X} E_x$.
\end{proof}

\begin{corl}
\label{cr:real-bundle}
If for each $x \in U_\al \subset X$, we define the real vector space
$$
E_{\R,x} = \R\mbox{-}\linspan \set{\d a^1_\al(x),\dots,\d a^p_\al(x)},
$$
then $E_\R := \biguplus_{x\in X} E_{\R,x}$ is the total space of a
real vector bundle over~$X$.
\end{corl}

\begin{proof}
We need only show that the transition functions are actually real
matrices. Since each $\d a^j_\bt$ is skewadjoint, taking the adjoint
(in $\End S_x$) of \eqref{eq:dak-combo} yields
$$
\d a^k_\al(x) = \sum_{j=1}^p \bar{c}^k_{j\al\bt}(x) \,\d a^j_\bt(x).
$$
By uniqueness of the coefficients, we conclude that
$\bar{c}^k_{j\al\bt} = c^k_{j\al\bt}$ for each $j,k,\al,\bt$.
\end{proof}

We conclude this Section by indicating that the functions $a_\al^j$
are not constant on sets with nonempty interior, and more importantly,
that the operator $\D$ is actually \emph{local}.

\begin{lemma}
\label{lm:nunca-denso}
Let $Y \subset U_\al$ be a level set of the function $a^j_\al$, for
some $j = 1,\dots,p$. Then $Y$ is closed and its interior $\intr Y$
is empty.
\end{lemma}

\begin{proof}
Clearly $Y$ is closed, since $a^j_\al \in C(X)$. Suppose that
$\intr Y$ were nonempty; then there would be a nonzero element
$f \in \A$ such that $\supp f \subset \intr Y$. Let $\la \in \R$ be
the value of $a^j_\al$ on the level set~$Y$, so that
$\la f = a^j_\al f$. Taking commutators with $\D$ gives
$$
\la\,[\D,f] = a^j_\al\,[\D,f] + [\D, a^j_\al]\,f
= a^j_\al\,[\D,f] + f\,[\D, a^j_\al],
$$
using the first order condition. Therefore, $f(y)\,\d a^j_\al(y) = 0$
for all $y \in Y$. This contradicts $\supp f \subset U_\al$, since
by definition $U_\al \subseteq \set{x \in X : \d a^j_\al(x) \neq 0}$
for each~$j$.
\end{proof}

\begin{corl}
\label{cr:nunca-denso}
For each $\al = 1,\dots,n$, let $a_\al \: U_\al \to \R^p$ be the
mapping with components $a^j_\al$, $j = 1,\dots,p$. Then any
level set of the mapping $a_\al$ is a closed set with empty interior.
\end{corl}

\begin{proof}
Any such level set for $a_\al$ is the intersection of level sets for
the several $a^j_\al$.
\end{proof}

\begin{corl}
\label{cr:local-opr}
If $a \in \A$, then $\supp(\d a) \subseteq \supp a$.
\end{corl}

\begin{proof}
Suppose that $Y := (X \less \supp a)$ is nonempty, otherwise there is
nothing to prove. Then $Y$ is a nonvoid open subset of~$X$, and the
function $a$ vanishes on its closure $\ol Y$. Arguing as in the proof
of Lemma~\ref{lm:nunca-denso}, with $\la = 0$ and $a^j_\al$ replaced
by~$a$, we see that $\d a(y) = 0$ for all $y \in Y$.
\end{proof}

\begin{lemma}
\label{lm:nearly-open}
If $V \subseteq U_\al$ is open, then $a_\al(V)$ has nonempty interior
in $\R^p$.
\end{lemma}

\begin{proof}
The level sets of each $a^j_\al$ on $U_\al$ are closed with no
interior, so no $a^j_\al$ is constant on $V$, or on any open subset
of~$V$. Then the range of $a^1_\al$ contains a nontrivial interval
(i.e., not a point), and so it contains an open interval
$(s,t) \subset \R$. Let $V_1 = V \cap (a^1_\al)^{-1}((s,t))$, an open
subset of~$V$. Likewise, since $a^2_\al$ is not constant on $V_1$, we
can find an open $V_2 \subseteq V_1$ which $a^j_\al$ maps onto an open
subinterval of $a^j_\al(V_1)$ for $j = 1,2$; and so on. After
$p$~steps, we obtain an open subset $V_p \subseteq V$ that $a_\al$
maps onto an open rectangle in~$\R^p$.
\end{proof}

\begin{corl}
\label{cr:open-elbows}
Let $x \in U_\al$ be such that $x$ is neither a local maximum nor
minimum of any of the functions $a^j_\al$, $j = 1,\dots,p$. Then there
is an open neighbourhood of~$x$ on which $a_\al$ is an open map.
\end{corl}

\begin{proof}
The hypothesis says that $a_\al(x)$ is not an endpoint of any interval
in $a^j_\al(U_\al)$ for any $j$. Thus we can find an open rectangle
$a_\al(x) \subset R \subseteq a_\al(U_\al)$, whereby $a_\al^{-1}(R)$
is an open neighbourhood of~$x$, such that every point
$y \in a^{-1}_\al(R)$ also satisfies the hypothesis of the corollary.
\end{proof}

\begin{corl}
\label{cr:nearly-open}
If $B \subseteq a_\al(U_\al) \subseteq \R^p$ has empty interior, then
$a_\al^{-1}(B)\cap U_\al$ has empty interior also.
\end{corl}

\begin{proof}
If $a_\al^{-1}(B) \cap U_\al$ had an interior point, then so too would
$a_\al(a_\al^{-1}(B) \cap U_\al) = B$.
\end{proof}

\section{A Lipschitz functional calculus}
\label{sec:lip-open}

\begin{defn}
Let $(\A,\H,\D)$ again be a spectral triple whose algebra $\A$ is
(unital and) commutative and complete. From now on, we shall say that
$(\A,\H,\D)$ is a \textit{spectral manifold of dimension~$p$} if it is
$QC^\infty$, $p^+$-summable, and satisfies the metric, first order,
finiteness, absolute continuity, orientability, irreducibility and
closedness conditions, that is, all postulates of
Subsection~\ref{ssc:geom-cond} except perhaps Conditions
\ref{cn:pdual} and~\ref{cn:real}.
\end{defn}

Throughout this section, we shall assume that $(\A,\H,\D)$ is a
spectral manifold of dimension~$p$. We shall use without comment the
open cover $\{U_\al\}$ of $X = \spec(\A)$ provided by
Proposition~\ref{pr:lin-indp}, and the vector bundle $E \to X$
afforded by Proposition~\ref{pr:gotta-bundle}. As indicated earlier,
we shall also assume that the $a_\al^j$ appearing
in~\eqref{eq:Hoch-cycle}, for $j \neq 0$, are selfadjoint.

Our next task is to develop a Lipschitz version of the functional
calculus. For each $\al = 1,\dots,n$, we shall denote the joint
spectrum of $a^1_\al,\dots,a^p_\al$ by $\Delta_\al$, and shall write
$a_\al := (a_\al^1,\dots,a^p_\al)$ as a continuous mapping from
$U_\al$ to~$\R^p$.

We recall Nachbin's extension \cite{Nachbin} of the Stone--Weierstrass
approximation theorem to subalgebras of differentiable functions.

\begin{prop}[Nachbin \cite{Nachbin}]
\label{pr:Nachbin}
Let $U$ be an open subset of~$\R^p$, and let $\B \subset C^r(U,\R)$
for $r \in \{1,2,\dots\}$. A necessary and sufficient condition for
the algebra generated by $\B$ to be dense in $C^r(U,\R)$, in the $C^r$
topology, is that the following conditions hold:
\begin{enumerate}
\item
For each $x \in U$, there exists $f \in \B$ such that $f(x) \neq 0$.
\item
Whenever $x,y \in U$ with $x \neq y$, there exists $f \in \B$ such
that $f(x) \neq f(y)$.
\item
For each $x \in U$ and each nonzero tangent vector $\xi_x \in T_x U$,
there exists $f \in \B$ such that $\xi_x(f) \neq 0$.
\end{enumerate}
\end{prop}

In particular, the real polynomials on $\R^p$, restricted to $U$, are
$C^r$-dense in $C^r(U,\R)$; and thus also, the complex-valued
polynomials are $C^r$-dense in $C^r(U) = C^r(U,\C)$.

\begin{lemma}
\label{lm:Lip-est}
Let $(\A,\H,\D)$ be a spectral manifold of dimension~$p$, and let
$Y$ be a compact subset of~$U_\al$ with nonempty interior. Let
$a \in \A$ with $\supp a \subset Y$, and suppose there is a
bounded $C^1$ function $f\: L \to \C$ such that
$a = f(a^1_\al,\dots,a^p_\al)$, where $L$ is a bounded open subset 
of~$\R^p$ with $\Delta_\al \subset L$. Then there are positive
constants $C'_\al$ and $C_Y$, independent of $a$ and~$f$, such that
\begin{equation}
C'_\al\,\|[\D,a]\| \leq \|df\|_\infty \leq C_Y\,\|[\D,a]\|
\label{eq:Lip-est}
\end{equation}
where
$\|df\|_\infty^2 := \sup_{t\in\Delta_\al} \sum_j |\del_j f(t)|^2$.
\end{lemma}

\begin{proof}
By Proposition~\ref{pr:Nachbin}, we may approximate $f$ by polynomials
$p_k$ such that $p_k \to f$ and $\del_j p_k \to \del_j f$ for
each~$j$, uniformly on $\Delta_\al$. The first order condition shows
that \eqref{eq:Dcomm-poly} holds for each~$p_k$. In the limit, the 
proof of Proposition~\ref{pr:mult-cinfty} yields
\begin{equation}
[\D,a] = [\D, f(a^1_\al,\dots,a^p_\al)]
= \sum_{j=1}^p \del_j f(a^1_\al,\dots,a^p_\al) \,[\D,a^j_\al].
\label{eq:local-expan}
\end{equation}

Since $\supp a \subset Y$, Corollary~\ref{cr:partn-unity} provides an
element $\phi \in \A$ such that $0 \leq \phi \leq 1$ and
$\supp\phi \subset U_\al$ while $\phi(x) = 1$ for $x \in Y$; and
consequently, $\phi a = a$. Moreover, $\phi[\D,a] = [\D,a]$, on
account of Corollary~\ref{cr:local-opr}. The elements
$g_\al^{jk} \in \A$ defined by \eqref{eq:metr-defn} form a positive
definite matrix of functions on~$U_\al$; again let $g_{\al,ij}$ denote
the entries of its inverse matrix. The proof of
Proposition~\ref{pr:covec-span} shows that $\phi\,g_{\al,ij}$ lies
in~$\A$, for all~$i,j$.

For $x \in U_\al$, we compute that
\begin{align*}
\tr\biggl( \sum_j g_{\al,ij} \,\d a^j_\al \,\d a \biggr)(x)
&= \sum_{j,k=1}^p \del_k f(a^1_\al,\dots,a^p_\al)(x) g_{\al,ij}(x)
\tr\bigl( \d a^j_\al(x) \,\d a^k_\al(x) \bigr)
\\
&= -C_p^{-1} \sum_{j,k=1}^p \del_k f(a^1_\al,\dots,a^p_\al)(x)
g_{\al,ij}(x) g_\al^{jk}(x)
\\
&= -C_p^{-1}\, \del_i f(a^1_\al,\dots,a^p_\al)(x).
\end{align*}

The trace is defined pointwise in the endomorphism algebra
$\Ga(U_\al, \End S)$. Let $\|\.\|_x$ be the operator norm in
$\End S_x$ induced by the hermitian pairing \eqref{eq:prng-defn}
on this algebra. Then
\begin{align}
|(\del_i f)(a^1_\al,\dots,a^p_\al)(x)|
&= C_p \biggl| \tr\biggl( \sum_j g_{\al,ij} \,\d a^j_\al \,\d a
\biggr)(x) \biggr|
\leq \biggl| \sum_j
\bigpairing{g_{\al,ij}\,\d a^j_\al}{\phi\,\d a}(x) \biggr|
\nonumber \\
&\leq \sum_j \bigpairing{\phi\,g_{\al,ij}\,\d a^j_\al}
{\phi\,g_{\al,ij}\,\d a^j_\al}^{1/2}(x)
\bigpairing{\d a}{\d a}^{1/2}(x)
\nonumber \\
&\leq \sum_j p\, 
\|\phi(x)g_{\al,ij}(x)\,\d a^j_\al(x)\|_x \,\|\d a(x)\|_x
\nonumber \\
&\leq \sum_j p\, \|g_{\al,ij}(x)\,\d a^j_\al(x)\|_x
\,\|\d a(x)\|_x \leq B_i(x)\, \|\d a(x)\|_x \,,
\label{eq:Lip-est2}
\end{align}
where $B_i(x)$ is independent of $a$ and $f$ and is bounded on~$Y$.
Since $\|\d a(x)\|_x \leq \|[\D,a]\|$ by Lemma~\ref{lm:oper-endo},
taking $C_Y := \max_i \sup\set{B_i(x) : x \in Y}$ yields the second
inequality in~\eqref{eq:Lip-est}.

On the other hand, for any $x \in U_\al$, the estimate
\begin{align*}
\|\d a(x)\|_x
&\leq \sum_j |(\del_j f)(a^1_\al,\dots,a^p_\al)(x)| \,
\|\d a^j_\al(x)\|_x
\\
&\leq \sum_j \|\d a^j_\al\|_{\End S}
\,|(\del_j f)(a^1_\al,\dots,a^p_\al)(x)|
\\
&\leq \sum_j \|[\D,a^j_\al]\| \,\|df\|_\infty,
\end{align*}
again by Lemma~\ref{lm:oper-endo}, shows that
$\|[\D,a]\| \leq (C'_\al)^{-1} \|df\|_\infty$ if we take
$(C'_\al)^{-1} := \sum_j \bigl\|[\D,a^j_\al]\bigr\|$; this gives the
first inequality in \eqref{eq:Lip-est}.
\end{proof}

We need to remove the hypothesis that $a\in\A$ has compact support in
$Y$. This is possible because the proof of Lemma \ref{lm:Lip-est}
relies on pointwise estimates.

\begin{corl}
\label{cr:local-lipest}
Let $a \in \A$ and let $Y \subset U_\al$ be any compact subset such
that $a|_Y = (f \circ a_\al)|_Y$ for some $C^1$ function $f$ defined
and bounded in a neighbourhood of $a_\al(Y) \subset \R^p$. Then there
are positive constants $C'_\al$ and $C_Y$, independent of $a$ and~$f$,
such that
\begin{equation}
C'_\al\,\|[\D,a]\|_Y \leq \|df\|_{Y,\infty} \leq C_Y\,\|[\D,a]\|
\label{eq:local-lipest}
\end{equation}
where
$\|df\|_{Y,\infty}^2 := \sup_{t\in a_\al(Y)} \sum_j |\del_j f(t)|^2$,
and
$\|[\D,a]\|_Y := \sup_{x\in Y} \|[\D,a](x)\|_x$.
\end{corl}

\begin{proof}
Choose $\phi \in \A$ with $\phi(x) = 1$ for all $x \in Y$ and
$\supp\phi \subset U_\al$. Then for all $x \in Y$, the proof of
Lemma~\ref{lm:Lip-est} gives us
$$
|(\del_i f \circ a_\al)(x)|
\leq \sum_j p\, \|g_{\al,ij}(x)\,\d a^j_\al(x)\|_x \,\|\d a(x)\|_x
=
\sum_j p\, \|\phi(x) g_{\al,ij}(x)\,\d a^j_\al(x)\|_x \,\|\d a(x)\|_x.
$$
Taking suprema over $x \in Y$ yields the second inequality (note that
$\|[\D,a]\|_Y\leq\|[\D,a]\|$, as a result of Lemma~\ref{lm:oper-endo}).
\end{proof}

Denote by $A_\D$ the completion of $\A$ in the norm
$\|a\|_\D := \|a\| + \|[\D,a]\|$.

\begin{lemma}
\label{lm:local-Lip}
Let $Y \subset U_\al$ be a compact set on which $a_\al\: Y \to \R^p$
is one-to-one. Then for all $b \in \A_D$ there exists a unique bounded
Lipschitz function $g \: a_\al(Y) \to \C$ such that
$b|_Y = g \circ a_\al|_Y$.
\end{lemma}

\begin{proof}
Since $a_\al$ is one-to-one on~$Y$, the functions 
$a^1_\al,\dots,a^p_\al$ separate the points of~$Y$, so there is a 
unique bounded continuous function $g \in C(a_\al(Y))$ such that
$b|_Y = g(a^1_\al,\dots,a^p_\al)$.

Choose $b_k = g_k(a^1_\al,\dots,a^p_\al) \in \A$ where the $g_k$ are
smooth functions defined on a neighbourhood of $a_\al(Y)$, such that
$$
b_k\bigr|_Y \to b\bigr|_Y  \word{and} 
[\D,b_k]\bigr|_Y \to [\D,b]\bigr|_Y.
$$
This is possible since $A_\D$ is the completion of~$\A$ and functions
of the form $f \circ a_\al$, with $f$ smooth, lie in $\A$ and separate
the points of~$Y$.
Therefore,
$$
\|[\D,b_j] - [\D,b_k]\|_Y
= \sup_{x\in Y} \|[\D,b_j](x) - [\D,b_k](x)\|_x \to 0
\as j,k \to \infty,
$$ 
and $\{g_k\}$ is a Cauchy sequence in the norm
$f \mapsto \sup_{y\in Y} |f(y)| + \|df\|_{Y,\infty}$, by 
Corollary~\ref{cr:local-lipest}. Hence there is a bounded Lipschitz
function $h\: a_\al(Y) \to \C$ such that $h := \lim_k g_k$ in this
norm. Thus,
$$
(h \circ a_\al)\bigr|_Y = \lim_k (g_k \circ a_\al)\bigr|_Y
= \lim_k b_k\bigr|_Y = b\bigr|_Y.
$$
Since $b|_Y = g \circ a_\al$, it follows that
$(g - h) \circ a_\al = 0$ and we have established that $g$ is
Lipschitz on $a_\al(Y)$, with Lipschitz constant bounded above by
$C_Y\,\|[\D,a]\|$.
\end{proof}

\begin{prop}
\label{pr:cts-inverse}
Suppose that $B \subseteq U_\al$ is such that $a_\al \: B \to \R^p$ is
one-to-one. Then the map $a_\al^{-1} \: a_\al(B) \to B$ is continuous
for the metric topology of~$B$ (and thus also for its weak$^*$
topology).
\end{prop}

\begin{proof}
Choose $x,y \in B$ and let $Y$ be any weak$^*$-compact subset of~$B$
with $x,y \in Y$. If $a \in \A$, let $f_a$ be the unique Lipschitz
function on $a_\al(Y)$ with $a\bigr|_Y = f_a(a^1_\al,\dots,a^p_\al)$.
Write $t := a_\al(x)$, $s := a_\al(y)$ in~$\R^p$. We now find, using
Corollary~\ref{cr:local-lipest} and Lemma~\ref{lm:local-Lip}, that
\begin{align}
d(x,y) &= \sup\set{|a(x) - a(y)| : a \in \A,\ \|[\D,a]\| \leq 1}
\nonumber \\
&= \sup\set{|(f_a \circ a_\al)(x) - (f_a \circ a_\al)(y)|
: a \in \A,\ \|[\D,a]\| \leq 1}
\nonumber \\
&= \sup\set{|f_a(t) - f_a(s)| : a \in \A,\ \|[\D,a]\| \leq 1}
\nonumber \\
&\leq \sup\set{C_Y\,\|[\D,a]\| \,|t - s| 
: a \in \A,\ \|[\D,a]\| \leq 1}
\nonumber \\
&= C_Y\, |t - s|,
\label{eq:one-way}
\end{align}
where $|t - s|$ is the Euclidean distance between $t$ and $s$
in~$\R^p$.
\end{proof}

\begin{corl}
\label{cr:twopt-dists}
If $Y \subset U_\al$ is compact and if $x,y \in Y$ satisfy
$a_\al(x)\neq a_\al(y)$, then
$$
d(x,y) \leq C_Y |a_\al(x) - a_\al(y)|.
\eqno \qed
$$
\end{corl}

Note that if $Y,Y'$ are compact subsets of~$U_\al$ with 
$Y \subseteq Y'$, then $C_Y \leq C_{Y'}$ by the proof of
Lemma~\ref{lm:Lip-est}. Thus, in the previous Corollary, the minimal 
value of $C_Y$ is $C_{\{x,y\}}$.

\begin{corl}
\label{cr:homeo-on-image}
If $B \subseteq U_\al$ is such that $a_\al \: B \to \R^p$ is
one-to-one, then the map $a_\al|_B$ is a homeomorphism onto its image
(for either the metric or the weak$^*$ topology of~$B$).
\end{corl}

\begin{proof}
The map $a_\al\: B \to a_\al(B) \subset \R^p$ is continuous for the
weak$^*$ topology on~$B$ since each $a_\al^j$ lies in $\A$, and thus
it is also continuous for the metric topology on~$B$. The estimate
\eqref{eq:one-way} shows that the inverse map
$a_\al^{-1}\: a_\al(B) \to B$ is continuous for the metric topology
on~$B$, and \textit{a posteriori} for the weak$^*$ topology.
\end{proof}

\begin{corl}
\label{cr:locally-open}
If $V \subseteq U_\al$ is an weak$^*$-open subset such that
$a_\al \: V \to \R^p$ is one-to-one, then $a_\al(V)$ is an open
subset of~$\R^p$.
\end{corl}

\begin{proof}
By Lemma~\ref{lm:nearly-open}, the set $a_\al(V)$ has nonempty 
interior in~$\R^p$. Now $a_\al(V) \less \intr a_\al(V)$ is the 
boundary of $\intr a_\al(V)$ in the relative topology of $a_\al(V)$.
Since $a_\al$ is a homeomorphism from $V$ onto $a_\al(V)$, it cannot 
map interior points to this boundary. Thus, since $V$ is open, this
boundary is empty, and hence $a_\al(V) = \intr a_\al(V)$ is open 
in~$\R^p$.
\end{proof}

\begin{lemma}
\label{lm:local-tops}
On any subset $V \subset U_\al$ on which $a_\al$ is one-to-one, the
weak$^*$ and metric topologies coincide.
\end{lemma}

\begin{proof}
By Lemma~\ref{lm:local-Lip}, the restriction of any function
$a \in \A$ to a compact subset $Y \subset V$ may be written as
$a|_Y = (f_a \circ a_\al)|_Y$, where $f_a$ is a bounded Lipschitz
function on $a_\al(Y)$. We need only show that convergence of a
sequence $V \ni x_m \to y$ for the weak$^*$ topology implies the
convergence $x_m \to y$ in the metric topology. Choose such a
weak$^*$-convergent sequence $\{x_m\}$; without loss of generality, we
may suppose that it is contained in a compact subset $Y \subset V$
such that each $x_m \in Y$ and hence also $y \in Y$.

Weak$^*$ convergence of the sequence $x_m$ means that
$|a(x_m) - a(y)| \to 0$ for all $a \in \A$. This implies that
$$
\frac{|a(x_m) - a(y)|}{\|[\D,a]\|} \to 0  \as  m \to \infty
$$
for all $a \in \Dom\d$ with $\d a = [\D,a] \neq 0$. By
Lemma~\ref{lm:local-Lip}, this is equivalent to
$$
\frac{|(f_a\circ a_\al)(x_m) - (f_a\circ a_\al)(y)|}{\|[\D,a]\|} \to 0
\as  m \to \infty
$$
for all such $a$. Using Corollary~\ref{cr:local-lipest}, this is
equivalent to
$$
|a_\al(x_m) - a_\al(y)| \to 0   \as  m \to \infty.
$$
The metric convergence follows:
\begin{align*}
d(x_m,y) &= \sup\set{|a(x_m) - a(y)| : \|[\D,a]\| \leq 1}
\\
 &= \sup\set{|(f\circ a_\al)(x_m) - (f\circ a_\al)(y)|
 : \|df\|_\infty \leq C_Y }
\\
&\leq C_Y\, |a_\al(x_m) - a_\al(y)| \to 0  \as  m \to \infty,
\end{align*}
on invoking Corollary~\ref{cr:local-lipest} once more.
\end{proof}

\section{Point-set properties of the local coordinate charts}
\label{sec:point-sets}

We must analyze the possible failure of $a_\al$ to be one-to-one, by
using some point-set topology. Some of this extra effort is due to the
arbitrariness of the Hochschild cycle~$\cc$. For example, consider the
manifold $\Sf^2$, with volume form
$x \,dy \w dz + y \,dz \w dx + z \,dx \w dy$. Each $U_\al$ (the subset
where $a^0_\al$ is nonzero, and $da^1_\al$, $da^2_\al$ are linearly
independent) consists of two open hemispheres, missing only an
equator. For instance, there is a chart domain consisting of
$\Sf^2 \less \{z = 0\}$ with local coordinates $(x,y)$; but these
``coordinates'' are actually two-to-one on this domain. This rather
simple problem can be addressed by restricting the coordinate maps to
the connected components of the $U_\al$, thereby increasing the number
of charts. Ultimately, in the next section, we shall show that this 
is the worst that can happen.

We begin by considering some set-theoretic properties of the map
$a_\al\: U_\al \to \R^p$ for a fixed~$\al$, with a view to showing
that $a_\al$ is at worst finitely-many-to-one on an open dense subset
of~$U_\al$. (Eventually, this behaviour will be improved to local
injectivity.)

The linear functional
$$
\mu_\D(a) := \ncint a\,\Dreg^{-p}
$$
on~$\A$ extends to a continuous linear functional on~$A$, and by the
Riesz representation theorem it is given by a regular Borel measure
that we also denote by~$\mu_\D$.

\begin{lemma}
\label{lm:Borel-meas}
The measure $\mu_\D$ has no atoms.
\end{lemma}

\begin{proof}
By \cite[Lemma~14, p.~408]{Royden}, we can decompose $X$ as a 
disjoint union $X = X' \uplus C$ where $C$ is countable, its closure 
$\ol{C}$ is the support of the atomic part of~$\mu_\D$, and $X'$ has 
no isolated points. Recalling that
$\braket{\xi}{\eta} = \mu_\D(\pairing{\xi}{\eta})$ for
$\xi,\eta \in \H_\infty$, we get a corresponding Hilbert-space
decomposition
$$
\H = L^2(X,S) = L^2\bigl(X',S|_{X'}\bigr) \oplus L^2\bigl(C,S|_C\bigr).
$$
If $\supp\xi = \{x\} \subseteq C$, then $a\xi = a(x)\xi$ for any 
$a \in \A$, so that $\C\xi$ would be a $1$-dimensional
subrepresentation of $\pi(A)$ and thus $\pi(A)$ would contain a
nontrivial projector, contradicting Corollary~\ref{cr:spec-conn}.
\end{proof}

We can write $\mu_\D$ as a sum of finite measures $\mu_{\D,\al}$
concentrated on each $U_\al$~by
$$
\ncint a\,\Dreg^{-p} = \ncint \Ga^2\,a\,\Dreg^{-p}
= \sum_{\al=1}^n \ncint
\Ga\, a a^0_\al \,\d a_\al^1 \dots \d a_\al^p \,\Dreg^{-p}
=: \sum_{\al=1}^n \mu_{\D,\al}(a)
$$
since the third expression depends only on the skew-symmetrization of
the $\d a_\al^j$, by the proof of Proposition~\ref{pr:lin-indp}, and
each skew-symmetrized $\Ga'_\al$ vanishes off the respective~$U_\al$.

We can transfer these measures to the several $a_\al(U_\al)$ by 
setting
$$
\la_{\D,\al}(f) := \mu_{\D,\al}(f(a_\al^1,\dots,a_\al^p))
\word{for}  f \in C_c(a_\al(U_\al)).
$$
Each $\la_{\D,\al}$ is a nonatomic regular Borel measure on 
$a_\al(U_\al) \subset \R^p$. Its Lebesgue decomposition
$$
\la_{\D,\al} = \la_{\D,\al}^s + \la_{\D,\al}^\ac
$$
provides measures $\la_{\D,\al}^s$ and $ \la_{\D,\al}^\ac$ that are
respectively singular and absolutely continuous with respect to the
Lebesgue measure on $a_\al(U_\al)$. The singular part $\la_{\D,\al}^s$
is concentrated on a set of Lebesgue measure zero, whereas the
absolutely continuous part $\la_{\D,\al}^\ac$ has support
$\ol{a_\al(U_\al)}$.

\begin{defn}
Define $n_\al \: U_\al \to \{1,2,\dots,\infty\}$ by
$$
n_\al(y) := \#\bigl( a_\al^{-1}(a_\al(y)) \cap U_\al \bigr),
$$
where $\#$ denotes cardinality (all infinite cardinals being treated
simply as~$\infty$).
\end{defn}

\begin{lemma}
\label{lm:infty-mult-nonfat}
The set $n_\al^{-1}(\infty)$ of infinite-multiplicity points has empty
interior in $U_\al$.
\end{lemma}

\begin{proof}
Let $\spec_\ac(a_\al) \subseteq \Delta_\al$ be the absolutely
continuous joint spectrum of $a_\al = (a_\al^1,\dots,a_\al^p)$,
regarded as $p$ commuting selfadjoint operators on~$\H$.

Over $\spec_\ac(a_\al)$, the multiplicity function $m_\al$ of (the
representation of) $a_\al$ is $L^1$ with respect to Lebesgue measure.
Without changing the measure class, we may take $m_\al$ to be
precisely the multiplicity of $a_\al$, namely,
$m_\al(t) = N\,n_\al(a_\al^{-1}(t) \cap U_\al)$ for
$t \in \spec_\ac(a_\al)$, where $N$ is the rank of the bundle~$S$. 
This is well defined since $n_\al$ is constant on
$a_\al^{-1}(t) \cap U_\al$.

Suppose that $V \subset a_\al(n_\al^{-1}(\infty))$ is a Borel subset
of $a_\al(U_\al)$. If $V$ had positive Lebesgue measure, it would
follow that
\begin{equation}
\int_{\spec_\ac(a_\al)} m_\al(t) \,d^pt \geq \int_V m_\al(t) \,d^pt
= \infty.
\label{eq:mult-too-big}
\end{equation}
However, by \cite{Voiculescu}, discussed also in
\cite[Sec.~IV.2.$\dl$]{Book}, there is an equality
$$
\int_{\spec_\ac(a_\al)} m_\al(t) \,d^pt = C'_p (k_p^-(a_\al))^p,
$$
where $C'_p$ is a constant depending only on~$p$, and Voiculescu's
modulus $k_p^-(a_\al)$ is a positive number which measures the size of
the joint absolutely continuous spectrum. This number is finite,
since by \cite[Prop.~IV.2.14]{Book}, see also
\cite{CareyPRSone,Polaris} for similar estimates:
$$
k_p^-(a_\al) \leq C_p \max_j\|[\D, a^j_\al]\|
\biggl(\, \ncint \Dreg^{-p} \biggr)^{1/p} < \infty.
$$
This contradicts \eqref{eq:mult-too-big}; consequently, $V$ is a
Lebesgue nullset. Since any subset of~$a_\al(U_\al)$ on which
$\la_{\D,\al}^s$ is concentrated also has Lebesgue measure zero, and
elsewhere $m_\al$ is finite, we conclude that
$a_\al(n_\al^{-1}(\infty))$ is a Lebesgue nullset, and in particular
it has empty interior. Now Lemma~\ref{lm:nearly-open} entails that
$n_\al^{-1}(\infty)$ has empty interior in $U_\al$.
\end{proof}

The remainder of this section proves that the topological structure
of~$U_\al$ is sufficiently nice for us to deploy our geometric tools
in the following section. These tools will prove the local injectivity
of $a_\al\: U_\al \to \R^p$.

{}From now until Theorem~\ref{thm:lip-mfld}, we fix $\al$ and work
solely within $U_\al$ using the relative weak$^*$ topology. Thus, all
closures, interiors and boundaries are taken in this relative weak$^*$
topology, unless specifically noted. Also, for
$E \subseteq \Delta_\al \subset \R^p$, ``$a_\al^{-1}(E)$'' will mean
$a_\al^{-1}(E) \cap U_\al$. \textit{Thus, up until
Theorem~\ref{thm:lip-mfld}, we restrict the universe to~$U_\al$.}

\begin{defn}
\label{df:mult-k}
Consider the following subsets of $U_\al$, for $k = 1,2,\dots$:
\begin{align}
D_k &:= n_\al^{-1}(\{1,\dots,k\}),
\nonumber \\
E_k &:= \intr D_k,
\nonumber \\
N_{k+1} &:= U_\al \less D_k = n_\al^{-1}(\{k+1,\dots,\infty\}),
\nonumber \\
W_k &:= \intr(n_\al^{-1}(k)).
\label{eq:mult-k}
\end{align}
\end{defn}

\begin{rmk}
If $Z \subset U_\al$, its closure in the relative (weak$^*$) topology
of $U_\al$ is ${\ol Z}^X \cap U_\al$. Elements of this closure are
limits of sequences in~$Z$, since $X$ is metrizable because the
$C^*$-algebra $A$ is assumed to be separable.
\end{rmk}

\begin{lemma}
\label{lm:uppercts}
The set $N_{k+1}$ is open and the set $D_k$ is closed in $U_\al$.
\end{lemma}

\begin{proof}
If $D_k$ is finite, there is nothing to prove; otherwise, choose any
convergent sequence $\{x_i\} \subset D_k$ with $x_i \to x \in U_\al$.
Each element $x_i$ has multiplicity $\leq k$, and so
$a_\al^{-1}(a_\al(\{x_i\}))$ consists of at most $k$ sequences with at
most $k$ limit points $y_1,\dots,y_m$, $m \leq k$, since $U_\al$ is
Hausdorff. Thus $x$ is one of the $y_j$, say $x = y_1$, and
$a_\al(y) = a_\al(x)$ for $y \in \ol D_k$ if and only if
$y \in \{y_1,\dots,y_m\}$. (Notice that $U_\al$ can contain no
isolated points, by Corollary~\ref{cr:spec-conn}.) Hence the limit
point $x$ has multiplicity $m \leq k$, and so $x \in D_k$.
\end{proof}

Notice that $D_k$ is nonempty for some finite~$k$, since
$n_\al^{-1}(\infty) \neq U_\al$. The next Proposition shows that the
awkward possibility of $n_\al^{-1}(\infty)$ being dense cannot occur.
We require a preparatory lemma.

\begin{lemma}
\label{lm:infty-dense}
If $n_\al^{-1}(\infty)$ were dense in $U_\al$, then each $N_{k+1}$
would be an open dense subset of $U_\al$, every neighbourhood of an
infinite-multiplicity point in $n_\al^{-1}(\infty)$ would contain
elements of $D_k$ for arbitrarily large~$k$, and every neighbourhood
of a finite-multiplicity point in some $D_k$ would contain infinitely
many points in $n_\al^{-1}(\infty)$.
\end{lemma}

\begin{proof}
The first statement is clear, since
$n_\al^{-1}(\infty) \subseteq N_{k+1}$ for each finite~$k$.

By Lemma~\ref{lm:infty-mult-nonfat}, the union
$\bigcup_{k=1}^\infty D_k$ is dense in~$U_\al$ and, by the proof of 
Lemma~\ref{lm:uppercts}, no infinite-multiplicity point can be the 
limit of a sequence contained in any~$D_k$.

The last statement is just the assumed density
of~$n_\al^{-1}(\infty)$.
\end{proof}

\begin{prop}
\label{pr:infty-mult-lean}
The subset $n_\al^{-1}(\infty)$ is nowhere dense in~$U_\al$.
\end{prop}

\begin{proof}
Let $Y$ be a compact subset of $U_\al$ with nonempty interior.
Suppose, \textit{arguendo}, that $n_\al^{-1}(\infty) \cap \intr Y$ is
dense in $\intr Y$. Since $a_\al$ cannot then be one-to-one on~$Y$, we
can choose a finite-multiplicity value $t \in a_\al(\intr Y)$ and two
distinct points $y,y' \in Y$ such that $a_\al(y) = a_\al(y') = t$. If
necessary, we can add a small compact neighbourhood of~$y'$ to~$Y$.

Using Lemma~\ref{lm:infty-dense}, we can find a sequence
$\{t_m\} \subset a_\al(Y \cap n_\al^{-1}(\infty))$ such that
$|t - t_m| < \eps_m$ for all $m$, where $\eps_m \to 0$ with
$0 < \eps_m < (2C_Y)^{-1} d(y,y')$ for all~$m$, where $C_Y$ is the
constant appearing on the right hand side of the
estimate~\eqref{eq:local-lipest}.

Suppose that for all~$m$, all points $z \in a_\al^{-1}(t_m) \cap Y$
satisfy $d(y,z) \leq C_Y \eps_m$. Then
$$
d(y',z) \geq \bigl| d(y,y') - d(y,z) \bigr|
\geq d(y,y') - C_Y \eps_m > C_Y \eps_m.
$$
Consequently, on applying Corollary~\ref{cr:twopt-dists} to the pair
of points $y',z \in Y$, we obtain
\begin{equation}
C_Y \eps_m < d(y',z) \leq C_Y \bigl| a_\al(y') - a_\al(z) \bigr|
= C_Y |t - t_m| < C_Y \eps_m,
\label{eq:toobig-toosmall}
\end{equation}
a contradiction. On the other hand, if there is some~$m$ and some
$z \in a_\al^{-1}(t_m) \cap Y$ such that $d(y,z) > C_Y \eps_m$, then
we reach the same impasse on replacing $y'$ by $y$ in 
\eqref{eq:toobig-toosmall}.

We conclude that no such sequence $t_m \to t$ can exist, so that in
particular $y$ has a neighbourhood excluding $n_\al^{-1}(\infty)$.
That is to say, $n_\al^{-1}(\infty) \cap \intr Y$ is not dense in
$\intr Y$. By the arbitrariness of~$Y$, $n_\al^{-1}(\infty)$ is
nowhere dense in~$U_\al$.
\end{proof}

\begin{corl}
\label{cr:infty-mult-lean}
The map $a_\al$ is finitely-many-to-one on an open dense subset
of~$U_\al$.
\qed
\end{corl}

\begin{lemma}
\label{lm:bdry-k}
Within $U_\al$, the following relations hold:
$$
n_\al^{-1}(k) \less W_k \subset \del W_k \cup \del N_{k+1}, \qquad
n_\al^{-1}(k) \less \ol W_k \subset \del N_{k+1}, \qquad
\del D_k = \del N_{k+1},
$$
and thus $U_\al = \intr D_k \cup \ol N_{k+1}$. Moreover,
$$
\del W_k \subseteq D_k  \word{and}
n_\al^{-1}(\infty) \cap \bigcup_{k=1}^\infty \ol W_k = \emptyset.
$$
\end{lemma}

\begin{proof}
Consider the set $H_k := U_\al \less 
\bigl( D_{k-1} \cup \ol W_k \cup \ol N_{k+1} \bigr)$.
It is open in $U_\al$, and $H_k \subset n_\al^{-1}(k)$. However,
$H_k \cap \intr(n_\al^{-1}(k)) = H_k \cap W_k = \emptyset$, entailing
$H_k = \emptyset$. Therefore
$U_\al = D_{k-1} \cup \ol W_k \cup \ol N_{k+1}$, and so
$$
n_\al^{-1}(k) \less W_k \subset D_{k-1} \cup \del W_k \cup \ol N_{k+1}.
$$
However, both $D_{k-1}$ and $N_{k+1}$ are disjoint from
$n_\al^{-1}(k)$ by definition, so
$$
n_\al^{-1}(k) \less W_k \subset \del W_k \cup \del N_{k+1}.
$$
Taking the closure of $W_k$ in $U_\al$ gives the second relation. The
third relation follows on recalling that
$U_\al \less N_{k+1} = D_k$.

It is clear that $\del W_k \subseteq D_k$, since $W_k \subset D_k$ 
and $D_k$ is closed. The last relation follows from 
$\ol W_k \subset D_k$, since $n_\al^{-1}(\infty)$ is disjoint from
each~$D_k$.
\end{proof}

\begin{lemma}
\label{lm:w-and-d}
The interiors of $D_k$ and of $\bigcup_{j=1}^k \ol W_j$ coincide:
$E_k = \intr \bigl( \bigcup_{j=1}^k \ol W_j \bigr)$.
\end{lemma}

\begin{proof}
Take $x \in E_k$. Then every open neighbourhood of $x$ in $D_k$ must
meet some $W_j$ with $j \leq k$, so that
$x \in \bigcup_{j=1}^k \ol W_j$. Hence
$E_k \subseteq \bigcup_{j=1}^k \ol W_j$.

On the other hand, since each $D_j$ is closed in $U_\al$, we get
$\bigcup_{j=1}^k \ol W_j \subseteq \bigcup_{j=1}^k D_j
= D_k$. The inclusion
$\intr \bigl( \bigcup_{j=1}^k \ol W_j \bigr) \subseteq E_k$
follows at once.
\end{proof}

These topological preliminaries now allow us to prove three basic 
results about the mapping $a_\al$ and the set~$U_\al$.

\begin{prop}
\label{pr:first-cover}
There exists a weak$^*$-open cover $\{Z_j\}_{j\geq 1}$ of
$\bigcup_{k=1}^\infty W_k \subseteq U_\al$ such that
$a_\al \: Z_j \to a_\al(Z_j) \subset \R^p$ is a homeomorphism and an 
open map, for each $j \geq 1$.
\end{prop}

\begin{proof}
First observe that $W_1$ is an open subset of $U_\al$ such that
$a_\al \: W_1 \to a_\al(W_1)$ is an open homeomorphism. This follows
from Corollaries \ref{cr:homeo-on-image} and~\ref{cr:locally-open} and
the definition of~$W_1$.

Choose $x \in W_k$, $k > 1$. Since $a_\al^{-1}(a_\al(x))$ consists of
$k$ points of $W_k$ and $W_k$ is Hausdorff, we may choose $k$ disjoint
open subsets of $W_k$, each containing precisely one of the preimages
$x_1,\dots,x_k$ of $a_\al(x)$. Call these sets $V_1,\dots,V_k$, and
suppose that for some $j = 1,\dots,k$, the map $a_\al \: V_j \to \R^p$
is not one-to-one. Then there exist $z_1,z_2 \in V_j$ with
$z_1 \neq z_2$, such that $a_\al(z_1) = a_\al(z_2)$. Again we can
separate these two points, and thereby suppose that only one of
$z_1,z_2$ lies in $V_j$. If for every $x\in W_k$ we can repeat this
process finitely often to obtain neighbourhoods of each of the $x_j$
on which $a_\al$ is one-to-one, we are done.

On the other hand, if in this manner we always find pairs of distinct
points $z_{1m} \neq z_{2m}$ from successively smaller neighbourhoods
of~$x_j$, we thereby obtain two sequences, both converging to~$x_j$.
We may alternate the labelling of each pair, if necessary, so that
both sequences are infinite. Then for $i = 1,2$, $a_\al(z_{im})$
converges to $a_\al(x_j)$ for $j = 1,\dots,k$, and so in any open
neighbourhood $N$ of $a_\al(x_j)$ there are infinitely many such
$a_\al(z_{im})$. Now take a neighbourhood $V_l$ of each $x_l$, for $l
\neq j$, such that $a_\al(V_l)$ maps onto $N$ (this is possible; we
started with disjoint neighbourhoods of the $x_l$ which map onto a
neighbourhood of $a_\al(x_j)$). Then there is a sequence $\{y_{lm}\}
\subset V_l$ mapping onto $a_\al(z_{im})$. A quick count now shows
that each element of the sequence $a_\al(z_{im})$ has (at least) $k+1$
preimages in $W_k$: contradiction.

Hence, for all $x\in W_k$ we may find a neighbourhood $V$ of $x$ in
$W_k$ such that $a_\al \: V \to a_\al(V)$ is one-to-one, and thus is 
an open homeomorphism for the weak$^*$ topology on~$V$, by
Corollary~\ref{cr:homeo-on-image}. In this way we obtain a open cover
of the (locally compact, metrizable) set $\bigcup_{k=1}^\infty W_k$ of
the desired form; now let $\{Z_j\}_{j\geq 1}$ be an enumeration of a
countable subcover. Notice that we have also shown that each $Z_j$ is
included in some $W_k$, with $k = k(j)$.
\end{proof}

The next Proposition is the most crucial consequence of the Lipschitz
functional calculus. It allows us to extend the homeomorphism
$a_\al \: Z_j \to a_\al(Z_j)$ to the closure $\ol Z_j$ \emph{as a
homeomorphism}. It is essential in the next Proposition that we use
the relative topology of~$U_\al$.

\begin{prop}
\label{pr:local-homeos}
For each open set $Z_j$ as in Proposition~\ref{pr:first-cover}, the
function $a_\al$ extends to a homeomorphism (for both the metric and
weak$^*$ topologies) $a_\al \: \ol Z_j \to a_\al(\ol Z_j)$. The same
is true when $Z_j$ is replaced by any open set $V \subset U_\al$ on
which $a_\al$ is one-to-one. In particular, the weak$^*$ and metric
topologies agree on $\ol V$.
\end{prop}

\begin{proof}
Take any $b \in \A$. Then by Lemma~\ref{lm:local-Lip}, for
any compact $Y \subset Z_j$ there is a unique bounded Lipschitz
function $g \: a_\al(Y) \to \C$ such that $b|_Y = g \circ a_\al|_Y$.
Now we can cover $Z_j$ by open sets which are interiors of compact
subsets, and so obtain many function representations of $b$ on these
sets. By uniqueness, they agree on overlaps, and $b = g \circ a_\al$
for each of these local representations. Hence we arrive at a single
function $g \: a_\al(Z_j) \to \C$ with
$b|_{Z_j} = g \circ a_\al|_{Z_j}$. The function $g$ is bounded since
$b$ is bounded, but might be only locally Lipschitz (since
$a_\al^{-1}$ might only be locally Lipschitz).

To proceed, suppose first that $g$ is a $C^1$-function on
$a_\al(Z_j)$. The proof of Lemma~\ref{lm:Lip-est}, in particular
\eqref{eq:Lip-est2}, and Corollary~\ref{cr:local-lipest} guarantee
that for any subset $Y$ of $Z_j$ with compact closure contained in
$\ol Z_j$,
$$
\sup_{x\in Y} |\del_jg(a_\al(x))|
\leq \sup_{x\in\ol Y} B_j(x) \,\|[\D,b](x)\|_x < \infty.
$$
The finiteness of the right hand side follows because
$\|[\D,b](x)\|_x \leq \|[\D,b]\|$, and $B_j$ is continuous on all of
$U_\al$, so that $B_j$ is bounded on~$\ol Y$.

To see that $g$ extends to a locally Lipschitz function on
$a_\al(\ol Z_j)$, we argue as follows. If
$t \in a_\al(\ol Z_j)$, and if $\{t_n\}$ and
$\{t'_n\}$ are two sequences in $a_\al(Z_j)$ such that $t_n \to t$
and $t'_n \to t$ in $a_\al(U_\al)$, then
$$
|g(t_n) - g(t'_n)| \leq C\,|t_n - t'_n|,  \word{where}
C = \sup_{x\in\ol Y} B_j(x) \,\|[\D,b]\|.
$$
For the estimate, since $Z_j \subseteq W_k$ for a suitable~$k$, it is
enough take $Y$ to be the union of the $k$ preimages of the sequences
$\{t_n\}$ and $\{t'_n\}$ and their limits, which is compact in
$U_\al$. Thus $\tilde g(t) := \lim_{n\to\infty} g(t_n)$ is well
defined, and coincides with $g(t)$ whenever $t \in a_\al(Z_j)$
already. The upshot is a bounded continuous function
$\tilde g \: a_\al(\ol Z_j) \to \R$.

Its Lipschitz norm on any compact subset
$Y \subset \ol Z_j$ satisfies
$\|d\tilde g\|_Y \leq C_Y\,\|[\D,b]\|$, and so $\tilde g$ is locally 
Lipschitz. The continuity of $\tilde g$ and $b$ yields
$b = \tilde g \circ a_\al$ over the set $\ol Z_j$.

For an arbitrary $b \in \A$, we may remove the assumption that $g$ be
$C^1$ on $a_\al(Z_j)$ by approximating $b$ by a sequence $\{b_r\}$ in
the norm $a \to \|a\| + \|[\D,a]\|$, where
$b_r|_{Z_j} = g_r \circ a_\al|_{Z_j}$ with each $g_r$ being~$C^1$. For
any subset $Y \subset Z_j$ with compact closure, we get an estimate
$$
\sup_{x\in\ol Y} |\del_jg_r(a_\al(x)) - \del_jg_s(a_\al(x))|
\leq \sup_{x\in\ol Y} B_j(x) \,\|[\D, b_r - b_s]\|,
$$
and thus $\{g_r\}$ is a Cauchy sequence in the Lipschitz norm of
$a_\al(\ol Y)$. Each $g_r$ extends to a locally Lipschitz
function $\tilde g_r$ on $a_\al(\ol Z_j)$, and
these converge uniformly on compact subsets to a locally Lipschitz
function $\tilde g$ satisfying $b = \tilde g \circ a_\al$ over
$\ol Z_j$.

Since we can thereby extend the function representation for all
functions $b \in \A$, we conclude that $a_\al$ separates points of
$\ol Z_j$, and thus it is one-to-one on this set. By
Corollaries \ref{cr:homeo-on-image} and~\ref{cr:locally-open},
$a_\al\: \ol Z_j \to \R^p$ is an open map, which is a
homeomorphism onto its image. By Lemma~\ref{lm:local-tops}, the
weak$^*$ and metric topologies agree.
\end{proof}

Here is a first and critical consequence of
Proposition~\ref{pr:local-homeos}.

\begin{lemma}
\label{lm:pancakes}
For each $k\geq 1$, the set $W_k$ is a disjoint union of $k$ open
subsets $W_{1k},\dots,W_{kk}$, on each of which $a_\al$ is injective,
such that $a_\al(W_{jk}) = a_\al(W_{j'k})$ for $j,j' = 1,\dots,k$.
\end{lemma}

\begin{proof}
Choose any point $x\in W_k$, and choose disjoint open sets
$V_j \subset W_k$, for $j = 1,\dots,k$, each containing precisely one
of the preimages of $x$ in $a_\al^{-1}(a_\al(x))$. We may suppose by
Proposition~\ref{pr:first-cover} that $a_\al$ is one-to-one on each
$V_j$, and on replacing the $V_j$ by the components of
$a_\al^{-1}\bigl( \bigcap_j a_\al(V_j) \bigr)$, we may suppose also
that the various $V_j$ have the same image $a_\al(V_j)$ in $\R^p$. By
Proposition~\ref{pr:local-homeos}, $a_\al$ extends to a homeomorphism
on each $\ol V_j$.

Suppose first that there exists $z\in \ol V_1 \cap \ol V_l \cap W_k$
for some $l \geq 2$. Then, since $a_\al$ is a homeomorphism on each
$\ol V_j$, the set
$a_\al^{-1}(a_\al(z)) \cap \bigcup_{j=1}^k \ol V_j$ consists of fewer
than $k$ points. So there must be some
$z' \in a_\al^{-1}(a_\al(z)) \cap W_k$ with $z' \notin \ol V_j$ for
$j = 1,\dots,k$. The open set $W_k \less \bigcup_{j=1}^k \ol V_j$
includes an open neighbourhood $O$ of~$z'$. Now $a_\al^{-1}(a_\al(O))$
is an open set in $W_k$ meeting each $V_j$, since $z \in \del V_1$ and
$a_\al(V_1) = a_\al(V_j)$ for each~$j$. Consider a sequence
$\{t_m\} \subset a_\al(O) \cap a_\al(V_j)$ with
$t_m \to a_\al(z) = a_\al(z')$. Then $a_\al^{-1}(\{t_m\})$ meets each
$V_j$ and~$O$. Hence for $m$ sufficiently large,
$n_\al(a_\al^{-1}(t_m)) > k$, which is impossible within $W_k$.

We conclude that the relatively closed subsets $\ol V_j \cap W_k$ of
$W_k$ are disjoint. Since $W_k$ is metrizable (by
Remark~\ref{rk:metrizable}) and is thus a normal topological space,
there are disjoint open sets $U_1,\dots,U_k \subset W_k$ with
$\ol V_j \cap W_k \subset U_j$.

Now choose a point $z$ on the boundary of $V_1$ lying in $W_k$ (if
there is no such point, then $V_1$ is a union of connected components
of~$W_k$, and so too are the other $V_j$, and we are done). Choose an
open neighbourhood $U \subset W_k$ of $z$ on which $a_\al$ is
one-to-one. Note that $a_\al$ is one-to-one on $U \cap V_1$, on $U$,
and on $V_1$. Let $U' := U_1 \cap U$; we claim that
$a_\al\: U' \cup V_1 \to \R^p$ is one-to-one. For if not, there would
exist $y \in U' \less V_1$ and $y' \in V_1 \less U'$ such that
$a_\al(y) = a_\al(y')$. However, $y \notin V_1$ and
$a_\al^{-1}(a_\al(y'))$ consists of $k$ points already, so this forces
$y \in V_j \cap a_\al^{-1}(a_\al(y'))$ for some $j > 1$; otherwise $y$
would have multiplicity at least $k + 1$, contradicting $y \in W_k$.
However, $y \in U' \less V_1 \subset U_1$, which forbids $y \in V_j$
for $j > 1$. The upshot is that $a_\al^{-1}(a_\al(U' \cup V_1))$ is a
union of $k$ disjoint open sets on each of which $a_\al$ is
one-to-one. It is clear that $U$ may be chosen such that their images
under $a_\al$ coincide.

By this argument, we may continue this process by taking a boundary
point of $U' \cup V_1$ within $W_k$, finding a neighbourhood on which
$a_\al$ is one-to-one, and deducing that $a_\al$ is one-to-one on the
union. In this way we cover the entire connected component of $W_k$ in
which $x$ lies. Thus $W_k$ is a disjoint union of $k$ open subsets
$W_{1k},\dots,W_{kk}$, and by construction we see that
$a_\al(W_{jk}) = a_\al(W_{j'k})$ for all~$j,j'$.
\end{proof}

\begin{corl}
\label{cr:little-mfld}
Each set $W_{jk}$ is a smooth manifold with the coordinate map
$a_\al|_{W_{jk}}$.
\qed
\end{corl}

\begin{rmk}
If there is only a single nonempty $W_k$, and $U_\al = W_k$,
Lemma~\ref{lm:pancakes} shows that we are in the situation of the
$2$-sphere described at the beginning of this section, where $k = 2$.
\end{rmk}

Consider now $D_k$, the subset of multiplicity at most~$k$. Clearly,
$\bigcup_{j=1}^k W_j \subseteq E_k = \intr D_k$.

If it were true that for all $k$ one could find some $m > k$ with
$D_k \subseteq \intr D_m$, then one could deduce from Lemmas
\ref{lm:bdry-k} and \ref{lm:w-and-d} that
$U_\al = \intr\bigl(\bigcup_k \ol W_k \bigr) \cup n_\al^{-1}(\infty)$.
For lack of such a guarantee (at present), we name the following
subsets where the multiplicities may be troublesome. We shall
eventually show that these subsets are empty.

\begin{defn}
\label{df:bad-points}
Consider the following subsets of $U_\al$:
$$
B_k := \set{x\in D_k \less E_k : x\notin E_m \text{ for all } m > k},
\quad
B_{(\al)} := \bigcup_{k\geq 1} B_k,  \qquad
C_k := E_k \less \biguplus_{j=1}^k W_j.
$$
The set $B_k$ consists of (some) boundary points of $D_k$, while $C_k$
collects interior points $x$ of $D_k$, if any, that lie on the
boundary of some $W_j$ with $j \leq k$. The points $C_k$ will be 
branch points of the ``branched manifold'' 
$\bigcup_{k\geq 1} E_k \subseteq U_\al$.
\end{defn}

\begin{lemma}
\label{lm:bad-mult}
If the multiplicity is bounded on $U_\al$, then
$n_\al^{-1}(\infty) = B_{(\al)} = \emptyset$. If the multiplicity is
unbounded, then $x \in D_k \less E_k$ lies in $B_k$ if and only if
$x = \lim_{n\to\infty} x_n$ for some sequence $\{x_n\}$ satisfying
$n_\al(x_n) \to \infty$.
\end{lemma}

\begin{proof}
The first statement is obvious: if the multiplicity is bounded, by $m$
say, then $D_m = U_\al$ and every $D_k$ is included in $D_m = E_m$.

So suppose that the multiplicity is unbounded, and take $x\in B_k$.
Then $x \in D_m \less E_m$ for all $m > k$. Since
$\del D_m = \del N_{m+1}$, for each $m$ we can find a sequence
$\{y_j^{(m)}\} \subset N_{m+1}$ with $y_j^{(m)} \to x$. By a diagonal
argument, we can now construct a sequence $\{y_j\}$ converging to $x$
and such that $n_\al(y_j) \to \infty$.

Conversely, suppose that $n_\al(x) = k$ and that there is a sequence
$\{x_n\}$ with $x_n \to x$ and $n_\al(x_n) \to \infty$. Then
$x \notin E_k = \intr D_k$ since any sequence converging to an
$x \in E_k$ would have multiplicity eventually bounded by~$k$. Thus
$x \in D_k \less E_k$ and $x \notin E_m$ for $m > k$.
\end{proof}

\begin{lemma}
\label{lm:bad-sets}
Within $U_\al$, the set $B_{(\al)}$ is closed with empty interior.
\end{lemma}

\begin{proof}
If $\{x_j\} \subset B_{(\al)}$ is a sequence with a limit
$z \in U_\al$, choose for each $j$ a sequence $\{x_{jn}\}$ of
unbounded multiplicity such that $x_{jn} \to x_j$. Passing to a
subsequence if necessary, the sequence $\{x_{jj}\}$ converges to $z$
and $n_\al(x_{jj}) \to \infty$. Hence $B_{(\al)}$ is closed.

By Lemma~\ref{lm:bad-mult}, any $x \in B_{(\al)}$ lies in the 
boundary of $\bigcup_k E_k$, and thus $B_{(\al)}$ has empty interior.
\end{proof}

\begin{corl}
\label{cr:bad-sets}
The set $\bigcup_{k\geq 1} E_k$ is an open dense subset of $U_\al$.
Moreover,
$$
U_\al
= \bigcup_{k=1}^\infty E_k \cup B_{(\al)} \cup n_\al^{-1}(\infty).
\eqno\qed
$$
\end{corl}

Within $E_k=\intr D_k$ we have possible branch points where 
two or more `sheets' of one or more $W_j$, $j\leq k$, may 
have common boundary. We characterise these points in $C_k$ next.

\begin{lemma}
\label{lm:ugly-sets}
If $x \in C_k$, then $n_\al(x) < k$.
\end{lemma}

\begin{proof}
Let $x \in C_k = E_k \less \biguplus_{j=1}^k W_j$ and suppose that
$n_\al(x) = k$. Since
$U_\al = E_k \cup \ol N_{k+1} = D_k \cup N_{k+1}$ by
Lemma~\ref{lm:bdry-k}, then not all neighbourhoods of~$x$ can meet
$\ol N_{k+1}$, since $x \in E_k$ entails $x \notin \del D_k$. (Again
we are taking all closures and boundaries in~$U_\al$.)

Thus there is a neighbourhood of~$x$ disjoint from $\ol N_{k+1}$.
Since $n_\al(x) = k$, $x$ cannot lie in
$\bigcup_{j=1}^{k-1} \ol W_j \subseteq D_{k-1}$. Therefore,
$x \in E_k \less \bigcup_{j=1}^{k-1} \ol W_j = \intr \ol W_k$. Nor
does $x$ lie in any $\ol W_k \cap \ol W_j$ for $j < k$, since $D_j$ is
closed and $n_\al(x) = k$. Hence, in the open set
$\intr \ol W_k \less \bigcup_{j=1}^{k-1} (\ol W_k \cap \ol W_j)$
we can find a neighbourhood $V$ of~$x$ consisting of points of
multiplicity~$k$ only. Then $a_\al^{-1}(a_\al(V))$ is the union of $k$
neighbourhoods of the $k$ preimages of~$x$. But that would imply that
$x \in W_k$, contradicting $x \in C_k$.
\end{proof}

\begin{lemma}
\label{lm:ugly-points}
If $x \in C_k$, then either $x \in \intr \ol W_j$ for some
$j \in \{1,\dots,k\}$; or else there is some set of (at least two)
indices $j_1,\dots,j_l \in \{n_\al(x),\dots,k\}$ such that 
$x \in \bigcap_{r=1}^l \del W_{j_r}$ and also
$x \in \intr\bigl( \bigcup_{r=1}^l \ol W_{j_r} \bigr)$.
\end{lemma}

\begin{proof}
Observe that $C_k \subseteq \bigl( \bigcup_{j=1}^k \ol W_j \bigr)
\less \bigl( \bigcup_{j=1}^k W_j \bigr)
\subseteq \bigcup_{j=1}^k \del W_j$. Since $x \in \del W_j$ for
some~$j$, every neighbourhood of~$x$ intersects $U_\al \less W_j$.

Since $x \in E_k$, it follows from Lemma~\ref{lm:w-and-d} that $x$ has
a neighbourhood $V$ contained in $\bigcup_{j=1}^k \ol W_j$. Suppose,
then, that $x \notin \intr \ol W_j$ for any $j$. Then there is some
$j_1$, with $n_\al(x) \leq j_1 \leq k$ on account of 
Lemma~\ref{lm:uppercts}, such that 
$V \less \ol W_{j_1} \neq \emptyset$. It follows that
$V \less \ol W_{j_1} \subseteq \bigcup_{i\neq j_1} \ol W_i$.

Now $x \in \del W_{j_2}$ for some $j_2 \neq j_1$, with
$n_\al(x) \leq j_2 \leq k$, too. Moreover, if
$V \less (\ol W_{j_1} \cup \ol W_{j_2})$ is nonempty, then there
is some other $j_3$ such that $x \in \del W_{j_3}$, and so on.
Eventually we reach $j_l$ such that
$x \in \bigcap_{r=1}^l \del W_{j_r}$ and 
$V$ contained in $\bigcup_{r=1}^l \ol W_{j_r}$. 
\end{proof}

We now have a fairly detailed picture of the sets $U_\al$. There may
be two `bad' subsets, both of which are limit sets of sequences of
unbounded multiplicity. Away from the `bad' subsets, $X$ looks like a
branched manifold, with branchings at points of~$C_k$, since for
$m > k$ there may be several `sheets' of $W_m$ having boundary with
$W_k$, or more generally with $n_\al^{-1}(k)$, and this boundary has
multiplicity~$k$. Ultimately we must both deal with the `bad' subsets
and show that there is actually no branching within $E_k$.

\section{Reconstruction of a differential manifold}
\label{sec:spec-mfld}

To obtain a clearer view of the chart domain $U_\al$, we require the
unique continuation properties of (Euclidean) Dirac-type operators.
The next subsection analyzes the local properties of the operator
$\D$, and in particular shows that $\D$ is locally of Dirac type.

We continue to work within a fixed $U_\al$, taking closures and 
boundaries in the relative topology of $U_\al$.

\subsection{Local structure of the operator $\D$}
\label{ssc:Dirac-opr}

In order to analyze the operator $\D$ restricted to sections of the
bundle $S$ over $W_{jk}$ (recall that $\D$ is local, by
Corollary~\ref{cr:local-opr}), we need to identify the smooth
sections. The easiest way to accomplish this is to employ our
$C^\infty$ functional calculus.

By construction, the vector bundle $E$ is trivialized by the sections
$[\D,a_\al^j]$ over $W_{jk} \subset U_\al$. Any element of~$\CDA$
determines a local section of the bundle
$\End S|_{W_{jk}} \to W_{jk}$, which is likewise trivialized. Over
$W_{jk}$ there is an (involutive) algebra subbundle $C_{jk}$ of
$\End S|_{W_{jk}}$ such that $T \in \CDA$ if and only if
$T|_{W_{jk}} \in \Ga(W_{jk}, C_{jk})$. This bundle decomposes over
$W_{jk}$ as a Whitney sum of trivial matrix bundles; compare
Corollary~\ref{cr:full-rank}:
\begin{equation}
C_{jk} \simeq \bigoplus_{r=1}^s W_{jk} \x M_{n_r}(\C).
\label{eq:local-Cliff}
\end{equation}
Now define local sections $e_r \in \Ga(W_{jk}, \End S)$ by
$e_r(x) := 1_{n_r}$, for $r = 1,\dots,s$ and $x \in W_{jk}$, which are
the minimal central projectors in this decomposition of $C_{jk}$.

The bundle $S \to X$ is locally trivial, by the Serre--Swan theorem,
although \textit{a priori} the subsets $W_{jk}$ need not be
trivializing chart domains for~$S$. However, since $X$ is compact
each $W_{jk}$ can be covered by finitely many such chart domains. In 
any case, we may write $S|_{W_{jk}} = \bigoplus e_r S|_{W_{jk}}$.

\begin{lemma}
\label{lm:locally-smooth}
For all $k\geq 1$ and $j = 1,\dots,k$, the set $W_{jk}$ is a smooth
manifold of dimension~$p$, by Corollary~\ref{cr:little-mfld}, and
$E_\R|_{W_{jk}}$ is isomorphic to the cotangent bundle of $W_{jk}$.
Moreover, the algebra of restrictions of elements of $\A$ to $W_{jk}$
is isomorphic to a subalgebra of $C_b^\infty(a_\al(W_{jk}))$.
\end{lemma}

\begin{proof}
Define a bundle morphism
$\rho \: T^*W_{jk} \to E_\R|_{W_{jk}}$ by giving the corresponding map
$\rho_*$ on sections:
\begin{equation}
\rho_*\biggl(  \sum_{r=1}^p b_{r\al} \,da^r_\al \biggr)
:=  \sum_{r=1}^p b_{r\al} \,[\D,a^r_\al].
\label{eq:cotg-trans}
\end{equation}
This is a well-defined bundle isomorphism since the two local bases 
of sections $\{da^1_\al, \dots, da^p_\al\}$ and
$\{[\D,a^1_\al], \dots, [\D,a^p_\al]\}$ determine the trivial bundles
$T^*W_{jk}$ and $E_\R|_{W_{jk}}$, respectively.

It suffices to show that each $b \in \A$ can be written as a smooth
bounded function of the $a_\al$ over $W_{jk}$. We already know that
$b|_{W_{jk}} = f \circ a_\al$ for a unique bounded locally Lipschitz
function~$f$. From Corollary~\ref{cr:local-sum} we get
$[\D,b] = \sum_i b_i\,[\D,a^i_\al]$, where each $b_i$ is continuous on
$W_{jk}$ and $cb_i \in \A$ for $c \in \A$ compactly supported within
$W_{jk}$. If $K \subset W_{jk}$ is compact, by taking
$c = \phi_K \in \A$ with $\supp\phi_K \subset W_{jk}$ and $\phi_K = 1$
on~$K$, we get $b_i|_K = f_{i,K} \circ a_\al$ from the Lipschitz
functional calculus, and the uniqueness of these representatives shows
that they agree on overlaps, yielding a single function representation
$b_i = f_i \circ a_\al$ with $f_i$ a locally Lipschitz function on
$a_\al(W_{jk})$.

Choosing a sequence of $C^1$ functions $g_k$ converging to $f$ in the
Lipschitz norm on every compact subset of $a_\al(W_{jk})$, and using
the uniqueness of the coefficients $b_i$, we see that $f_i = \del_i f$.
Thus $\del_i f$ is bounded and locally Lipschitz for $i = 1,\dots,p$,
and so $f$ is in fact $C^1$. A straightforward induction now shows
that $f$ is actually~$\Coo$.
\end{proof}

\begin{rmk}
\label{rk:locally-smooth}
Since $\H_\infty \simeq q\A^m$ is a finitely generated projective 
module over $\A$, with $q \in M_m(\A)$, it follows that over any open
subset $U \subseteq W_{jk}$ for which $S|_U \to U$ is trivial, the 
local sections $\xi \in \Ga(U,S)$ may be regarded as column vectors 
with entries $f_i \circ a_\al|_U$ for $i = 1,\dots,N$, where  
$f_i \in C_b^\infty(a_\al(U))$. In other words, the sections in 
$\H_\infty$ have smooth coefficient functions over $W_{jk}$.
\end{rmk}

\begin{lemma}
\label{lm:diffl-oper}
For each $k \geq 1$ and $j = 1,\dots,k$, the operator $\D$ is an
elliptic first order differential operator on $\H_\infty|_{W_{jk}}$.
\end{lemma}

\begin{proof}
By Remark~\ref{rk:locally-smooth}, we can regard $\H_\infty|_{W_{jk}}$
as a subspace $\Ga_{\infty,\al}(W_{jk},S)$ of $\Gaoo(W_{jk},S)$. In
like manner, if $\Omega^1(W_{jk})$ denotes the $\A$-module of smooth
$1$-forms on~$W_{jk}$, then
$$
\Omega^1(W_{jk}) \ox_\A \H_\infty|_{W_{jk}}
\simeq \Ga_{\infty,\al}(W_{jk}, T^*_\C W_{jk} \ox S).
$$

Choose any connection $\nabla \: \Ga_{\infty,\al}(W_{jk},S)
\to \Ga_{\infty,\al}(W_{jk}, T^*_\C W_{jk} \ox S)$, and define an
$\A$-linear map
$\hat c \: \Ga_{\infty,\al}(W_{jk}, T^*_\C W_{jk} \ox S)
\to \Ga_{\infty,\al}(W_{jk},S)$ by setting
$\hat c(da \ox \xi) := [\D,a]\,\xi$ for $a \in \A|_{W_{jk}}$ and
$\xi \in \H_\infty|_{W_{jk}}$, extending by $\A$-linearity. Then,
recalling that $\D$ is a local operator so that $\D$ maps
$\H_\infty|_{W_{jk}}$ to $\H_\infty|_{W_{jk}}$, we get
$$
(\D - \hat c \circ \nabla)(a\xi)
= [\D,a]\,\xi + a\,\D\xi - [\D,a]\,\xi - a(\hat c \circ \nabla)\xi
= a(\D - \hat c \circ \nabla)\xi,
$$
and thus $\D = \hat c \circ \nabla + B$ with
$B \in \Gaoo(W_{jk}, \End S)$. Since $\hat c \circ \nabla$ is a
first-order differential operator, so too is~$\D$. Let $\sg_\D$ denote
its principal symbol.

For ellipticity, suppose that $x \in W_{jk}$ and $v \in T^*_x W_{jk}$
are such that $\sg_\D(x,v)$ is not invertible. Since $\D$ is of
first order, it may be evaluated as $\sg_\D(x,v) = [\D,a](x)$, for
any $a = a^* \in \A$ such that $da(x) = v$. Since we assume the map
$\sg_\D(x,v) \in \End S_x$ is not invertible, there is some 
$\xi \in \H_\infty$ with $\xi(x) \neq 0$ such that
$$
0 = \sg_\D(x,v)\xi(x) = \sum_{j=1}^p a_j(x)\,[\D,a^j_\al](x) \xi(x)
= \sum_{j=1}^p a_j(x) \xi^j(x),
$$
where we have set
$\xi^j(x) := [\D,a^j_\al](x)\xi(x) := \d a^j_\al(x)\xi(x)$. Now if
any $\d a^j_\al(x)$ had a zero eigenvector in $S_x$, then the
representative $\Ga'_x$ in $\End S_x$ of
$\d a^1_\al(x) \wyw \d a^p_\al(x) \in \La^p E_x$ would have a zero
eigenvector. (By Corollary~\ref{cr:local-sum}, this multivector spans
$\La^p E_x$.) But in that case $\Ga'_x = 0$, since
$\dim \La^p E_x = 1$, contrary to the proof of
Proposition~\ref{pr:lin-indp}. Therefore,
$\d a^1_\al(x),\dots,\d a^p_\al(x)$ are invertible in $\End S_x$, as
well as linearly independent. We conclude that the vectors
$\xi^j(x) \in S_x$ are linearly independent, which forces
$a_1(x) =\cdots= a_p(x) = 0$. This in turn entails $v = da(x) = 0$,
since each $a_j(x) = \del_jf_a(a_\al(x))$ when $a = f_a(a_\al(x))$.
Hence $\D$ is elliptic.
\end{proof}

\begin{lemma}
\label{lm:central-elt}
For all $a,b \in \A|_{W_{jk}}$, the operator
$$
[[\D^2, a], b] = [\D,a]\,[\D,b] + [\D,b]\,[\D,a]
$$
is a central element of the algebra $\CDA\bigr|_{W_{jk}}$.
\end{lemma}

\begin{proof}
We need only consider the special case $a = b$, in view of the
polarization identity
$$
[\D,a]\,[\D,b] + [\D,b]\,[\D,a]
= [\D, a + b]\,[\D, a + b] - [\D,a]\,[\D,a] - [\D,b]\,[\D,b].
$$
The operator $[\D^2, [\D,a]\,[\D,a]]$ is a differential operator of at
most second order over $W_{jk}$; it is of at most first order if and
only if~\cite{BerlineGV}:
$$
\bigl[\bigl[ [\D^2, [\D,a]\,[\D,a]], b \bigr], c \bigr] = 0,
\word{for all} b,c \in \A.
$$
Again we simplify to $b = c$; using the first order condition, we
find that
$$
\bigl[\bigl[ [\D^2, [\D,a]\,[\D,a]], b \bigr], b \bigr]
= 2 \bigl[ [\D,b]\,[\D,b], [\D,a]\,[\D,a] \bigr].
$$
Thus $[\D^2,[\D,a]\,[\D,a]]$ is of first order if and only if
$[\D,a]\,[\D,a]$ commutes with all $[\D,b]\,[\D,b]$ in
$\CDA|_{W_{jk}}$. Similarly, the operator $[\D^2,[\D,a]]$ is first
order if and only if $[\D,a]$ commutes with all $[\D,b]\,[\D,b]$ in
$\CDA|_{W_{jk}}$. To show that $[\D,b][\D,b]$ is central, it therefore
suffices to show that for all $a \in \A$ the operators
$[\D^2,[\D,a]\,[\D,a]]$ and $[\D^2,[\D,a]]$ are of first order at most.

For $T = [\D,a]\,[\D,a]$ or $[\D,a]$, the operator
$\Dreg^{-1} [\D^2, T]$, regarded now as a pseudodifferential operator
on $W_{jk}$, has order $(-1 + \mathrm{order}[\D^2,T])$. The regularity
condition entails that in both cases, the operator
$\Dreg^{-1} [\D^2, T]$ is bounded (see, e.g., \cite{CareyPRStwo}) and
thus of order at most zero. We conclude that $[\D^2,T]$ is a
differential operator of order at most one. Hence the commutators
$\bigl[ [\D,b]\,[\D,b], [\D,a] \bigr] \bigr|_{W_{jk}}$ and
$\bigl[ [\D,b]\,[\D,b], [\D,a]\,[\D,a] \bigr] \bigr|_{W_{jk}}$ vanish
for all $a,b \in \A$.
\end{proof}

\begin{lemma}
\label{lm:metric-smooth}
For all $a,b \in \A$ with compact support in $W_{jk}$, the densely
defined operator $\bigl[\D, [\D,a]\,[\D,b] + [\D,b]\,[\D,a] \bigr]$
extends to a bounded operator on~$\H$.
\end{lemma}

\begin{proof}
The operator in question maps $\H_\infty$ to itself, so by
Proposition~\ref{pr:lin-bdd} we need only show that it is $\A$-linear
on this domain. This follows from the first order condition together
with the centrality of $[\D,a]\,[\D,b] + [\D,b]\,[\D,a]$ in
$\CDA|_{W_{jk}}$; for if $c \in \A$, then
$$
\bigl[ \bigl[ \D, [\D,a]\,[\D,b] + [\D,b]\,[\D,a] \bigr], c \bigr]
= \bigl[ [\D,c], [\D,a]\,[\D,b] + [\D,b]\,[\D,a] \bigr] = 0,
$$
since we may assume without loss of generality that
$\supp c \subset W_{jk}$ also.
\end{proof}

\begin{rmk}
Though we do not need it here, it is interesting and perhaps useful
that for all $a,b \in \A$, the operator $[[\D^2,a], [\D,b]]$ is
bounded over $W_{jk}$. This is proved just as in the previous Lemma.
When we later globalize these results, we shall likewise see that this
operator is globally bounded.
\end{rmk}

\begin{prop}
\label{pr:cliff-dirac}
The space $\H_\infty|_{W_{jk}}$ carries a nondegenerate representation
of the algebra $\Gaoo(W_{jk},C)$ of smooth sections of an algebra
bundle $C = \bigoplus_{r=1}^s \Cliff(T^*W_{jk},g_r)$, which is the
Whitney sum of complex Clifford-algebra bundles defined by finitely
many Euclidean metrics $g_r$ on~$T^*W_{jk}$.
\end{prop}

\begin{proof}
The bundle isomorphism $\rho\: T^*W_{jk} \to E_\R|_{W_{jk}}$ 
determined by~\eqref{eq:cotg-trans} extends to an isomorphism
$\La^\8\rho \: \La^\8 T^*W_{jk} \to \La^\8 E_\R|_{W_{jk}}$. The map
\eqref{eq:cotg-trans} and the action of $E|_{W_{jk}}$ on
$\H_\infty|_{W_{jk}}$ together determine an action, also called 
$\rho_*$, of the local $1$-forms $\Omega^1(W_{jk})$ on
$\H_\infty|_{W_{jk}}$. If $\eta = \sum_{l} b_{l\al} \,da^l_\al$ and
$\zeta = \sum_{m} c_{m\al} \,da^m_\al$ are two such $1$-forms, then
$$
\rho_*(\eta)\,\rho_*(\zeta) + \rho_*(\zeta)\,\rho_*(\eta)
= \sum_{l,m}  b_{l\al} c_{m\al} \,
\bigl( [\D,a^l_\al]\,[\D,a^m_\al] + [\D,a^m_\al]\,[\D,a^l_\al] \bigr)
$$
is central in~$\CDA|_{W_{jk}}$, and has bounded commutator with $\D$
over $W_{jk}$.

As was noted after \eqref{eq:local-Cliff}, the algebra
$\CDA|_{W_{jk}}$ has minimal central projectors $e_1,\dots,e_s$; for 
all $a,b \in \A$ and each $r = 1,\dots,s$, the maps
$a\,db \mapsto e_r a\,[\D,b]\,e_r = a\,[\D,b]\,e_r$ make sense over
$W_{jk}$.

Since $[\D,a]\,[\D,b] + [\D,b]\,[\D,a]$ is central in
$\CDA|_{W_{jk}}$, over $W_{jk}$ it decomposes as a sum of scalar
matrices:
\begin{equation}
[\D,a]\,[\D,b] + [\D,b]\,[\D,a]
=: \bigoplus_{r=1}^s -2g_{r,\al}(da,db) \,1_{n_r},
\label{eq:multi-Cliff}
\end{equation}
where each $g_{r,\al}$ is a symmetric bilinear form on
$\Omega^1(W_{jk})$ with values in $\A|_{W_{jk}}$. Note that, for $a$
real, the $g_{r,\al}(da,da)$ are nonnegative since the operator
$[\D,a]^* [\D,a]$ is positive; and at each $x \in W_{jk}$ the matrix
with entries $g_{r,\al}(da_\al^j,da_\al^k)(x)$ is positive definite:
compare the discussion after \eqref{eq:metr-defn}. Thus, each
$g_{r,\al}$ is a positive definite Euclidean metric on $1$-forms. For
each~$r$, the map $\rho_*$ defines an action of the Clifford-algebra
bundle $\Cliff(T^*W_{jk}, g_{r,\al})$.
\end{proof}

\begin{corl}
\label{cr:D-is-Dirac}
Over the set $W_{jk}$, the operator $\D$ is, up to the addition of an
endomorphism of $S|_{W_{jk}}$, a direct sum of Dirac-type operators
with respect to the several metrics $g_1,\dots,g_s$.
\end{corl}

\begin{proof}
The symbol of $\D$ at $(x,v) \in T^*W_{jk}$ is given by $[\D,f](x)$
for any smooth function $f$ such that $df(x) = v$. Also,
\eqref{eq:multi-Cliff} says that $[\D,f]$ is Clifford multiplication
by $df$ in each summand of~$C_{jk}$.
\end{proof}

\subsection{Injectivity of the local coordinates}
\label{ssc:ucp-good-charts}

In this section we shall show that the sets $C_k$, $B_k$ for
$k\geq 1$, and $n_\al^{-1}(\infty)$ are all empty. For that, we shall
use the weak unique continuation property of Dirac-type operators
\cite{BoossMW}, as well as the strong unique continuation
property~\cite{Kim}. Both unique continuation properties have long
histories, for which we refer the reader to the papers cited. The
precise statements we require are as follows.

\begin{thm}[Weak Unique Continuation Property
{\cite[Thm.~2.1]{BoossMW}}]
\label{thm:wucp}
Let $\D$ be an operator of Dirac type acting on sections of a vector
bundle $V$ over a smooth manifold $M$ (that need not be compact), and
suppose that $\D\xi = 0$ for some $\xi \in \Ga(M,V)$. If $\xi$ is zero
on an open set $U \subset M$, then $\xi$ vanishes on the whole
connected component containing the open set~$U$.
\end{thm}

\begin{rmk}
Actually, in \cite{BoossMW} the weak unique continuation property is
shown to be stable under a large family of possibly nonlinear
perturbations of order zero. This allows us to employ unique
continuation for nonzero eigenvalues $\la$, on replacing $\D$ by
$\D - \la$.
\end{rmk}

\begin{thm}[Strong Unique Continuation Property {\cite[Cor.~2]{Kim}}]
\label{thm:sucp}
Let $U$ be a connected open subset of $\R^p$, with $p \geq 3$, and let
$D = \sum_{j=1}^p \ga_j \,\del/\del x^j$ be the constant-coefficient
Dirac operator on~$U$. Let $V \in L^r(U, M_m(\C))$ where
$r = (3p - 2)/2$, and suppose that $(D + V)\xi = 0$ for some
$\xi \in L^2(U,\C^m)$ with $D\xi \in L^2(U,\C^m)$. If
$\int_{|x-x_0|<\eps} |\xi(x)|^2 \,d^px = O(\eps^N)$ for all~$N$, for
some $x_0 \in U$, then $\xi$ is identically zero on~$U$.
\end{thm}

\begin{rmk}
The case $p = 2$ was proved by Carleman~\cite{Carleman} for
$V \in L^\infty(U)$. The corresponding result for $p = 1$ and
$V \in L^1(U)$ can be proved directly from the explicit solution for
$(\D + V)\xi = 0$. The various integrability constraints on $V$ will
not concern us, as we will be working with $V \in L^\infty(U)$, where
$U$ is a bounded open set in~$\R^p$.
\end{rmk}

A basic result we require is the smoothness of eigenspinors of the
Euclidean Dirac operator, a consequence of its ellipticity only: see,
for example, \cite[Thm.~III.5.4]{LawsonM}.

\begin{lemma}
\label{lm:smooth-eigs}
Let $U \subset \R^p$ be open and bounded, let $\Shat \to U$ be its
spinor bundle, and let $\Dhat = \sum c(dx^j)\,\del_j + \widehat V$
denote the constant coefficient Dirac operator on $U$ perturbed by
$\widehat V \in L^\infty(U, \End\Shat)$. If $s \in L^2(U,\Shat)$
satisfies $\Dhat s = \la s$ with $\la \in \R$, then $s$ is a smooth
section.
\qed
\end{lemma}

\begin{lemma}
\label{lm:derivatives-agree}
Let $V,W \subset U_\al$ be open subsets such that $a_\al$ is injective
over $\ol V$ and $\ol W$ separately. If
$a \in\A$ satisfies $a|_V = f \circ a_\al$ and $a|_W = g \circ a_\al$,
for smooth functions $f$ and~$g$, then for each multi-index
$K = (k_1,\dots,k_p)$, the partial derivatives $\del^K f$ and
$\del^K g$ extend continuously to $a_\al(\ol V)$ and
$a_\al(\ol W)$ respectively. At each
$x \in \ol V \cap \ol W$, the partial derivatives agree:
$\del^Kf(a_\al(x)) = \del^Kg(a_\al(x))$ for all~$K$.
\end{lemma}

\begin{proof} 
The case $|K| = 0$ is the equality $f(a_\al(x)) = g(a_\al(x))$ found
in Proposition~\ref{pr:local-homeos}. Over $V$, $\del_kf \circ a_\al$
agrees with the coefficient function $a_k$ of the expansion
$[\D,a] = \sum_{k=1}^p a_k\,[\D,a_\al^k]$.

Since $a_k \in \A$, the proofs of Proposition~\ref{pr:local-homeos}
and Lemma~\ref{lm:locally-smooth} show that $\del_k f$ extends to a
continuous, bounded and locally Lipschitz function on
$a_\al(\ol V)$. Hence the extension of $f$ is actually
$C^1$ on $a_\al(\ol V)$. Replacing $a = f \circ a_\al$
by $a_k = \del_kf \circ a_\al$, and the $a_k$ by the coefficients
$a_{kl}$ in the expansion of $[\D, a_k]$, and so on, we see that all
the $\del^K f$ extend in like manner to $\ol V$, and
all these extensions are~$C^N$ there. Similar comments apply to the
$\del^K g$.

Now an application of Proposition~\ref{pr:local-homeos} at each stage
of the iteration yields $\del^K f(a_\al(x)) = \del^K g(a_\al(x))$ for
any $x \in \ol V \cap \ol W$.
\end{proof}

\begin{lemma}
\label{lm:sections-agree}
Let $V,W \subset U_\al$ be open subsets such that $a_\al$ is injective
over $\ol V$ and $\ol W$ separately. Let $U \subseteq U_\al$ be an
open subset over which $S|_U$ is trivial. For any $\xi \in \H_\infty$,
we can simultaneously write
$$
\xi\bigr|_{\ol V \cap U} = (f_1 \circ a_\al, \dots, f_N \circ a_\al),
\qquad
\xi\bigr|_{\ol W \cap U} = (g_1 \circ a_\al, \dots, g_N \circ a_\al),
$$
for $f_j,g_j \in \Coo_b(a_\al(U))$, $j = 1,\dots,N$. Then, at any
$x \in \ol V \cap \ol W \cap U$, these component functions (and their 
derivatives) agree:
\begin{equation}
\del^K f_j(a_\al(x)) = \del^K g_j(a_\al(x)),
\word{for all $K$ and}  j = 1,\dots,N.
\label{eq:components-agree}
\end{equation}
\end{lemma}

\begin{proof}
Since $S|_U$ is trivial, we may identify elements of
$\H_\infty|_U \subseteq \Gaoo(U,S)$ with column vectors in $\A^N|_U$.
The restrictions of their component functions to $\ol V \cap U$ and
$\ol W \cap U$ respectively can in turn be represented as
$f_j \circ a_\al$, respectively $g_j \circ a_\al$, using the Lipschitz
functional calculus. The uniqueness of the locally Lipschitz functions
thereby obtained, together with Lemma~\ref{lm:derivatives-agree},
yields~\eqref{eq:components-agree}.
\end{proof}

We now begin the process of tidying up $U_\al$ by showing that the
subsets $C_k$ are empty, in the two Propositions which follow. The
techniques and notation in the next proof will be reused when we
banish the remaining undesirable subsets.

\begin{prop}
\label{pr:reg-open}
For each $\al$ and all $m \geq 1$, the equality $\intr\ol W_m = W_m$
holds.
\end{prop}

\begin{proof}
Suppose $W_m$ is not empty, and that $\intr\ol W_m \neq W_m$. Then
all $x \in \del W_m$ which lie in the interior of $\ol W_m$ have
multiplicity $n_\al(x) < m$ by Lemma~\ref{lm:ugly-sets}.

In particular, the result is proved for $m = 1$, and we can now assume
that $m \geq 2$. Thus there are at least two `sheets' $W_{mi}$ and
$W_{mj}$, for some $i,j \in \{1,\dots,m\}$, for which some
$x \in \ol W_m \less W_m$ is a common boundary point within~$U_\al$.

Fix such an $x \in \ol W_m \less W_m$ and choose an open 
neighbourhood $U$ of~$x$ such that $S|_U$ is trivial. Then we may, 
and shall, regard sections of $S|_U$ as functions $\xi\: U \to \R^N$.

We choose, once and for all, a numbering of the `sheets' $W_{mj}$ of
$W_m$. For $i,j = 1,\dots,m$, define switching maps
$\vf_{ij} \: W_{mj} \to W_{mi}$ by
$\vf_{ij}(x) := a_\al^{-1}(a_\al(x)) \cap W_{mi}$; these maps are
homeomorphisms among these $m$ subsets of $W_m$, and they permute
the $m$-element sets $a_\al^{-1}(a_\al(x))$ for each~$x$. Write
$a_{\al,j}^{-1}$ for the homeomorphism $a_\al(W_m) \to W_{mj}$ inverse
to $a_\al$, and note that
$$
\vf_{ij} \circ a_{\al,j}^{-1} = a_{\al,i}^{-1},
\word{for}  i,j = 1,\dots,m.
$$

Now consider the bundle $\Shat \to a_\al(\ol W_m)$ of rank $N$,
obtained by pulling back $S \to \ol W_m$ via the map
$a_{\al,1}^{-1} \: a_\al(\ol W_m) \to \ol W_{m1}$. Here we are using 
Proposition~\ref{pr:local-homeos}. Observe that since 
$S|_{U \cap \ol W_{m1}}$ is trivial, so too is
$\Shat|_{a_\al(U \cap \ol W_{m1})}$. Thereby, sections of 
$\Shat|_{a_\al(U \cap \ol W_{m1})}$ are $N$-tuples of functions over
$a_{\al,1}(U) = a_\al(U) = a_{\al,j}(U)$.

Given a section $\xi \in \Ga(S|_{U \cap \ol W_m})$, the two pullbacks
$\xi \circ a_{\al,i}^{-1}$ and $\xi \circ a_{\al,j}^{-1}$ are 
different sections of $\Shat$ corresponding to the restrictions 
of~$\xi$ to $\ol W_{mi}$ and $\ol W_{mj}$ respectively. By
Lemma~\ref{lm:sections-agree}, these two pullback sections agree at
points of $\ol W_m \less W_m$.

We define an operator $\Dhat$ on smooth sections of 
$\Shat|_{a_\al(U \cap \ol W_m)}$ by
$$
\Dhat\,\hat\xi := \bigl(\D(\hat\xi \circ a_\al|_{\ol W_{m1}})\bigr)
\circ a_{\al,1}^{-1},  \word{for}
\hat\xi \in \Gaoo(a_\al(\ol W_m), \Shat).
$$
The operator $\Dhat$ is well defined, and gives a first order
differential operator on $\Gaoo(a_\al(\ol W_m), \Shat)$. Indeed, by
Lemma~\ref{lm:diffl-oper},
\begin{align*}
\Dhat\,\hat\xi
&= \biggl( \sum_{j=1}^p [\D,a^j_\al]\,(\del_j \hat\xi) \circ a_\al
+ V(\hat\xi \circ a_\al) \biggr) \circ a_{\al,1}^{-1}
\\
&= \sum_{j=1}^p (a_{\al,1}^{-1})^* c(da^j_\al)\,\del_j\hat\xi
+ ((a_{\al,1}^{-1})^*V)\,\hat\xi,
\\
&= \sum_{j=1}^p c(dx^j) \,\del_j\hat\xi
+ ((a_{\al,1}^{-1})^*V)\,\hat\xi,
\end{align*}
where $c$ denotes Clifford multiplication coming from
Proposition~\ref{pr:cliff-dirac} and $V \in \Ga(\ol W_m, \End S)$ is
smooth, because $V$ maps $\H_\infty$ to $\H_\infty$ and is uniformly
bounded. That last statement holds because $V$ is $\A$-linear, and the
proof of Proposition~\ref{pr:lin-bdd} shows that the norm of $V$ over
$\ol W_m$ is determined on a finite generating set of $\H_\infty$.
Thus $\Dhat$ is a bounded perturbation of (a possible direct sum of
copies of) a constant coefficient Dirac operator on an open subset
of~$\R^p$.

Take two sheets $W_{mi}$, $W_{mj}$ as above, with 
$x \in \ol W_{mi} \cap \ol W_{mj}$. Observe that $a_\al(W_m)$ is 
open, by Lemma~\ref{lm:pancakes} and Corollary~\ref{cr:locally-open},
and that Proposition~\ref{pr:local-homeos} shows that
$a_\al(\intr\ol W_m)$ is also open. Similarly, the set
$a_\al(U \cap \intr\ol W_m)$ is open.

Next, let $\xi \in \H_\infty$ be any eigenvector of $\D$, with
eigenvalue $\la$ say, and define a section over the open set
$a_\al(U \cap \intr\ol W_m)$ by
$$
\psi_{ij\la}(t) := \begin{cases}
0, &\text{if } t \in a_\al(U \cap \intr\ol W_m \less W_m), \\
\xi \circ a_{\al,i}^{-1}(t) - \xi \circ a_{\al,j}^{-1}(t),
&\text{if } t \in a_\al(U \cap W_m). \end{cases}
$$
This section is well defined because at any boundary point
$x \in \del\ol W_{mi} \cap \del\ol W_{mj}$ the two local sections
$\xi \circ a_{\al,i}^{-1}$ and $\xi \circ a_{\al,j}^{-1}$ agree, by
Lemma \ref{lm:sections-agree}. Since $\D\xi = \la\xi$, it is clear
that $\Dhat \psi_{ij\la} = \la \psi_{ij\la}$ on the given domain. Thus
by Lemma~\ref{lm:smooth-eigs}, $\psi_{ij\la}$ is a smooth section of
$\Shat$ over $a_\al(U \cap \intr\ol W_m)$.

Consider the Taylor expansion of the smooth function
$h_{ij\la} := \pairing{\psi_{ij\la}}{\psi_{ij\la}}$ about some point
$t = a_\al(x) \in a_\al(U \cap \del\ol W_{mi} \cap \del\ol W_{mj})$.
According to Lemma~\ref{lm:derivatives-agree}, it satisfies
$$
h_{ij\la}(s)
= \sum_{|J|=0}^N \frac{\del^J h_{ij\la}(t)}{J!} \,(s - t)^J
+ R_{N+1}(t,s) = R_{N+1}(t,s),
$$
since the leading terms vanish up to the prescribed order~$N$. Since
$h_{ij\la}$ can be taken to have bounded derivatives of order $N+1$,
the remainder can be estimated by 
$|R_{N+1}(t,s)| \leq C\,|s - t|^{N+1}$. Thus, on a ball of radius
$\eps$ about $t = a_\al(x)$, the function
$h_{ij\la} = \pairing{\psi_{ij\la}}{\psi_{ij\la}}$ is bounded by
$C\eps^{N+1}$ for any~$N$. Hence for each such $\eps$, we get an
estimate
$$
\int_{B_\eps(a_\al(x))} \pairing{\psi_{ij\la}}{\psi_{ij\la}} \,d^p t
= O(\eps^k)
$$
for every $k \in \N$. Since we know that
$(\Dhat - \la)\psi_{ij\la} = 0$, and that $\Dhat - \la$ is a bounded
perturbation of the constant-coefficient Dirac operator, we may apply
the main result of \cite{Kim}, cited here as Theorem~\ref{thm:sucp},
to deduce that $\psi_{ij\la}$ vanishes on a neighbourhood of
$a_\al(x)$ in $a_\al(U \cap \intr\ol W_m)$.

Hence the restrictions of each eigenvector $\xi$ of $\D$ to the two
sheets $W_{mi}$ and $W_{mj}$ yield equal sections of~$\Shat$, namely
$\xi \circ a_{\al,i}^{-1} = (f_1,\dots,f_N)$ and
$\xi \circ a_{\al,j}^{-1} = (g_1,\dots,g_N)$, over the open set
$a_\al(U \cap \intr\ol W_m)$. Thus each $f_r = g_r$ on this set, and
so the sections
$$
\xi\bigr|_{W_{mi}} = (f_1,\dots,f_N) \circ a_\al\bigr|_{W_{mi}},
\qquad
\xi\bigr|_{W_{mj}} = (g_1,\dots,g_N) \circ a_\al\bigr|_{W_{mj}}
$$
satisfy $\xi|_{W_{mi}} = \xi|_{W_{mj}} \circ \vf_{ij}$ over
$a_\al(U \cap \intr\ol W_m)$. Hence, for any two eigenvectors
$\xi,\eta$ of $\D$, the function $\pairing{\xi}{\eta}$ is constant on
the fibres $a_\al^{-1}(t) \cap U \cap (W_{mi} \cup W_{mj})$ for all
$t \in a_\al(U \cap W_m)$.

Since the eigenvectors of $\D$ span $\H$ and are contained in
$\H_\infty$, and since $\H_\infty$ is a full right $\A$-module, the
algebra $\A$ is densely generated by functions of the form
$\pairing{\xi}{\eta}$. Hence, no function in $\A$ can distinguish
points in $a_\al^{-1}(t) \cap U \cap (W_{mi} \cup W_{mj})$. Since
$m \geq 2$ and we have assumed that $W_m$ is nonempty, we have 
reached a contradiction.
\end{proof}

Next we tackle the other possible ``branch points'' in $E_k$. We shall
assume, without loss of generality, that $W_m \neq \emptyset$ for each
$m \geq 1$. If some $W_k$ is actually empty, we may omit it and
renumber these subsets without affecting the arguments below.

\begin{prop}
\label{lm:Wclosures-disjoint}
For each $\al$ and all $k,l = 1,\dots,m$ with $k \neq l$, the sets
$\ol W_k \cap E_m$ and $\ol W_l \cap E_m$ are disjoint. Hence $E_m$ is
the disjoint union of its subsets $W_j$, that is,
$$
E_m = W_1 \uplus\cdots\uplus W_m .
$$
\end{prop}

\begin{proof}
We proceed by induction. First of all,
$E_1 = \intr(n_\al^{-1}(1)) = W_1$. Suppose then that
$E_k = \biguplus_{j=1}^k W_j$ for $k < m$. Using
Proposition~\ref{pr:reg-open} and Lemma~\ref{lm:ugly-points}, we 
obtain
$$
C_m = \bigcup_{k=1}^{m-1} \del W_k \cap \del W_m \cap E_m,
$$
and $\bigcup_{k=1}^{m-1} \del W_k \cap \del W_m \cap E_m \subseteq
\intr\bigl( \bigcup_{k=1}^{m-1} \ol W_k \cap \ol W_m \cap E_m \bigr)$.
It follows that
$$
\del\ol W_m \cap E_m = \del W_m \cap E_m
\subseteq \bigcup_{j=1}^{m-1}(\del W_j \cap E_m).
$$
(Notice that $\del\ol W_m = \ol W_m \less \intr\ol W_m
= \ol W_m \less W_m = \del W_m$ on account of
Proposition~\ref{pr:reg-open}.)

Since $E_{m-1} = \biguplus_{j=1}^{m-1} W_j$ by the inductive
hypothesis, we obtain
$E_{m-1} \cap \ol W_j \cap \ol W_k = \emptyset$ for 
$j \neq k \in \{1,\dots,m-1\}$. Thus,
$$
\del E_{m-1} \cap E_m \subseteq \del W_m \cap E_m
\subseteq \bigcup_{j=1}^{m-1} \del W_j \cap E_m,
$$

We require an open image of an open neighbourhood of $x \in C_m$. 
Such an $x$ lies in $\intr\bigl( \bigcup_{r=1}^l \ol W_{j_r} \bigr)$ 
for some indices $j_1,\dots,j_l$. We shall show that
$a_\al\bigl( \intr\bigl( \bigcup_{r=1}^l (\ol W_{j_r} \cup \ol W_m)
\cap E_m \bigr) \bigr)$ is open in~$\R^p$. Since we deal with finite
unions, we may suppose here that actually
$\bigcup_{r=1}^l \ol W_{j_r} = \ol W_k$ for a single~$k$. The reader
may replace $\ol W_k$, $\ol W_{kl}$ below by
$\bigcup_{r=1}^l \ol W_{j_r}$ and $\bigcup_{r=1}^l \ol W_{j_r,i_r}$
and check that the result still holds for these more general unions.

Within $\intr((\ol W_k \cup \ol W_m) \cap E_m)$, choose a sheet
$W_{kl}$ and a sheet $W_{mn}$. Then 
$a_\al \: W_{kl}\cup W_{mn} \to a_\al(W_{kl}\cup W_{mn}) \subset \R^p$
is a homeomorphism onto its image for the topology of~$\R^p$. This 
homeomorphism extends to the closure in~$U_\al$, and \textit{a 
fortiori} to the closure in~$E_m$. Thus 
$a_\al \: E_m \cap (\ol W_{kl} \cup \ol W_{mn})
\to a_\al(E_m \cap (\ol W_{kl} \cup \ol W_{mn}))$ is a homeomorphism.

Next choose any other pair of sheets $W_{kl'}$ and $W_{mn'}$ within
$\intr((\ol W_k \cup \ol W_m) \cap E_m)$. Then, taking closures in 
$E_m$, we get
$$
a_\al(\ol W_{kl} \cup \ol W_{mn}) = \ol{a_\al(W_{kl} \cup W_{mn})}
= \ol{a_\al(W_{kl'} \cup W_{mn'})}
= a_\al(\ol W_{kl'} \cup \ol W_{mn'}),
$$
on account of Lemma~\ref{lm:pancakes}. In fine, the image
$a_\al((\ol W_{kl} \cup \ol W_{mn}) \cap E_m)$ does not depend on 
the choices of sheets in $W_k$ and~$W_m$. 

Consequently, the set
$a_\al\bigl( \intr(\ol W_{kl} \cup \ol W_{mn}) \cap E_m \bigr)$ is
open and is independent of the choice of sheets. Since each pair of
sheets yields the same image, we see that
$a_\al\bigl( \intr(\ol W_k \cup \ol W_m) \cap E_m \bigr)$ is open.

Now let $x \in \del W_k \cap \del W_m \cap E_m$. Choose a 
neighbourhood $U$ of~$x$ with
$U \subseteq \intr(\ol W_k \cup \ol W_m) \cap E_m$ such that $S|_U$
is trivial. As in Proposition~\ref{pr:reg-open}, we choose sheets
$W_{mi}$, $W_{mj}$ of~$W_m$ such that 
$x \in \del W_k \cap \del W_{mi} \cap \del W_{mj}$. This is possible
since $m > 1$. [In the case of a more general union 
$\bigcup_{r=1}^l \ol W_{j_r}$, each sheet $W_{mi}$ of~$W_m$ with
$x \in \del W_{mi} \cap \bigcap_{r=1}^l \del W_{j_r}$ has
$x \in \del W_{mi} \cap \del W_{j_r,i_r}$ for all sheets 
$W_{j_r,i_r}$ of $W_{j_r}$ with $x \in \del W_{j_r,i_r}$. Hence we 
can find two sheets of~$W_m$ meeting all such sheets of $W_{j_r}$.
This happens because $n_\al(x) \leq j_1 <\cdots< j_r < m$.]

We pull back the bundle $S|_U$ to a bundle $\Shat$ over
$a_\al\bigl( \intr(\ol W_k \cup \ol W_m) \cap E_m \cap U \bigr)$.
Let $\xi \in \H_\infty$ any eigenvector of $\D$ of eigenvalue $\la$
and define a section of $\Shat$ by the formula
$$
\psi_{ij\la}(t) := \begin{cases}
0, &\text{if } t \in a_\al(\ol W_k \cap E_m \cap U),
\\
\xi \circ a_{\al,i}^{-1}(t) - \xi \circ a_{\al,j}^{-1}(t),
&\text{if } t \in a_\al(W_m \cap U). \end{cases}
$$
Just as in Proposition~\ref{pr:reg-open}, this section
is well defined since at any point of
$(\ol W_k \cup \ol W_m) \cap E_m \cap U$, the local sections
$\xi \circ a_{\al,i}^{-1}$ and $\xi \circ a_{\al,j}^{-1}$ agree, and
$(\Dhat - \la)\psi_{ij\la} = 0$.

Now $\psi_{ij\la}$ vanishes on the set $a_\al(\ol W_k\cap E_m \cap U)$
which has nonempty interior, so the weak unique continuation property
for Dirac operators, Theorem~\ref{thm:wucp}, says that either
$\psi_{ij\la}$ is identically zero on
$a_\al\bigl( \intr(\ol W_k \cup \ol W_m) \cap E_m \cap U \bigr)$, or
$a_\al(\ol W_k \cap E_m \cap U)$ is disconnected from 
$a_\al(W_m \cap U)$.

If there is at least one eigenvector $\xi$ such that this 
$\psi_{ij\la}$ is not identically zero, then there is no boundary, 
and we are done for this possible intersection (i.e.,
$\ol W_k \cap E_m = W_k \cap E_m$). Otherwise, for all eigenvectors
$\xi$, the corresponding $\psi_{ij\la}$ vanishes identically.
Therefore, for all such $\xi$, the equality
$\xi|_{W_{mi}} = \xi|_{W_{mj}} \circ \vf_{ij}$ holds and moreover,
given two eigenvectors $\xi,\eta$, the function $\pairing{\xi}{\eta}$ 
cannot distinguish points of the fibre $a_\al^{-1}(t)$ for any
$t \in a_\al(W_m \cap U)$.

Again, since $\A$ is generated by such functions and they separate 
the points of~$X$, we have reached a contradiction.

Hence, $\del W_k \cap \del W_m \cap E_m = \emptyset$. Repeating the
argument for any other $\del W_j \cap \del W_m \cap E_m$, or more
generally for $\bigcap_{r=1}^l \del W_{j_r} \cap \del W_m \cap E_m$,
completes the inductive step. The conclusion follows.
\end{proof}

At this point, if the multiplicity is bounded, $U_\al$ is the 
disjoint union of finitely many~$W_k$.

\begin{prop}
\label{pr:not-so-bad}
The equality $n_\al^{-1}(\infty) \cup B_{(\al)} = B_1$ holds.
\end{prop}

\begin{proof}
Let $x \in n_\al^{-1}(\infty) \cup B_{(\al)}$. Then there is a 
sequence $\{x_m\}_{m\geq 1} \subset U_\al$ such that $x_m \to x$ and
$n_\al(x_m) = m \to \infty$. Here we suppose, without loss of 
generality, that each $W_m$ is nonempty.

Fix, for now, any compact set $K \subset U_\al$ such that $x \in K$
and $x_m \in K$ for each~$m$. Put
\begin{align*}
\eps_m &:= |a_\al(x_m) - a_\al(x)| \to 0,
\\
e_m
&:= \max\set{d(x_m,x''_m) : x''_m \in a_\al^{-1}(a_\al(x_m)) \cap K}.
\end{align*}
Choose $x'_m \in a_\al^{-1}(a_\al(x_m)) \cap K$ such that 
$d(x_m,x'_m) = e_m$. Then
$$
e_m = d(x_m,x'_m) \leq d(x_m,x) + d(x'_m,x)
\leq C_{Y_m} \eps_m + C_{Y'_m} \eps_m,
$$
where the last estimate comes from Corollary~\ref{cr:twopt-dists}
---note that (eventually) $a_\al(x_m) \neq a_\al(x)$ since
$n_\al(x_m) \neq n_\al(x)$--- with $Y_m = \{x_m,x\}$ and
$Y'_m = \{x'_m,x\}$. Then, after perhaps passing to a subsequence, we 
are faced with two possibilities.
\begin{enumerate}
\item
There is some $\delta > 0$ such that $e_m \geq \delta$ for all~$m$.
Then, since $e_m/\eps_m \leq C_{Y_m} + C_{Y'_m}$, we see that
$C_{Y_m} + C_{Y'_m} \to \infty$. However, this cannot happen since
$C_{Y_m} + C_{Y'_m}$ is bounded by $2C_K$.
\item
Otherwise, $e_m \to 0$. In this case, for each $m$ we choose 
any $x''_m \in a_\al^{-1}(a_\al(x_m)) \cap K$. Then the sequence
$\{x''_m\}$ converges to~$x$ since
$$
d(x''_m,x) \leq d(x''_m,x_m) + d(x_m,x) \leq e_m + C_{Y_m} \eps_m
\leq e_m + C_K \eps_m \to 0.
$$
Thus, any preimage of $\{a_\al(x_m)\}_{m\geq 1}$ in~$K$ converges 
to~$x$.
\end{enumerate}

Continuing with the second case, observe now that
$\{x_m\}_{m\geq 1}$ determines a subset $U_{I,\al} \subset U_\al$ as 
follows. First fix a numbering $W_{k1},\dots,W_{kk}$ of the 
$k$~sheets of each~$W_k$. For each $k \geq 1$, choose
$i(k) \in \{1,\dots,k\}$, write $I := \{i(k)\}_{k\geq 1}$ and put
$$
U_{I,\al} := \biguplus_{k\geq 1} W_{k,i(k)}.
$$
The set $U_{I,\al}$ determined by $\{x_m\}_{m\geq 1}$ is given by 
taking, for each~$m$, the sheet of $W_m$ in which the point $x_m$~lies.

Then $a_\al \: U_{I,\al} \to \R^p$ is one-to-one and open and is a
homeomorphism onto its image. This image is the same for each
$U_{I,\al}$, by Lemma~\ref{lm:pancakes}.

Each $y \in n_\al^{-1}(\infty) \cup B_{(\al)}$ is the limit of a 
sequence of unbounded multiplicity, and thus is contained in the 
boundary of some $U_{J,\al}$. In particular, let
$y \in a_\al^{-1}(a_\al(x))$, and suppose that $y \in \del U_{J,\al}$.
Then the sequence $\{y_m\}_{m\geq 1} := a_\al^{-1}\bigl( a_\al\bigl(
\{x_m\}_{m\geq 1} \bigr) \bigr) \cap U_{J,\al}$ converges to~$y$.

Define a new compact subset $K'$ of~$U_\al$ by 
$K' := K \cup \{y\} \cup \{y_m\}_{m\geq 1}$. On running our initial 
argument again, and observing that $\{y_m\}_{m\geq 1}$ is a preimage 
of $\{a_\al(x_m)\}_{m\geq 1}$ contained in~$K'$, we deduce that 
$y_m \to x$. Therefore, $y = x$ since $K'$ is a Hausdorff space.

Hence, the point $x \in n_\al^{-1}(\infty) \cup B_{(\al)}$ is such 
that $a_\al(x)$ has only one preimage in~$U_\al$, so that $x \in B_1$.
\end{proof}

It is ironic that our final task is to remove a (possible) set of
points of multiplicity one.

\begin{prop}
\label{pr:pretty-good}
The set $B_1$ is empty. Hence
$$
U_\al = \biguplus_{k=1}^\infty W_k
= \biguplus_{k=1}^\infty \biguplus_{j=1}^k W_{kj}.
$$
\end{prop}

\begin{proof}
The previous Proposition shows that
$$
U_\al = W_1 \uplus B_1 \uplus \biguplus_{k=2}^\infty W_k.
$$
The set $W_1 \uplus B_1 = D_1$ is closed and its boundary is $B_1$.
Thus elements of $B_1$ are limits of sequences of multiplicity~$1$ 
and simultaneously limits of sequences of unbounded multiplicity,
other possibilities being excluded.

Any open neighbourhood $V$ in $U_\al$ of the closed set $D_1$ will 
thus contain points of arbitrarily high multiplicity, if $B_1$ is 
nonempty. Suppose, then, that $B_1 \neq \emptyset$. 

Choose any $U_{I,\al}$ as described in the proof of 
Proposition~\ref{pr:not-so-bad}. Then $\del U_{I,\al} = B_1$ and
$a_\al \: U_{I,\al} \to \R^p$ is one-to-one, open, and a 
homeomorphism onto its image. By Proposition~\ref{pr:local-homeos}, 
$a_\al \: \ol U_{I,\al} \to \R^p$ is also a homeomorphism onto its
image.

Take $V := W_1 \cup B_1 \cup \bigcup_{k>r} W_{k,i(k)}
\subseteq \ol U_{I,\al}$ for some $r \geq 2$. Then $V$ is open in 
$\ol U_{I,\al}$, since $V = \ol U_{I,\al}
\cap \bigl( W_1 \cup B_1 \cup \bigcup_{k>r} W_k \bigr)$ and the 
complement of $W_1 \cup B_1 \cup \bigcup_{k>r} W_k$ is 
$\biguplus_{k=1}^r W_k$ which is closed in $U_\al$ (its boundary is 
empty). Actually, since $\biguplus_{k=1}^r W_k$ is also open 
in~$U_\al$, it is a union of open connected components. Thus 
$\ol U_{I,\al} = V \uplus \biguplus_{k=1}^r W_{k,i(k)}$ expresses
$\ol U_{I,\al}$ as a union of two mutually disconnected pieces. In
the relative topology of $\ol U_{I,\al}$, we get
$B_1 \subseteq \ol{V\less B_1} \subseteq \ol V = V$ (since $V$ is a
component), so that
$B_1 \subseteq \intr \ol{V\less B_1} \subseteq V$. We conclude that
$B_1 \subseteq \intr \ol U_{I,\al}$.

Now we choose $x \in B_1$ and a neighbourhood $U$ of~$x$ such that
$S|_U$ is trivial. Choose two sets of sheets $U_{I,\al}$ and
$U_{J,\al}$ with $x \in \del U_{I,\al} \cap \del U_{J,\al}$.

Observe that, by the above argument and
Corollary~\ref{cr:locally-open}, the set $a_\al(U \cap \ol U_{I,\al})$
is open. Let $\xi$ be any eigenvector of $\D$ with eigenvalue $\la$;
define a section of
$\Shat \to a_\al(U) = a_\al(U \cap \ol U_{I,\al})$ by
$$
\psi_{IJ\la}(t) 
:= \begin{cases} 0, &\text{if } t \in a_\al(W_1 \cup B_1),
\\
\xi \circ a_{\al,I}^{-1}(t) - \xi \circ a_{\al,J}^{-1}(t),
&\text{if } t \in a_\al(\ol U_{I,\al} \less (W_1 \cup B_1)),
\end{cases}
$$
where $a_{\al,I}^{-1} \: a_\al(\ol U_{I,\al}) \to \ol U_{I,\al}$ and 
similarly for $a_{\al,J}^{-1}$. As before, $\psi_{IJ\la}$ is a 
well-defined eigenvector for $\Dhat$ on $a_\al(U)$.

The weak unique continuation property for Dirac operators now says 
that any such $\psi_{IJ\la}$ is identically zero, or $W_1$ is 
disconnected from $a_\al(\ol U_{I,\al} \less W_1)$. If there is any 
$\psi_{IJ\la}$ which is not identically zero, for each 
$x \in B_1 \subseteq U$, we are done.

Otherwise every $\psi_{IJ\la}$ vanishes identically, and thus
$\xi|_{U\cap U_{I,\al}}
= \xi|_{U\cap U_{J,\al}} \circ a_{\al,J}^{-1} \circ a_{\al,I}$.
In that case, no function of the form $\pairing{\xi}{\eta}$, with
$\xi,\eta$ being eigenvectors of~$\D$, can separate points of the
fibre $a_\al^{-1}(t)$ for $t \in a_\al(U\cap U_{I,\al})$. Since
these functions generate $\A$, we have reached a contradiction.
\end{proof}

\begin{thm}
\label{thm:lip-mfld}
The space $X$ is a compact topological manifold.
\end{thm}

\begin{proof}
By Propositions 
\ref{pr:reg-open}, \ref{lm:Wclosures-disjoint}, 
\ref{pr:not-so-bad} 
and~\ref{pr:pretty-good}, we now know that
$$
X = \bigcup_{\al=1}^n U_\al
= \bigcup_{\al=1}^n \bigcup_{j,k} W_{jk,\al}.
$$
This is a weak$^*$-open cover, and since $X$ is compact, it has a
finite subcover. Since $a_\al$ is one-to-one on each $W_{jk}$,
now renamed $W_{jk,\al}$, Corollary \ref{cr:homeo-on-image} says that
$a_\al \: W_{jk,\al} \to a_\al(W_{jk}) \subset \R^p$ is a
homeomorphism, and is an open map to $\R^p$ by Corollary 
\ref{cr:locally-open}. On any overlap
$V = W_{jk,\al} \cap W_{j'k',\bt}$ of our finite subcover, there exist
locally Lipschitz functions $g^j \: a_\al(V) \to \R$ for
$j = 1,\dots,p$, such that
$a_\bt^j \bigr|_V = g^j \circ a_\al \bigr|_V$. This follows from
Lemma~\ref{lm:local-Lip} (and the selfadjointness of each
$a^j_\bt$). Thus, the transition functions
$a_\bt \circ a_\al^{-1} \: a_\al(V) \to a_\bt(V)$ are given by
$$
a_\bt(a_\al^{-1}(\xi)) = (g^1(\xi),\dots,g^p(\xi)),
\word{for all} \xi \in a_\al(V).
$$
Thus they are all locally Lipschitz, and in particular continuous,
which is what we need.
\end{proof}

\begin{rmk}
Since we now know that we can obtain suitable charts for the
coordinate functions $a_\al$, we may suppose henceforth that the
$a_\al$ are actually one-to-one on each $U_\al$. This can be achieved
by renumbering the charts, whereby each former $W_{jk,\al}$
participating in the finite open cover of~$X$ is relabelled as some
$U_\bt$. (We have not imposed any requirement that the elements
$a^j_\al$ of~$\A$ be distinct, as $\al$ varies.)
\end{rmk}

\begin{prop}
\label{pr:metr-wstar}
Let $(\A,\H,\D)$ be a spectral manifold of dimension $p$. Then the
metric and weak$^*$ topologies on $X = \spec(\A)$ agree.
\end{prop}

\begin{proof}
This follows from Theorem~\ref{thm:lip-mfld} and
Lemma~\ref{lm:local-tops} since each $x \in X$ has a neighbourhood
on which some $a_\al$ is one-to-one.
\end{proof}

\begin{prop}
\label{pr:impl-funct}
Let $(\A,\H,\D)$ be a spectral manifold of dimension~$p$, and let
$b \in \A$. Then there is a $\Coo$ function $g \: a_\al(U_\al) \to \C$
such that $g(a_\al(x)) = b(x)$ for all $x \in U_\al$.
\end{prop}

\begin{proof}
Since $a_\al \: U_\al \to \R^p$ is (now) one-to-one and open, and so a
homeomorphism onto its image, Lemma~\ref{lm:local-Lip} says that there
is a unique locally Lipschitz function $g\: a_\al(U_\al) \to \C$ such
that $b|_{U_\al} = g \circ a_\al$. Since $b \in \A$, it is also true
that
$$
\phi\,[\D,b] = \phi \sum_{j=1}^p b_j\,[\D,a^j_\al],
$$
where $\phi \in \A$ is any function with $\supp\phi \subset U_\al$,
and $\phi b_j \in \A$. Now by choosing a sequence $\{g_{(k)}\}$ of
$C^1$ functions converging in the Lipschitz norm to~$g$ on
$\supp\phi$, we find that
$$
\phi\,[\D,b] = \phi \sum_{j=1}^p (g_j \circ a_\al)\,[\D,a^j_\al],
$$
where $g_j(\xi) := \lim_k \del_j g_{(k)}(\xi)$. The linear
independence of the $[\D, a^j_\al]$ over $U_\al$ gives us the
uniqueness of the coefficients in this expansion of $[\D,b]$; hence
$\phi b_j = \phi g_j \circ a_\al$. Since this holds for any
$\phi \in \A$ supported in $U_\al$, we see that $g_j$ is a continuous
function defined on all of $a_\al(U_\al)$ that does not depend
on~$\phi$, and that $b_j\bigr|_{U_\al} = g_j \circ a_\al$. Since
$b_j \in \A$, our Lipschitz functional shows that each $g_j$ is
locally Lipschitz, and therefore that $g$ is actually $C^1$, with
$\del_j g = g_j$ on $a_\al(U_\al)$. Repeating this argument with $b$
replaced by any $b_j$ then shows that $g$ is~$C^2$, and so on. Thus
$g$ is $\Coo$ by induction.
\end{proof}

All elements of $\A$ can now be written as smooth functions of
finitely many elements $a^j_\al$, and the smooth manifold structure of
$X$ is easier to describe.

\begin{prop}
\label{pr:genrs}
Let $(\A,\H,\D)$ be a spectral manifold of dimension $p$. The unital
algebra generated by the $np$ functions $a^j_\al$ is dense in the
unital Fr\'echet algebra $\A$. In the norm topology, this algebra is
dense in $C^*$-algebra $A$.
\end{prop}

\begin{proof}
Let $\{\phi_\al\}$ be a partition of unity subordinate to the finite
open cover $\{U_\al\}$, and let $b \in \A$. For each $\al$, define
$g_\al \: a_\al(U_\al) \to \C$, as in Proposition~\ref{pr:impl-funct},
so that $\phi_\al b = g_\al(a^1_\al,\dots,a^p_\al)$ on $U_\al$. Then
\begin{equation}
b = \sum_{\al=1}^n \phi_\al b
= \sum_{\al=1}^n g_\al(a^1_\al,\dots,a^p_\al).
\label{eq:coord-funct}
\end{equation}
Since $X$ is compact and Hausdorff, and the algebra $\A$ separates
points of $X = \spec(\A)$, for any two points $x,y \in X$ we can find
$b \in \A$ such that $b(x) = 1$ and $b(y) = 0$. Since $b$ is of the
form \eqref{eq:coord-funct}, it follows that the $np$ functions
$a^j_\al$ themselves separate the points of~$X$. Since the unital
algebra they generate (over~$\C$) is closed under complex conjugation,
the Stone--Weierstrass theorem shows that it is dense in the
$C^*$-algebra $C(X) = A$.

Now consider the closure $\A_0$ in the Fr\'echet algebra $\A$ of this
subalgebra $\C[1,a^1_1,\dots,a^p_n]$. By
Proposition~\eqref{pr:mult-cinfty}, $\A_0$ contains all functions of
the form $h(a^1_1,\dots,a^p_n)$ where $h$ is~$\Coo$, and in particular
it contains all elements of the form \eqref{eq:coord-funct}. Thus
$\A_0 = \A$, as required.
\end{proof}

\begin{lemma}
\label{lm:smooth-trans}
Each mapping $a_\al \circ a_\bt^{-1}
: a_\bt(U_\al \cap U_\bt) \to a_\al(U_\al \cap U_\bt)$ is~$\Coo$.
\end{lemma}

\begin{proof}
Proposition~\ref{pr:impl-funct} shows that for each $\al$ and each
$j = 1,\dots,p$, there exists a $\Coo$ function
$f^j \: a_\bt(U_\al \cap U_\bt) \to \R$ such that
\begin{equation}
a^j_\al(x) = f^j(a^1_\bt,\dots,a^p_\bt)(x)
\word{for}  x \in U_\al \cap U_\bt.
\label{eq:trans-fn}
\end{equation}
For $\xi = a_\bt(x)$ with $x \in U_\al \cap U_\bt$, it follows that
$a_\al \circ a_\bt^{-1}(\xi) = a_\al(x) = (f^1(\xi),\dots,f^p(\xi))$,
and so $a_\al\circ a_\bt^{-1}$ is $\Coo$ on $U_\al \cap U_\bt$.
\end{proof}

\begin{thm}
\label{th:first-result}
Let $(\A,\H,\D)$ be a spectral manifold of dimension $p$. Then
$X = \spec(\A)$ is a smooth manifold and $\A = \Coo(X)$. Moreover, $X$
is orientable.
\end{thm}

\begin{proof}
Theorem \ref{thm:lip-mfld} and Lemma~\ref{lm:smooth-trans} together
establish that $X$ is a smooth manifold. If $f \in \Coo(X)$ with
$\supp f \subset U_\al$, then by definition
$g_\al = f \circ a_\al^{-1} \: a_\al(U_\al) \to \C$ is a smooth
mapping. Now $f = g_\al \circ a_\al$ lies in $\A$ by the multivariate
$\Coo$-functional calculus of Proposition~\ref{pr:mult-cinfty}. More
generally, if $f \in \Coo(X)$ we can write $f = \sum_\al f h_\al$
where $\{h_\al\}$ is a finite partition of unity in $\Coo(X)$
subordinate to $\{U_\al\}$. Therefore, $\Coo(X) \subseteq \A$.

Conversely, if $b \in \A$, then by Lemma~\ref{lm:partn-unity} we can
choose a finite partition of unity $\{\phi_\al\}$ in~$\A$, also
subordinate to $\{U_\al\}$, and by Proposition~\ref{pr:impl-funct}
we can write $b \phi_\al = g(a^1_\al,\dots,a^p_\al)$ where $g$ is a
smooth function defined on $a_\al(U_\al)$. Thus,
$b = \sum_\al b \phi_\al$ lies in $\Coo(X)$. In fine, $\A = \Coo(X)$.

The real vector bundle $E_\R \to X$ of Corollary~\ref{cr:real-bundle}
is now seen to have smooth transition functions. The local formula
\eqref{eq:cotg-trans} extends to a map on sections
\begin{equation}
\rho_*\biggl( \sum_\al \phi_\al \sum_j b_{j\al} \,da^j_\al \biggr)
:= \sum_\al \phi_\al \sum_j b_{j\al} \,[\D,a^j_\al],
\label{eq:cotg-trans-bis}
\end{equation}
that determines a bundle morphism $\rho \: T^*X \to E_\R$. Since two
local coordinate bases of $1$-forms $\{da^1_\al, \dots, da^p_\al\}$
and $\{da^1_\bt, \dots, da^p_\bt\}$ are related on $U_\al \cap U_\bt$,
according to \eqref{eq:trans-fn}, by
$$
da^j_\al = \sum_{k=1}^p \del_k f^j(a^1_\bt,\dots,a^p_\bt) \,da^k_\bt
$$
and $[\D,a^j_\al]$ is expressed in terms of the $[\D,a^k_\bt]$ with
the same coefficients, by \eqref{eq:local-expan}, it follows that the
map \eqref{eq:cotg-trans-bis} is well-defined and that $\rho$ is an
isomorphism of vector bundles over~$X$. We also get the corresponding
isomorphism $\La^\8\rho \: \La^\8 T^*X \to \La^\8 E_\R$.

Recall the skewsymmetrization $\Ga' = \sum_\al \Ga'_\al$ of~$\Ga$
introduced in~\eqref{eq:skew-Gamma}. We have established in
Section~\ref{sec:cotg-bdl} that $\Ga'$ is a nowhere vanishing section
of the complex vector bundle $\La^p E \to X$. A small but necessary
adjustment now yields a nowhere vanishing section of
$\La^p E_\R \to X$. Since each $[\D,a^j_\al]$ is skewadjoint,
we see that, with $\w$ denoting skewsymmetrization,
$$
\bigl( [\D,a^1_\al] \wyw [\D,a^p_\al] \bigr)^*
= (-1)^{p(p+1)/2} \,[\D,a^1_\al] \wyw [\D,a^p_\al].
$$
Note that the sign can also be written as $(-1)^{\piso{(p+1)/2}}$.
Now $\Ga'$ differs from $\Ga$ by adding selfadjoint junk terms, see
Lemma~\ref{lm:metr-junk}, so that $\Ga'$ is also selfadjoint. We may
therefore rewrite each coefficient $a^0_\al$ in~\eqref{eq:skew-Gamma}
as $a^0_\al = i^{\piso{(p+1)/2}} \tilde a^0_\al$, where
$\tilde{a}^0_\al$ is a selfadjoint element in $\A$ and is thus a
\textit{real} function on~$X$. Therefore, $i^{-\piso{(p+1)/2}} \Ga'$
is a nonvanishing smooth section of $\La^p E_\R \to X$, so it equals
$\La^\8\rho_*(\nu)$, where $\nu$ is a nonvanishing $p$-form on~$X$.
The de~Rham cohomology class $[\nu]$ of this volume form confers the
desired orientation on~$X$.
\end{proof}

The following Corollary is true because it holds locally, as proved
earlier in Lemmas \ref{lm:diffl-oper} and~\ref{lm:central-elt}.

\begin{corl}
\label{cr:global-diffl-oper}
The operator $\D$ is a first order elliptic differential operator
on~$\Gaoo(X,S)$. For all $a,b \in \A$, the operator
$$
[[\D^2, a], b] = [\D,a]\,[\D,b] + [\D,b]\,[\D,a]
$$
is a central element of the algebra $\CDA$.
\qed
\end{corl}

\subsection{Riemannian structure of the spectral manifold}
\label{ssc:one-metric}

Next, we show that the algebra $\CDA$, acting as operators on
$\H_\infty$, is the carrier of a Clifford action for a unique
Riemannian metric on~$X$.

Similar techniques to those used to construct $E$ allow us to show
that the representation of the algebra $\CDA$ is irreducible,
according to the following Proposition.

\begin{prop}
\label{pr:no-endo-proj}
Let $e$ be a projector in $\Ga(X,\End S)$ such that
$e\H_\infty \subseteq \H_\infty$. If $e$ commutes with $\CDA$, then
$e = 0$ or~$1$.
\end{prop}

\begin{proof}
Note first that $[\D,a] = e\,[\D,a]\,e + (1-e)\,[\D,a]\,(1-e)$ for any
$a \in \A$, since $e$ commutes with $[\D,a]$. Consider 
$$
B := e\D(1 - e) + (1 - e)\D e
$$
as a linear map of $\H_\infty$ to itself. Since $ea = ae$ for 
$a \in \A$, we find that
$$
[B,a] = e\,[\D, a]\,(1 - e) + (1 - e)\,[\D, a]\,e = 0,
$$
so that $B$ is $\A$-linear. By Proposition~\ref{pr:lin-bdd}, the
operator $B$ extends to a bounded operator on~$\H$, which is
selfadjoint. The same is also true of $e$, which is $\A$-linear too.

As a bounded selfadjoint perturbation of~$\D$, the operator $\D - B$
is also selfadjoint with $\Dom(\D - B) = \Dom\D$: see, for instance,
\cite[App.~A]{CareyP}. If $\xi \in \H_\infty$, then
$$
(\D - B)e\xi = \D(e\xi) - (1 - e)\D e\xi = e \D(e\xi)
= e \D\xi - e\D(1 - e)\xi = e(\D - B)\xi,
$$
so that $[\D - B, e] = 0$ on $\H_\infty$, or equivalently,
$[\D, e]\xi = [B, e]\xi$ for $\xi \in \H_\infty$. Thus $[\D, e]$ 
extends to the bounded operator $[B, e]$, and therefore
\cite{BratteliRoI} we get $e \in \Dom\d$ with $\d e = [\D,e] = [B,e]$.

Since $e$ commutes with $\CDA$, we may now apply
Corollary~\ref{cr:comm-formula} to obtain
$$
\Ga\,[\D,e] = \half (-1)^{p-1} \sum_{\al=1}^n a_\al^0
\sum_{j=1}^p (-1)^{j-1} \d a^1_\al \dots
(\d e\,\d a^j_\al + \d a^j_\al\,\d e) \dots \d a^p_\al.
$$
Skewsymmetrizing in $\d e,\d a^1_\al,\dots,\d a^p_\al$ gives zero, so
for all $x \in X$ we can find $c_{j\al}(x)$ with
$$
\d e(x) = \sum_{j,\al} c_{j\al}(x)\,\d a^j_\al(x).
$$
The $c_{j\al}$ define bounded functions on $X$, and since $\d e$ and
each $\d a^j_\al$ are bounded operators on the Hilbert space $\H$ of
$L^2$-sections of~$S$ with respect to the measure~$\mu_\D$, these
functions are measurable and preserve $L^2$-sections. The endomorphism
$e$ commutes with all such functions and $\CDA$, and thus $e$ commutes
with $[\D,e]$; therefore, $e\,[\D,e] = e^2\,[\D,e] = e\,[\D,e]\,e = 0$,
and similarly $(1 - e)[\D,e] = -(1 - e)[\D,(1 - e)] = 0$, since
$e(\dl e)e = 0$ for any idempotent $e$ and derivation $\dl$. Thus
$[\D,e] = 0$. By irreducibility, $e$ is a scalar, and so it equals $0$
or~$1$.
\end{proof}

\begin{corl}
\label{cr:not-reducible}
There is no proper subbundle $\tilde S \subset S$ such that
$\CDA\,\Gaoo(X,\tilde S) \subseteq \Gaoo(X,\tilde S)$. That is,
$\H_\infty = \Gaoo(X,S)$ is irreducible as a $\CDA$-module.
\end{corl}

\begin{proof}
If $S$ were reducible, we could find a projector $e \in \Ga(X,\End S)$
with $e\,\H_\infty = \Gaoo(X,\tilde S)$ and $[e,\CDA] = 0$. By the
previous Proposition, such an $e$ is either $0$ or~$1$.
\end{proof}

\begin{prop}
\label{pr:global-cliff-dirac}
The Hilbert space $\H$ carries a nondegenerate representation of the
algebra $\Gaoo(X,C)$ of smooth sections of an algebra bundle
$C = \Cliff(T^*X,g)$, which is the complex Clifford-algebra bundle
defined by a Riemannian metric $g$ on~$X$.
\end{prop}

\begin{proof}
If $\eta = \sum_{\al,j} \phi_\al b_{j\al} \,da^j_\al$ and
$\zeta = \sum_{\bt,k} \phi_\bt c_{k\bt} \,da^k_\bt$ are two
$1$-forms in $\Omega^1(X)$, then \eqref{eq:cotg-trans-bis} yields
\begin{equation}
\rho_*(\eta)\,\rho_*(\zeta) + \rho_*(\zeta)\,\rho_*(\eta)
= \sum_{\al,\bt,j,k} \phi_\al \phi_\bt b_{j\al} c_{k\bt} \,
\bigl( [\D,a^j_\al]\,[\D,a^k_\bt] + [\D,a^k_\bt]\,[\D,a^j_\al] \bigr),
\label{eq:Cliff-repn}
\end{equation}
which is central in $\CDA$ and has bounded commutator with~$\D$.

By Corollary~\ref{cr:local-sum}, the finite products of the local
sections $\d a^j_\al = [\D,a^j_\al]$, which are restrictions
to~$U_\al$ of elements of~$\CDA$, determine a trivialization over
$U_\al$ of the bundle $\End S$. Thus there are algebra subbundles
$C_\al \to U_\al$ of $\End S|_{U_\al}$ for each~$\al$, such that
$T \in \CDA$ if and only if $T|_{U_\al} \in \Ga(U_\al, C_\al)$ for
all~$\al$. Over $U_\al$, the bundle $C_\al$ decomposes as a Whitney
sum of trivial matrix bundles, compare Corollary~\ref{cr:full-rank}:
\begin{equation}
C_\al \simeq \bigoplus_{i=1}^r U_\al \x M_{k_i}(\C).
\label{eq:semi-simplicity}
\end{equation}

Just as in the proof of Proposition~\ref{pr:cliff-dirac}, the local
sections $\tilde e_i \in \Gaoo(U_\al, \End S)$ given by
$\tilde e_i(x) := 1_{k_i}$, for $i = 1,\dots,r$, are the minimal
central projectors in this decomposition. For $a,b \in \A$ compactly
supported in $U_\al$, the map
$a\,db \mapsto \tilde e_i a\,\d b\,\tilde e_i = a\,\d b\,\tilde e_i$
makes sense over~$U_\al$.

Since $\d a\,\d b + \d b\,\d a$ is central in $\CDA$, it
decomposes over each $U_\al$ as a sum of scalar matrices:
$$
\d a\,\d b + \d b\,\d a
=: \bigoplus_{i=1}^r -2g_{i\al}(da,db) \,1_{k_i},
$$
where each $g_{i\al}$ is again a positive definite symmetric bilinear
form whose values this time are bounded smooth functions on~$U_\al$,
i.e., restrictions to $U_\al$ of elements of~$\A$. On any overlap
$U_\al \cap U_\bt$, we can write
\begin{align}
\d a\,\d b + \d b\,\d a
&= \sum_{j,k}a_{j\al} b_{k\al} \,
(\d a^j_\al\,\d a^k_\al + \d a^k_\al\,\d a^j_\al)
\nonumber \\
&= -2 \bigoplus_{i=1}^r \sum_{j,k} a_{j\al} b_{k\al}\,
g_{i\al}(da^j_\al,da^k_\al) \,1_{k_i}
\nonumber \\
&= -2 \bigoplus_{i=1}^r \sum_{j,k,n,m} a_{j\al} b_{k\al}
c^j_{m,\al\bt} c^k_{n,\al\bt}\, g_{i\bt}(da^m_\bt, da^n_\bt) \,1_{k_i}
\nonumber \\
&= -2 \bigoplus_{i=1}^r \sum_{n,m} a_{m\bt} b_{n\bt}\,
g_{i\bt}(da^m_\bt, da^n_\bt) \,1_{k_i}.
\label{eq:2-tensor}
\end{align}
Here the $c^j_{m,\al\bt}$ are the transition functions for $E_\R$,
defined in \eqref{eq:dak-combo}.

Some of the scalar components of $\d a\,\d b + \d b\,\d a$ in the
block-matrix decomposition \eqref{eq:2-tensor} might coincide, even
when $a,b$ run over all elements of~$\A$. To consolidate such blocks,
we relabel the decomposition \eqref{eq:semi-simplicity} as
\begin{equation}
C_\al = C_{\al,1} \oplus\cdots\oplus C_{\al,s},
\label{eq:global-decomp}
\end{equation}
where each $C_{\al,j}$ is a Whitney sum of matrix-algebra bundles
over~$U_\al$, in which the sections $\d a\,\d b + \d b\,\d a$ have
distinct components $-2g_{i\al}(da,db)$ in general. Let $N_j$ be the
rank of $C_{\al,j}$, so that $\sum_{j=1}^s N_j = \sum_{i=1}^r k_i$. On
comparing the scalar components of $\d a\,\d b + \d b\,\d a$ on any
overlap $U_\al \cap U_\bt$, we see that the number $s$ and the ranks
$N_1,\dots,N_s$ do not depend on~$U_\al$, so the block decomposition
\eqref{eq:global-decomp} is global. Let $e_j \in \CDA$ denote the
central projector given by the identity element of~$C_{\al,j}$. The
relation \eqref{eq:2-tensor} shows that each corresponding
$g_j := g_{j\al}$ is a globally defined symmetric $2$-tensor.

Applying Corollary~\ref{cr:not-reducible}, we find that there can only
be one such $2$-tensor $g := g_1$, and only one such global block
in~\eqref{eq:global-decomp}. Thus $g$ is a positive definite 
Euclidean metric on $T^*X$, that is to say (after transposing to the 
tangent bundle $TX$, if one prefers), a Riemannian metric on~$X$.
In view of \eqref{eq:Cliff-repn}, the map $\rho$ defines an action of
$C := \Cliff(T^*X, g)$ on $\H_\infty$, whose algebra of smooth
sections is precisely~$\CDA$.

Comparing now with \eqref{eq:prng-defn}, we see that
$$
\pairing{\d a}{\d b} = C_p \tr((\d a)^*\,\d b)
= C_p  \tr((\d a)^*\,\d b) = g(da,db)\,1
$$
for selfadjoint $a,b \in \A$. Consequently we set $C_p := N^{-1}$
where $N = \rank S$.
\end{proof}

\begin{corl}
\label{cr:D-is-Dirac-again}
The operator $\D$ is, up to the addition of an endomorphism of $S$, a
Dirac-type operator with respect to the metric~$g$.
\qed
\end{corl}

\begin{rmk}
The orientation on (the cotangent bundle of) the manifold $X$ is fixed
by~$\Ga$, according to Theorem~\ref{th:first-result}. Recall that the
chirality element $\ga$ for the Clifford algebra $\Cliff(T^*X, g)$
depends on a choice of orientation of $T^*X$ \cite{Polaris,LawsonM};
reversal of this orientation the replaces $\ga$ by $-\ga$. Having
chosen the orientation of $T^*X$, let $\ga$ be the chirality element
in~$C$.
\end{rmk}

\begin{lemma}
\label{lm:right-side-up}
The chirality element $\ga$ of~$C$ may be chosen so that
$\rho_*(\ga) = \Ga$.
\end{lemma}

\begin{proof}
In the representation $C$ of $\Cliff(T^*X, g)$, the chirality element
is given by
$$
\rho_*(\ga)
:= i^{\piso{(p+1)/2}} \sum_\al \sqrt{\det g} \,
\rho_*(da^1_\al) \dots \rho_*(da^p_\al)
= i^{\piso{(p+1)/2}} \sum_\al \sqrt{\det g} \,
[\D,a^1_\al] \dots[\D,a^p_\al].
$$
In $C$, skewsymmetrization of $i^{-\piso{(p+1)/2}}\Ga$ and of
$i^{-\piso{(p+1)/2}} \rho_*(\ga)$ yields nonvanishing sections of
$\La^p T^*X$. Consequently, $\Ga = f\,\rho_*(\ga)$ for some
nonvanishing real function $f$ on~$X$. Squaring gives
$1 = f^2\,\rho_*(\ga^2) = f^2$, and therefore $f = \pm 1$ since $X$ is
connected. We now fix the orientation on $T^*X$ for $g$ so that
$f = +1$, and thereby $\rho_*(\ga) = \Ga$.
\end{proof}

\section{Poincar\'e duality and spin$^c$ structures}
\label{sec:PD}

This Section uses Poincar\'e duality in $K$-theory to identify the
manifold $X$ as a spin$^c$ manifold.

We shall use the Kasparov intersection product \cite{KasparovTech},
but in fact only require its functorial properties and some results
for products with particular classes. For an executive summary of its
properties, see \cite[IV.A]{Book}. 

In order to identify the spectral triple $(\A,\H,\D)$ with that of the
Dirac operator on an irreducible bundle of spinors, we must in
particular show that the manifold $X = \spec(\A)$ is spin$^c$. By the
work of \cite{Karrer} and~\cite{Plymen}, when $\dim X$ is even this
amounts to the existence of a spinor bundle $S \to X$, carrying a
(pointwise) irreducible representation of the Clifford algebra bundle
$\Cliff(X)$; and likewise in the odd-dimensional case, using instead
irreducible representations of the algebra subbundle
$\Cliff^+(X) := \frac{1+\ga}{2} \Cliff(X)$.

The spin$^c$ condition can also be rephrased as the existence of a
Morita equivalence bimodule between $\Ga(X,\Cliff(X))$ and $C(X)$ if
$\dim X$ is even, respectively between $\Ga(X,\Cliff^+(X))$ and $C(X)$
if $\dim X$ is odd; this bimodule is provided by the sections of an
irreducible spinor bundle. Proposition~\ref{pr:spin-from-PD} below
shows that the existence of a class $\mu \in K^*(A \ox A)$
satisfying Poincar\'e duality in $K$-theory, Condition~\ref{cn:pdual},
implies that the manifold is spin$^c$. We do this by comparing the
Poincar\'e duality isomorphism coming from $\mu$ and the
``Riemannian'' Poincar\'e duality isomorphism described next,
which holds for any compact oriented manifold, spin$^c$ or not.

The class $\la \in KK(\Ga(\Cliff(X)) \ox C(X), \C)$ for a compact
oriented manifold $X$ is defined as follows. Choose any Riemannian
metric $g$, and consider $d + d^*$ on $\H = L^2(\La^\8 T_\C^*X,g)$.
The algebra $C(X)$ acts by multiplication operators $m(f)$ on~$\H$,
and with the phase $V$ of $d + d^*$ gives a Fredholm module $(\H,V)$
for $C(X)$. This module is even, since $d + d^*$ anticommutes with the
grading $\eps$ of differential forms by degree $\bmod 2$.

The Clifford algebra bundle, which is isomorphic as a vector bundle to
the exterior bundle $\La^\8 T^*X$, is likewise $\Z_2$-graded as
$\Cliff(X) = \Cliff^0(X) \oplus \Cliff^1(X)$. We need sections of
$\Cliff(X)$ to graded-commute with $d + d^*$ (respectively, with~$V$)
up to bounded (respectively, compact) operators and to graded-commute
with~$\eps$. For that \cite[Defn.~4.2]{KasparovEqvar}, we define the
action of covectors $v \in T^*X$, as usual in physics, by
$$
\tilde c(v)\,\omega := v \w \omega + i_g(v)\,\omega,
\word{for} \omega \in \La^\8 T_\C^*X,
$$
where $i_g$ denotes contraction with respect to~$g$. We find that
$$
\tilde c(v) \tilde c(w) + \tilde c(w) \tilde c(v) = + 2g(v,w).
$$
This differs from the action of covectors arising from the symbol
of $d + d^*$, since $c(v) := \sg_{d+d^*}(x,v)$ yields
$c(v)\,\omega = v \w \omega - i_g(v)\,\omega$. This action satisfies
$c(v) c(w) + c(w) c(v) = -2g(v,w)$. As is pointed out in
\cite{Polaris}, $\tilde c(v)$ comes from the action of the symbol of
$i(d - d^*)$. While these actions may generate non-isomorphic real
algebras, their complexifications are isomorphic.

The pair $(\H, V)$ is thus a $\Z_2$-graded Fredholm module carrying
the left action $\tilde c \ox m$ of $\Ga(\Cliff(X)) \ox C(X)$ on $\H$
and the trivial right action of~$\C$; we denote its (operator
homotopy) class by $\la \in KK(\Ga(\Cliff(X)) \ox C(X),\C)$. We use
the $KK$ notation to distinguish $\la$ from a $K$-homology class in
$K_0(\Ga(\Cliff(X)) \ox C(X))$ in order to stress that we are dealing
here with $\Z_2$-graded algebras.

It turns out that $KK(\C,\Ga(\Cliff(X)))$ is the Grothendieck group of
equivalence classes of $\Z_2$-graded $\Ga(\Cliff(X))$-modules. By
\cite{LawsonM}, there is a canonical isomorphism with the group of
ungraded modules for $\Ga(\Cliff^0(X))$, which is
$K_0(\Ga(\Cliff^0(X)))$.

The interest in the class $\la$ is the following special case of
\cite[Thm.~4.10]{KasparovEqvar}.

\begin{thm}
\label{th:Kas-iso}
Let $X$ be a compact boundaryless oriented Riemannian manifold. Then
$$
\8 \ox \la: KK^i(\C,\Ga(\Cliff(X))) \to KK^i(C(X),\C) \simeq K^i(C(X))
$$
is an isomorphism.
\qed
\end{thm}

We shall suppose in what follows that $X$ is a compact boundaryless
manifold on which a Riemannian metric $g$ is given, and we write,
somewhat sloppily, $\Cliff(X) := \Cliff(T^*X,g)$ for the corresponding
Clifford algebra bundle over~$X$. We also write
$\B = \Gaoo(X,\Cliff(X))$ for its algebra of smooth sections, $B$ for
the norm completion of~$\B$ (the continuous sections), and $B^0$ for
the even part of~$B$ in the natural $\Z_2$-grading. We also write
$B^\pm = \frac{1\pm\ga}{2} B$ in the odd case, and note that while
$B^+$ is isomorphic to $B^0$, $B^+$ contains both odd and even
elements of~$B$.

There are two useful representations of $B$ on $L^2(\La^*T^*X)$. We
denote by $\theta_-$ the representation coming from the symbol of
$d + d^*$, since $c(v)^2 = -g(v,v)$ for each $v \in T^*X$; and by
$\theta_+$ the representation coming from the symbol of $i(d - d^*)$,
since therein $\tilde c(v)^2 = +g(v,v)$ for all covectors. These two
representations graded-commute.

Let $\la \in KK(B \ox A, \C)$ be the class described above, so the
representation of the $\Z_2$-graded algebra $B$ making this a graded
Kasparov module is~$\theta_+$. We set
$$
\la_A := i_B^*\la \in K^0(A),
$$
where $i_B\: A \to B \ox A : a \mapsto 1_B \ox a$. This $\la_A$ is the
class (see \cite[Sec.~10.9]{HigsonR}) of the spectral triple
$(\Coo(X), L^2(\La^\8 T_\C^*X), d+d^*)$. Thus if $(\E,F)$ is any
Fredholm module representing $\la_A$, there is a unitary operator
$U \: L^2(\La^\8 T_\C^*X) \to \E$ that preserves gradings and an
operator homotopy $\{F_t\}$ such that, modulo degenerate Fredholm
modules,  $F_0 = F$ and
\begin{equation}
(U^*\E, U^* F_1 U) = (L^2(\La^\8 T_\C^*X), V),
\word{where}  d + d^* =: V\,|d + d^*|.
\label{eq:uni-equiv}
\end{equation}
Using $U$, we can transport both graded-commuting representations of
$B$ on $L^2(\La^\8 T_\C^*X)$ to~$\E$, and $\theta_+(B)$ will give
$(\E,F_1)$ the structure of a $(B \ox A)$-Kasparov module.

\begin{lemma}
\label{lm:no-graded-rep}
Let $(\A,\H,\D)$ be a spectral manifold of dimension~$p$. Then if $p$
is even, there is no faithful $\Z_2$-graded (by $\Ga$) representation
of $B$ on $\H$ graded-commuting with $\CDA = \pi_\D(\Omega^\8\A)$. If
$p$ is odd, there is no faithful representation of $B^+$ on $\H$
graded-commuting with~$\CDA$.
\end{lemma}

\begin{proof}
The algebra $B$ (and also the representation $\CDA$ of $\B$) is
$\Z_2$-graded by the parity of the number of $1$-form components in a
product $\omega = b_0\,c(db_1)\dots c(db_k)$, which coincides with the
grading given by $\Ga = \pi_\D(c)$. In the even-dimensional case, we
observe that a representation of $B$ on $\H$ graded-commuting with
$\CDA$ must commute with $\Ga$ since $\Ga$ is an even element of
$\CDA$, whereas any odd element of $B$ must be represented by an 
operator anticommuting with~$\Ga$. Any such representation must kill 
the odd elements of~$B$ and thus cannot be faithful.

In the odd-dimensional case, since $\Ga = 1$ is an odd element of
$\CDA$, because it is a sum of products of an odd number of $1$-forms,
such a representation of~$B$ must graded-commute with $\Ga = 1$, and 
so its $1$-form elements would anticommute with~$1$, which is absurd.
The same remains true for $B^+$ since it contains elements with 
nonzero odd components.
\end{proof}

\begin{rmk}
\label{rk:no-graded-rep}
In fact, the proof shows that even if we can (globally) split the
cotangent bundle into $r$- and $s$-dimensional subspaces with
$r + s = p$, and so write
$\Ga(\Cliff(X)) \simeq \Ga(\Cliff_r(X)) \hatox_A \Ga(\Cliff_s(X))$
(graded tensor product) then neither tensor factor can act
faithfully on $\H$ in such a way as to graded-commute with $\CDA$ and,
in the even case, be graded by~$\Ga$.
\end{rmk}

In all of what follows, the group $KK(\C, B)$ consists of formal
differences of $\Z_2$-graded right $B$-modules, while
$KK(\C, A) = K_0(A)$ consists of formal differences of right
$A$-modules: we consider only even representations of $A = C(X)$. By
\cite[Thm.~5.4]{KasparovTech}, we may regard elements of $K_1(A)$ as
elements of $KK(\Cliff_1, A)$, where $\Cliff_1$ is the
($2$-dimensional) graded complex algebra generated by a single odd
element $\ga$ with $\ga^2 = 1$.

\begin{prop}
\label{pr:spin-from-PD}
If $X$ is a compact boundaryless oriented manifold and $A = C(X)$,
then $X$ has a spin$^c$ structure if and only if there is a class
$\mu \in K^\8(A \ox A)$, represented by a spectral manifold
$(\Coo(X),\H,\D)$, for which
\begin{equation}
x \mapsto x \ox_A \mu : K_\8(A) \to K^\8(A)
\label{eq:AA-match}
\end{equation}
is an isomorphism.
\end{prop}

\begin{proof}
Suppose first that there exists a $\mu \in K^\8(A\ox A)$ such that
\eqref{eq:AA-match} is an isomorphism, and that $\mu$ is represented
by a spectral manifold (regarded as an unbounded Fredholm
module~\cite{HigsonR}). Write $\mu_A := i_A^* \mu \in K^\8(A)$, where
$i_A \: A \to A \ox A : a \mapsto 1_A \ox a$. Since $\8 \ox \la :
KK^i(\C,B) \to K^i(A)$ is an isomorphism by Theorem~\ref{th:Kas-iso},
and \eqref{eq:AA-match} is an isomorphism by hypothesis, there exist
classes $y \in KK^p(\C, B)$ and $x \in K_p(A)$ such that
$$
\mu_A = y \ox_B \la  \word{and}  \la_A = x \ox_A \mu.
$$

We shall treat explicitly below only the even case, employing the
suspension isomorphism $s \: K^{p+1}(C_0(\R) \ox A) \to K^p(A)$ to
handle the odd case. By Theorem~10.8.7 and Proposition~11.2.5
of~\cite{HigsonR}, for any Dirac type operator $\D$ on a manifold
there are identifications
$$
[\Dslash_\R] \ox_\C [\D] = [\Dslash_\R \hatox \D]  \word{and}
s([\Dslash_\R \hatox \D]) = [\D],
$$
where $\Dslash_\R = -i\,d/dx$ is the usual Dirac operator on~$\R$.
The analogous result for the $K$-theory suspension can also be found
in~\cite{HigsonR}.

Suppose, then, that $p$ is even. Let $x$ be represented by $(\E,0)$,
where $\E = \E_1 \oplus \E_2$ is a finitely generated projective
graded $A$-module, and denote the grading by
$\eps = \twobytwo{1}{0}{0}{-1}$. Note that the corresponding
$K$-theory class is $x = [\E_1] - [\E_2] \in K_0(A)$.

We shall represent the product $x \ox_A \mu$ using the unbounded
formalism of~\cite{Kucerovsky}. Using Lemma~\ref{lm:dense-proj}, we
suppose without loss of generality that $\E$ is a finitely generated
projective $\A$-module (rather than $A$-module). Let
$\nabla \: \E \to \E \ox_\A \Omega^1(X)$ be a connection, and define
$$
\hat c : \E \ox_\A \Omega^1(X) \ox_\A \H  \to  \E \ox_\A \H
: s \ox a\,db \ox \xi \mapsto \eps s \ox a\,[\D,b]\xi,
$$
for $s \in \E$, $a,b \in \A$, $\xi \in \H$. Then
$\Dhat := \hat c \circ (\nabla \ox 1) + \eps \ox \D$ is well defined
and essentially selfadjoint on $\E \ox_\A \H_\infty$, and the pair
$(\E \ox_\A \H_\infty, \Dhat)$ is an unbounded representative of the
Kasparov product $x \ox_A \mu$, by \cite[Thm.~13]{Kucerovsky}. 
Moreover, $\Dhat$ is a first-order Dirac-type operator since
$$
[\Dhat,a](s \ox \xi) = \eps s \ox [\D,a]\xi,
\word{for all} s \in \E, \ \xi \in \H_\infty, \ a \in \A.
$$

Now modulo degenerate Kasparov modules, there is a unitary $U$ such
that
\begin{equation}
(U(\E \ox_A \H), U \Dhat U^*) = (L^2(\La^\8T_\C^*X), \widehat{d+d^*}),
\label{eq:kasp-path}
\end{equation}
where the phase of $\widehat{d+d^*}$ is operator homotopic to the
phase of $d + d^*$. Using this unitary we may transport the two
graded-commuting representations $\theta_-$, $\theta_+$ of $B$ to
$\E \ox_A \H$.

Now any (unbounded) operator $\widetilde\D$ representing the product
$x \ox_A \mu$ on the module $\E \ox_\A \H_\infty$ has principal symbol
homotopic to 
$$
\tilde\sg_\D :  da \mapsto \eps \ox [\D,a]
: \Omega^1(X) \to \End_A(\E \ox_A \H).
$$
This makes sense, since we may work on smooth sections and can define
the principal symbol. By Lemma~\ref{lm:no-graded-rep} and
Remark~\ref{rk:no-graded-rep}, no nonzero $1$-form on $X$ can act
nontrivially on $\H$ so that it graded-commutes with $\CDA$ and
anticommutes with $\Ga = \pi_\D(c)$. Hence the representation
$\Ad U^* \circ \theta_+$ acts effectively on the first tensor factor
$\E$ of $\E \ox_A \H$; in this way, $\E$ becomes a left $B$-module.
This action is $A$-linear (since $A$ is central in $B$) and moreover
is adjointable. Indeed, for $s,t \in \E$, $\xi, \eta \in \H$ and
$b \in B$, we find that
\begin{align*}
\braket{\xi}{\pairing{bs}{t}\eta} 
&= \braket{bs \ox \xi}{t \ox \eta} = \braket{b(s \ox \xi)}{t \ox \eta}
\\
&= \braket{s \ox \xi}{b^*(t \ox \eta)}
= \braket{s \ox \xi}{b^*t \ox \eta}
= \braket{\xi}{\pairing{s}{b^*t}\eta} 
\end{align*}
which implies $\pairing{bs}{t} = \pairing{s}{b^*t}$ for the action
of~$B$ on~$\E$. Thus $B$ acts by endomorphisms of the $\Z_2$-graded
vector bundle $E$ for which $\E = \Gaoo(X,E)$.

Now we let $V \subseteq U_\al \subset X$ be an open set over which the
bundle $E$ is trivial, so that
$\E \cdot C_0(V) \simeq C_0(V)^{\rank E}$. Let $j\: C_0(V) \to A$ be
the homomorphism obtained by extending functions by zero, and let
$E_+$, $E_-$ be the subbundles of~$E$ whose sections are 
$\E_\pm := \half(1 \pm \eps)\E$. We claim that
\begin{equation}
j^*(x \ox_A \mu) = (\rank E_+ - \rank E_-)\, j^*\mu.
\label{eq:rank-count}
\end{equation}
Recall that the map $j^*$ just restricts the representation of~$A$ to
$C_0(V)$. To prove \eqref{eq:rank-count}, we begin by splitting the
Hilbert space $\E \ox_A \H$ as
$$
\E \ox_A \H = L^2(V, E \ox S) \oplus L^2(X \less V, E \ox S).
$$
If we let $P$ be the projector whose range is the first summand,
then we can replace $F_{\Dhat} := \Dhat(1 + \Dhat^2)^{-1/2}$ by 
$P F_{\Dhat} P + (1 - P) F_{\Dhat} (1 - P)$. This is valid because for
$a \in C_0(V)$ the relation
$$
a PF_{\Dhat}(1 - P) = P [a,F_{\Dhat}] (1 - P)
$$
holds and defines a compact operator. Now the representation of
$C_0(V)$ on the second summand is zero, and so the straight line path
$$
t \mapsto 
(1 - P)\biggl((1 - t)F_{\Dhat} + t\twobytwo{0}{1}{1}{0}\biggr)(1 - P),
$$
where the matrix $\twobytwo{0}{1}{1}{0}$ is defined relative to the
splitting $\pi_\D(\cc) = \twobytwo{1}{0}{0}{-1}$, gives an operator
homotopy from the second summand to a degenerate Kasparov module.

For the first factor we now observe that
$$
L^2(V, E \ox S) = L^2(V, (X \x \C^{\rank E}) \ox S)
= C_0(V)^{\rank E} \ox_{C_0(V)} L^2(V,S).
$$ 
Similarly, by choosing the trivial connection over $V$, we get
$\Dhat|_V = \hat c \circ (d|_V \ox 1) + \eps|_V \ox \D|_V$, where
$\eps^2 = 1$ on $C_0(V)^{\rank E}$. Thus,
\begin{align*}
[j^*(\E \ox_A \H, \Dhat)]
&= [(C_0(V)^{\rank E}, 0)] \ox_{C_0(V)} [(P\H, PF_\D P)]
\\
&= [(C_0(V)^{\rank E}, 0)] \ox_{C_0(V)} j^*[(\H, F_\D)]
\\
&= [(C_0(V)^{\rank E}, 0)] \ox_{C_0(V)} j^*\mu
\\
&= (\rank E_+ - \rank E_-)\, j^*\mu.
\end{align*}

Now if we choose $V \subset U_\al$ so that the closure of $V$ is
contractible, then every bundle on $X$ restricts to a trivial bundle
on $\ol V$. Then $V$ is (trivially) a spin$^c$ manifold, and so 
$$
j^*\la_A = \la_A\bigr|_V = [S_V^*] \ox_A [\Dslash]_V
= 2^{p/2} \,[\Dslash]_V
$$
where $[\Dslash]_V$ is the class of any Dirac operator on $V$, and
$S_V \to V$ is the (complex) spinor bundle over~$V$. Likewise, from
\eqref{eq:rank-count} we get
$$
j^*\la_A = j^*(x \ox_A \mu) = (\rank E_+ - \rank E_-)\, j^*\mu
= r (\rank E_+ - \rank E_-) \,[\Dslash]_V,
$$
where $r$ is the number of (pointwise) irreducible components of the
representation $\CDA$, obtained in Proposition~\ref{pr:cliff-dirac}.

Comparing these two expressions for $j^*\la_A$, we see that
$\rank E_+ - \rank E_- > 0$. Write
$E|_V = E_1 \oplus E_2 \oplus E_3$, where $\eps = +1$ on
$E_1 \oplus E_2$, $\eps = -1$ on $E_3$, and $E_2 \simeq E_3$ (since
these are trivial bundles over $V$, this just means that $E_2$ and 
$E_3$ have equal rank). Then the explicit formula for the Kasparov 
product shows that on $\E_2 \oplus \E_3$, $\Dhat$ decomposes as
$$
\Dhat = \twobytwo{\hat c_+(d \ox 1) + 1 \ox \D}{0}{0}
{-\hat c_+(d \ox 1) - 1 \ox \D},  \word{graded by}
\twobytwo{1 \ox \Ga}{0}{0}{-1 \ox \Ga},
$$
where $\hat c_+(s \ox a\,db \ox \xi) = s \ox a\,[\D,b]\xi$, i.e.,
just like $\hat c$ but with $\eps$ acting by~$1$. This displays the
restriction of the product to $\E_2 \oplus \E_3$ as the sum of a
Fredholm module for $A$ and its negative. Hence it is homotopic to a
degenerate module, and therefore
$$
j^*(x \ox_A \mu) = [(\E_1 \ox_A \H, \hat c_+(d \ox 1) + 1 \ox \D)]
= r(\rank E_1) \,[D]_V.
$$
Now $j^* \la_A = \la_{C_0(V)}$, by \cite[Prop.~10.8.8]{HigsonR}, and
since the $B$-module structure of $\la$ is locally defined, our
previous arguments apply to show that $\E_1$ carries a faithful
representation of $B$. This immediately implies that
$\rank E_1 \geq 2^{p/2}$, and we conclude that
$$
\rank E_1 = 2^{p/2}  \word{and}  r = 1.
$$

Therefore, the bundle $S$ underlying the module
$\H_\infty \simeq \Gaoo(X,S)$ is not only globally irreducible as a
$B$-module, but is also fibrewise irreducible. Thus, as shown
in~\cite{Plymen}, $S$ is the spinor bundle for a spin$^c$ structure
on~$X$. Our Propositions \ref{pr:cliff-dirac}
and~\ref{pr:global-cliff-dirac} and Corollaries \ref{cr:D-is-Dirac}
and~\ref{cr:D-is-Dirac-again} have established that $\D$ is a Dirac
operator for this spin$^c$ structure and a suitable metric, up to the
addition of an endomorphism of the spinor bundle. Hence
$\mu = [(\A,\H,\D)]$ is the fundamental class of the spin$^c$ 
manifold~$X$.

Conversely, if $X$ carries a spin$^c$ structure, we may construct the
spectral triple of any Dirac operator on its complex spinor bundle. By
\cite{BaumD,Book,HigsonR}, this class satisfies Poincar\'e duality in
$K$-theory.
\end{proof}

\begin{rmk}
The above proof uses only a little of the structure of the spectral
manifold representing $\mu$: we need the first-order property, we
require $\Ga$ to provide the grading on $\H$, and we need the
equality $\Ga = \pi_\D(c)$.
\end{rmk}

We summarize our results by stating the following theorem.

\begin{thm}
\label{th:Poincare-dual}
Let $(\A,\H,\D)$ be a spectral manifold of dimension $p$, and suppose
that $\mu := [(\A \ox \A, \H, \D)] \in K^\8(A \ox A)$ satisfies
Condition \ref{cn:pdual}. Then $X$ is a spin$^c$ manifold, $\H_\infty$
is the module of smooth sections of its complex spinor bundle, the
representation $\CDA$ of the associated Clifford algebra bundle is
irreducible, and thus $\D$ is, up to adding an endomorphism, a Dirac
operator on this spinor bundle. The class of $(\A,\H,\D)$ is thus the
fundamental class of the spin$^c$ manifold $X = \spec(\A)$ with the
spin$^c$ structure defined by the irreducible representation~$\CDA$.
\end{thm}

\begin{proof}
The assertion that $X$ is spin$^c$ follows from
Proposition~\ref{pr:spin-from-PD}, under the assumption that the
$K$-homology class of $(\A,\H,\D)$ in $K^\8(A\ox A)$ gives a
Poincar\'e duality isomorphism.

Moreover, Proposition~\ref{pr:global-cliff-dirac} establishes that
$\CDA$ is a representation on $\H_\infty = \Gaoo(X,S)$ of the
algebra $\Gaoo(X,C)$ for the complex Clifford-algebra bundle
$C = \Cliff(T^*X,g)$, for a specific Riemannian metric $g$ on~$X$.

Thus $\H_\infty$ is a spinor module, $S$ is the corresponding complex
spinor bundle, and up to the addition of an $\A$-linear endomorphism
of $\H_\infty$, $\D$ is the Dirac operator for this metric and spinor
bundle.
\end{proof}

\begin{corl}
\label{cr:measure-recovery}
The functional $\mu_\D$ is given by the Wodzicki residue:
$$
\mu_\D(a) = \ncint a \,\Dreg^{-p}
= \frac{N\Vol(\Sf^{n-1})}{p(2\pi)^p} \int_X a(x) \,d\vol_g.
\eqno \qed
$$
\end{corl}

\begin{rmk}
The only condition not used at all so far is the ``reality''
requirement, Condition~\ref{cn:real}. Since we have now at our
disposal a spin$^c$ structure, we only need to refine it to a spin
structure. For that, we refer to \cite{Plymen} and \cite{Polaris},
wherein it is shown that the spinor module for a spin structure is
just the spinor module for a spin$^c$ structure equipped with
compatible charge conjugation that is none other than the real
structure operator $J$ (acting on $\H_\infty$); and to
\cite{RennieComm}, where the spin structure is extracted, using~$J$,
from a representation of the real Clifford algebra of~$T^*X$. Thus, by
additionally invoking Condition~\ref{cn:real}, we may replace each
mention of `spin$^c$' by `spin' in the statement of 
Theorem~\ref{th:Poincare-dual}.
\end{rmk}

\section{Conclusion and outlook}
\label{sec:more-conds}

We close with several remarks on the hypotheses of the reconstruction
theorem, and some possible variants.

In \cite{RennieComm} and in \cite{Polaris}, under the additional
assumption that $\A = \Coo(X)$, which is now redundant, the Dirac
operator for the spinor bundle $S$ and the metric~$g$ is also
recovered from the given spectral triple, by minimizing the action
functional described in \cite{ConnesGrav}: see also \cite{Portia} for
that. Thus most of the statements of \cite{ConnesGrav} can be
recovered, if need be.

Poincar\'e duality in $K$-theory, expressed here as
Condition~\ref{cn:pdual}, enters our proof only to show that the
manifold carries a (preferred) spin$^c$ structure, and identify the
fundamental class of that spin$^c$ structure as the $K$-theory class
of a spectral triple. We can only make this identification once we
know that $\spec(\A)$ is a manifold.

It would be more economical to replace Condition~\ref{cn:pdual} with
the following one.

\begin{cond}[Poincar\'e duality II: Morita equivalence]
\label{cn:pmorita}
The $C^*$-module completion of $\H_\infty$ is a Morita equivalence
bimodule between $A$ and the norm completion of~$\CDA$.
\end{cond}

There are good reasons to prefer this Morita equivalence condition.
First, we could have employed it from the beginning to simplify the
structure of $\CDA$. This would slightly simplify the task of building
the manifold.

Also, if we wished to propose axiomatics for different kinds of
manifold: Riemannian, complex, K\"{a}hler, and so on, the Morita
equivalence axiom can be easily replaced with an appropriate
characterization of the behaviour of the Clifford algebra. Together
with specifying the behaviour of the volume form, $\pi_\D(c)$, and
adapting the reality condition, these axiomatics should be flexible
enough to cope with different types of manifolds. A framework for this
has been suggested in~\cite{FroehlichGR}.

Our proof may be adapted to deal with less smooth cases. If we start
from a $QC^k$ spectral triple, $k \in \N$ (see \cite{CareyPRSone} for
instance), the completion of $\A$ in the $QC^k$-topology is a Banach
algebra stable under holomorphic functional calculus. The $\Coo$
functional calculus will still work, as would a $C^k$ functional
calculus. The Lipschitz functional calculus so crucial to our results
only requires $QC^0$.

Additional smoothness (beyond $QC^0$) is needed for the following
items. To employ the Hochschild class of the Chern character to deduce
measurability, $QC^{\max(2,p-2)}$ is needed, although we do not
actually use this in our proof. To deduce that the noncommutative
integral given by the Dixmier trace defines a hypertrace on $\A$ needs
only $QC^0$, but to obtain a trace on $\CDA$ requires $QC^2$
\cite{CiprianiGS}; the question is further explored in
\cite[Sec.~6]{CareyRSS}. Finally, Lemma~\ref{lm:metric-smooth}
requires more than $QC^0$. Observe that the smoothness used in the
unique continuation arguments is a manifestation of ellipticity on
$\R^p$; we do not need to assume regularity of $\A$ for that, the
smoothness we use is that of the eigenvectors of an elliptic operator.

Reformulating the proof to deal with less smoothness would be somewhat
inelegant, as the metric $g(da,db)$ need only take values in bounded
measurable functions, so we would need to work with more algebras,
such as the bicommutant $A''$. Nevertheless, it appears feasible. In
this vein, one should replace $\H_\infty$ with $\Dom\D$, and similarly
modify the finiteness condition.

It is fascinating to reflect on just what input gave us a manifold.
While of course we needed everything, the role of the Dirac operator
is crucial in producing the coordinate charts. It would be of some
interest to understand this in the language of geometric topology.

Geometrically, the weak unique continuation property may be regarded
as a consequence of the local splitting $\D = n(\del_\xi + A)$ for any
direction~$\xi$ \cite{BoossW}. Thus Dirac operators know how to split
any neighbourhood into a product of a hypersurface and its normal.

The strong unique continuation property is better understood by
analogy with holomorphic functions and Cauchy--Riemann operators.
Essentially, Dirac type operators have ``sufficiently rigid''
eigenvectors to determine the local geometry in the neighbourhood of a
point.

The closedness property, our Condition~\ref{cn:closed}, does not
appear explicitly in the axiom scheme proposed by Connes
in~\cite{ConnesGrav}. However, it played a critical role in earlier
formulations of Poincar\'e duality at the level of chains, needed to
address colour symmetry in the Standard Model \cite[VI.4.$\ga$]{Book}
and we hope to have exemplified that it is still a useful tool.

The relation of this closedness condition with
Condition~\ref{cn:pdual}, that of Poincar\'e duality in $K$-theory,
deserves further scrutiny. We leave open the question of whether a
manifold could be reconstructed if one chooses to drop
Condition~\ref{cn:closed} in favour of $K$-theoretic Poincar\'e
duality alone. We suspect not; on the other hand, neither do we have a
counterexample.

\appendix

\section{Hermitian pairings on finite projective modules}

The finiteness condition calls for a Hermitian inner product on a
finite projective module with specific properties. We investigate the
nature of such Hermitian pairings here.

There is a small question of starting points to deal with. If we take
as given the scalar product on the Hilbert space, and try to find a
Hermitian pairing to satisfy the finiteness condition, we must tackle
a difficult existence problem. On the other hand, if we suppose that
there is a Hermitian pairing related to the scalar product as in
Condition~\ref{cn:finite} above, then we face a question of
characterisation. Here we adopt this latter point of view.

\begin{lemma}
\label{lm:inv-extend}
Let $A$ be a unital $C^*$-algebra and let $q \in M_m(A)$
be a projector, $q = q^* = q^2$. Suppose that $M \in q M_m(A) q$.
Then $M$ is invertible in $q M_m(A) q$ if and only if 
$M + 1 - q$ is invertible in $M_m(A)$.
\end{lemma}

\begin{proof}
If $M$ is invertible in $q M_m(A) q$ with inverse $M^{-1}$ satisfying 
$MM^{-1} = M^{-1}M = q$, then $M^{-1} + 1 - q$ is inverse to
$M + 1 - q$ in $M_m(A)$.

Conversely, suppose that $M + 1 - q$ is invertible, with inverse
$N$ in $M_m(A)$. Then $qNq$ is inverse to $M$ in $q M_m(A) q$, since
$N(M + 1 - q) = 1$ implies $q = NM = qNMq = qNqM$, and similarly
$(M + 1 - q)N = 1$ implies $q = MN = qMNq = MqNq$.
\end{proof}

\begin{lemma}
\label{lm:matr-pairing}
Let $A$ be a unital $C^*$-algebra, and $E = qA^m$ a finite projective
right $A$-module. Then every Hermitian pairing on $E$ is of the form
$$
\pairing{e}{f} = \sum_{j,k} e^*_j M_{jk} f_k
$$
where $e = (e_1,\dots,e_m)^T$ with each $e_j \in A$, $qe = e$;
similarly for $f$; and $M \in qM_m(A)q$ is positive.
\end{lemma}

\begin{proof}
Denote by $\eps \in A^m$ the column vector with $1$ in the $j$-th
spot, and zeroes elsewhere. We let
$x_j = q\eps_j = \sum_k q_{kj} \eps_k$, so we can write all $e \in E$
as $e = \sum_j x_j e_j$.

If $\pairing{\cdot}{\cdot}$ is a Hermitian pairing on $E$, then for
any $e \in E$
$$
0 \leq \pairing{e}{e} = \sum_{j,k} e_j^* \pairing{x_j}{x_k} e_k.
$$
Hence the matrix $M = [\pairing{x_j}{x_k}]_{jk} \in M_m(A)$ is
positive, by \cite[Proposition 1.20]{Polaris}. Next
\begin{align*}
(qM)_{mk} &= \sum_{j} q_{mj} \pairing{x_j}{x_k}
= \sum_{j} q_{mj} \bigpairing{\sum_lq_{lj} \eps_l}{x_k}
= \sum_{j,l} q_{mj} \pairing{\eps_l q_{lj}}{x_k}
\\
&= \sum_{j,l} q_{mj}q_{jl} \pairing{\eps_l }{x_k}
= \sum_{l} q_{ml} \pairing{\eps_l }{x_k}
= \pairing{x_m}{x_k} = M_{mk},
\end{align*}
and similarly $Mq = M$.
\end{proof}

\begin{lemma}
\label{lm:full-invertible}
Let $A$ be a unital $C^*$-algebra, and $E = qA^m$ a finite projective
right $A$-module. Let $M \in q M_m(A) q$ be positive, so that
\begin{equation}
\pairing{e}{f}_M := \sum_{j,k} e^*_j M_{jk} f_k,
\word{for all}  e,f \in E,
\label{eq:matr-pairing}
\end{equation}
defines a Hermitian pairing on $E$. Then the right~$C^*$ $A$-module
$\bigl(E, \pairing{\cdot}{\cdot}_M \bigr)$ is full if and only if 
$M$ is invertible in $q M_m(A) q$.
\end{lemma}

\begin{proof}
First suppose that the $C^*$~$A$-module 
$\bigl(E, \pairing{\cdot}{\cdot}_M \bigr)$ is full.
By \cite{RaeburnW}, the compact 
endomorphisms of the full right $A$-module $E$ fulfil
$\End_A^0(E) = \End_A(E) = q M_m(A) q$. This algebra is generated by 
the rank-one operators $\Theta^M_{e,f} \: E \to E$ defined by
$$
\Theta^M_{e,f}(g) := e \pairing{f}{g}_M = \Theta^q_{e,f}(Mg),
$$
where $\pairing{e}{f}_q := \sum_{j,k} e^*_j q_{jk} f_k$. The operators
$\Theta^q_{e,f}$ generate $q M_m(A) q$ as an algebra. If $M$ is not
invertible in $qM_m(A)q$, the operators $\Theta^q_{e,f}M$ and their
adjoints generate a proper two-sided ideal of it, contradicting
fullness.

Conversely, suppose that $M$ is invertible but that $E$ is not full.
Let $I$ be the closure of $\pairing{E}{E}_M$ in~$A$, a two-sided
ideal, so that $E$ is a full right $I$-module, and thus
$\End_A^0(E) = \End_I^0(E) = q M_m(I) q$. Then the algebra generated
by all $\Theta^M_{e,f} = \Theta^q_{e,f} M$ is $q M_m(I) q$. But the
invertibility of~$M$ entails that the $\Theta^q_{e,f} M$ generate all
of $q M_m(A) q$: contradiction.
\end{proof}

\begin{lemma}
\label{lm:smooth-pairing}
Let $\A$ be a dense pre-$C^*$-subalgebra of the unital $C^*$-algebra
$A$ with $1 \in \A$. If $\pairing{\cdot}{\cdot}$ is an $\A$-valued
Hermitian pairing on $\E = q\A^m$ making $\E$ full, with
$q = q^* = q^2 \in M_m(\A)$, then it is given by $M \in qM_m(\A)q$ as
in Lemma~\ref{lm:matr-pairing}.
\end{lemma}

\begin{proof}
By \cite[Lemma 2.16]{RaeburnW}, the Hermitian form
$\pairing{\cdot}{\cdot}$ on $\E$ has a canonical extension to the
completion $E = qA^m$. By Lemma~\ref{lm:matr-pairing}, this extension
is determined by a positive invertible $M \in qM_m(A)q$. Since
$\pairing{x_j}{x_k} \in \A$ for all $j,k$, we obtain 
$M \in qM_m(\A)q$.

Observe also that by using the stability under the holomorphic
functional calculus of $\A$, we find $M^{-1} \in qM_m(\A)q$ also.
\end{proof}

\begin{prop}
\label{pr::matr-pairing}
Let $(\A,\H,\D)$ be a spectral triple satisfying
Condition~\ref{cn:metr-dim},
Conditions~\ref{cn:qc-infty}--\ref{cn:first-ord} and
Condition~\ref{cn:irred} of subsection~\ref{ssc:geom-cond} (dimension,
regularity, finiteness, absolute continuity, first order,
irreducibility). Then its Hermitian pairing $\pairing{\cdot}{\cdot}$
on $\H_\infty = q\A^m$ is (a positive multiple of) the standard one.
\end{prop}

\begin{proof}
Let $M \in qM_m(\A)q$ be the positive invertible element such that
$\pairing{\xi}{\eta} = \pairing{\xi}{\eta}_M$ for
$\xi,\eta \in \H_\infty$. Then, for each $a \in \A$,
$$
\braket{\xi}{a\eta} = \braket{a^*\xi}{\eta}
= \ncint \pairing{a^*\xi}{\eta}_M \,\Dreg^{-p}
= \ncint \pairing{\xi}{M^{-1} a M\eta}_M \,\Dreg^{-p}
= \braket{\xi}{M^{-1} a M\eta}.
$$
Hence $[M,a] = 0$ for all $a \in \A$. 

Now since $\D$ is a selfadjoint operator on~$\H$, we obtain, for
$\xi,\eta \in \H_\infty$:
\begin{align}
0 &= \braket{\D\xi}{\eta} - \braket{\xi}{\D\eta}
= \ncint \bigl( \pairing{\D\xi}{\eta}_M 
- \pairing{\xi}{\D\eta}_M \bigr) \,\Dreg^{-p}
\nonumber \\
&= \ncint \bigl( \pairing{\D\xi}{M\eta}_q 
- \pairing{\xi}{M\D\eta}_q \bigr) \,\Dreg^{-p}
=: \bbraket{\D\xi}{M\eta} - \bbraket{M\xi}{\D\eta}
\label{eq:comm-quadform}
\end{align}
where $\bbraket{\xi}{\eta} := \braket{M^{-1}\xi}{\eta}$ defines a new
Hilbert space scalar product. Since $M^{-1}$ is bounded with bounded
inverse, this new scalar product $\bbraket{\cdot}{\cdot}$ is
topologically equivalent to the old one $\braket{\cdot}{\cdot}$, and 
so $\H$ is the completion of $\H_\infty$ with respect to either 
scalar product.

Now the right hand side of \eqref{eq:comm-quadform} is the quadratic
form defining the commutator $[\D,M]$ with respect to the scalar
product $\bbraket{\cdot}{\cdot}$. It vanishes on $\H_\infty$ and thus
$[\D,M] = 0$. Now the irreducibility condition implies that $M$ is (a
positive multiple of) the identity $q$ in $qM_m(\A)q$, represented by
a scalar operator on~$\H$.
\end{proof}

\section{Another look at the geometric conditions}

In this Appendix we consider the potential redundancy of, and the
relations between, our metric and irreducibility postulates,
Conditions~\ref{cn:metric} and~\ref{cn:irred}.

Our Condition~\ref{cn:metric} is unnecessarily strong as stated, since
the boundedness of the set
$\set{\eta(a) \in \A/\C\,1 : \|[\D,a]\| \leq 1}$ in the Banach space
$A/\C\,1$ ensures that the distance formula \eqref{eq:metric} not only
defines a metric, but one for which $X$ has finite diameter. However,
the reconstruction of a manifold requires merely a metric, with its
corresponding topology.

Indeed, with no additional assumptions, it is easy to see that the
formula~\eqref{eq:metric} defines a function $d$ satisfying all the
properties of a metric distance, except that there might exist states
$\phi$, $\psi$ for which $d(\phi,\psi) = \infty$.

Thus we could (provisionally) adopt the following condition:
\begin{quote}
For all states $\phi$, $\psi$ of~$A$, the set
$\set{|\phi(a) - \psi(a)| : \|[\D,a]\| \leq 1}$ is bounded.
\end{quote}
Our original Condition~\ref{cn:metric} asks for a uniform bound for
all $\phi,\,\psi$. As it turns out, this condition can be dispensed
with altogether.

\begin{prop}
\label{pr:metric-exists}
Let $(\A,\H,\D)$ be a spectral triple such that $\A$ is unital, $\A\H$
is dense in $\H$, and $\A$ has separable norm closure~$A$. Assume that
$[\D,a] = 0$ if and only if $a = \la\,1$ for some $\la \in \C$. Then
the formula
$$
d_\D(\phi,\psi) := \sup\set{|\phi(a) - \psi(a)| : \|[\D,a]\| \leq 1}
$$
defines a metric on the state space of~$A$.
\end{prop}

\begin{proof}
We need only show that $d(\phi,\psi) = \infty$ cannot occur for any
pair of states $\phi$, $\psi$.

Without loss of generality, we replace $\A$ by its completion in the
norm $a \mapsto \|a\| + \|[\D,a]\|$, since this will not change the
definition of the metric or the norm closure.

Suppose, then, that there is a pair of  states $\phi,\psi$ with
$d(\phi,\psi) = \infty$. There exists (using the separability of~$A$)
a sequence $\{a_N\}_{N\geq 1} \subset \A$ such that
$$
|\phi(a_N) - \psi(a_N)| > N  \word{and}  \|[\D,a_N]\| \leq 1
\quad\text{for all } N.
$$
On replacing $a_N$ by $a_N - \psi(a_N)\,1$ if necessary, we may assume
that $\psi(a_N) = 0$, and therefore
\begin{equation}
\|a_N\| \geq |\phi(a_N)| > N  \word{for all}  N \in \N.
\label{eq:aN-large}
\end{equation}
Set $u_N := a_N/\|a_N\|$, so that $\|u_N\| = 1$ for all~$N$.

Let $\xi \in \H$ and observe that $\D$ has compact resolvent: the
definition of spectral triple assumes $a(\D - \la)^{-1}$ is compact
for all $a \in \A$ and $\la \notin \spec(\D)$, and since the
representation of $\A$ is nondegenerate, $1 \in \A$ acts as the
identity on~$\H$. Choose and fix $\la \in \C \less \spec(\D)$. The
sequence $\{(\D - \la)^{-1} u_N \xi\}$ has a Cauchy subsequence, by
compactness of $(\D - \la)^{-1}$ and boundedness of the sequence
$\{u_N\xi\}$. Each $u_N$ lies in $\Dom \d$, so
$$
(\D - \la)^{-1} u_N \xi
= -(\D - \la)^{-1} [\D, u_N] (\D - \la)^{-1} \xi
  + u_N (\D - \la)^{-1} \xi.
$$
Now $[\D, u_N] \to 0$ in norm, so there is a subsequence $\{u_{N_j}\}$
such that $\{u_{N_j}(\D - \la)^{-1} \xi\}$ is Cauchy. Since every
$\zeta \in \Dom\D$ is of the form $\zeta = (\D - \la)^{-1} \xi$ for
some $\xi \in \H$, we see that for every $\zeta \in \Dom\D$ there is a
Cauchy subsequence $\{u_{N_j}\zeta\}$.

Let $\{\zeta_k\}_{k\geq 1}$ be an orthonormal basis of $\H$ 
consisting of eigenvectors of $\D$. By applying the above argument 
successively to $\xi_k := (\D - \la)\zeta_k$ for $k = 1,\dots,m$,
we can extract from $\{u_N\}$ a subsequence $\{u^{(m)}_j\}$ such that
$\{u^{(m)}_j\zeta_k\}_{j\geq 1}$ converges, for $k = 1,\dots,m$.
Inductively, we can choose a subsequence $\{u^{(m+1)}_j\}$ of
$\{u^{(m)}_j\}$ so that $u^{(m+1)}_j = u^{(m)}_j$ for $j \leq m$
and $\{u^{(m+1)}_j\zeta_{m+1}\}$ converges, too. In the end, we obtain
a subsequence $\{u^{(\infty)}_j\}$ of $\{u_N\}$ such that 
$u^{(\infty)}_j\zeta_k$ converges in $\H$ as $j \to \infty$, for each
basic eigenvector $\zeta_k$. We rename this subsequence to $\{u_N\}$,
noting that \eqref{eq:aN-large} still holds.

Note that the resulting subsequence cannot terminate, since if
$u_n = u_m$ for all $m,n \geq N_0$, then $a_m = \la_{mN_0} a_{N_0}$
for $n \geq N_0$, and
$\|[\D,a_m]\| = |\la_{mN_0}|\,\|[\D,a_{N_0}]\| \leq 1$, implying
$|\la_{mN_0}| \leq \|[\D,a_{N_0}]\|^{-1}$. But then
$$
|\phi(a_m)| = |\la_{mN_0}|\,|\phi(a_{N_0})|
\leq \frac{|\phi(a_{N_0})|}{\|[\D,a_{N_0}]\|},
$$
contradicting the unboundedness \eqref{eq:aN-large} of the sequence
$|\phi(a_m)|$.

If $\zeta$ is an eigenvector of $(\D - \la)^{-1}$ with eigenvalue
$\rho$, then
\begin{equation}
(\D - \la) u_N \zeta = [\D, u_N]\,\zeta + \rho^{-1} u_N \zeta.
\label{eq:domain-cvgce}
\end{equation}
The right hand side converges in~$\H$ and $[\D, u_N] \to 0$ in norm,
so $\lim_N u_N\zeta$ lies in $\Dom\D$ and is an eigenvector of 
$\D - \la$ with eigenvalue $1/\rho$.

Let $P_\rho$ denote the (finite rank) projector whose range is the
eigenspace of~$\D$ for the eigenvalue $\la + \rho$. Then
$$
\lim_{N\to\infty} \|P_\rho u_N P_\rho - u_N P_\rho\| = 0,
$$
the norm limit being appropriate since $P_\rho\H$ is
finite-dimensional. Let $V_\rho := \lim_N P_\rho u_N P_\rho$, a
finite-rank operator. Now \eqref{eq:domain-cvgce} entails
$$
\rho^{-1} V_\rho \zeta = \rho^{-1} \lim_N u_N\zeta
= (\D - \la)\lim_N u_N\zeta 
= \lim_N [\D, u_N]\zeta + \rho^{-1} \lim_N u_N\zeta,
$$
and thus $V_\rho\zeta - \lim_N u_N\zeta = \rho \lim_N [\D,u_N]\zeta$ 
for $\zeta \in P_\rho\H$. Therefore, if
$v := \bigoplus_{\rho\in\spec(\D-\la)} V_\rho$, then
$$
\lim_N \|u_N - v\| = \lim_N \biggl\|
\bigoplus_{\rho\in\spec(\D-\la)} \rho [\D, u_N] P_\rho \biggr\| = 0,
$$
thus the sequence $\{u_N\}$ is actually norm-convergent in $\L(\H)$.

In particular, since $[\D, u_N] \to 0$, the sequence $\{u_N\}$ is
Cauchy in the norm $a \mapsto \|a\| + \|[\D,a]\|$. Hence
$v = \lim_N u_N$ lies in~$\A$. Again using \eqref{eq:domain-cvgce}, we
see that $[\D, v] = \lim_N [\D, u_N] = 0$. Our hypothesis that
$[\D,a] = 0$ only for $a \in \C\,1$ now implies that $v =: \sg\,1$ is
a scalar.

We have already chosen the several $a_N$ (and so also the $u_N$) so
that $\psi(a_N) = \psi(u_N) = 0$. Therefore,
$\sg = \lim_N \psi(u_N) = 0$; we conclude that $u_N \to 0$ in norm.
However, each $u_N$ has norm~$1$, so we have arrived at a
contradiction.

In fine, $d_\D(\phi,\psi) < \infty$ holds for all $\phi$, $\psi$; thus
$d_\D$ is a metric on the state space of~$A$.
\end{proof}

Consider now the irreducibility, as given by Condition~\ref{cn:irred}.
This condition is imposed to guarantee a connected spectrum (in the
commutative case), and to ensure that the Hilbert space is not too
large.

However, already in the commutative case this irreducibility condition
has two aspects. The first is the connectedness of the spectrum: see
Lemma~\ref{lm:no-proj} and Corollary~\ref{cr:spec-conn}. The second is
the irreducibility of $\H_\infty$ as a $\CDA$-module, best controlled
using our proposed Condition~\ref{cn:pmorita} (the Morita-equivalence
version of the spin$^c$ condition). Inspection of the proof shows that
in the reconstruction of a manifold we invoke irreducibility both to
specify the Hermitian pairing and to deduce that the spectrum is
connected, with no nontrivial projectors existing in $\A$ or~$A$,
in the cited Lemma and Corollary. (Once the manifold has been found,
we invoke it once more to get irreducibility of $\H_\infty$ via
Proposition~\ref{pr:no-endo-proj}.)

A closer look shows that, without the irreducibility condition, the
rest of the proof holds with occasional adjustments to take account of
the possible nonconstancy of the rank of the bundle $S \to X$ (which
of course remains locally constant), provided we can guarantee the
existence of the metric defined by~$\D$.

A disconnected spectrum would be a disjoint union of closed components
$X = \biguplus_j X_j$. If there are only finitely many components
$X_j$, then they are clopen (i.e., both closed and open). However, if
there be infinitely many components, they would all be clopen if and
only if the space $X$ is locally connected. When $X$ is compact, this
of course forces the number of components to be finite.

If some component is not clopen, then its characteristic function is
not continuous and so does not lie in the $C^*$-algebra $C(X)$. In the
given (faithful, nondegenerate) representation $\pi$ of $C(X)$, we can
decompose $\Id_\H = \pi(1_{C(X)}) = \sum_j \pi(p_j)$ for a sum of
projectors $p_j\in C(X)$. Then this sum is finite by the separability
of $C(X)$. Hence there can only be finitely many clopen components of
$X$, say $X_1,\dots,X_n$, and $\bigcup_{j>n} X_j$ is also clopen.

Since $\D$ is a local operator, by Corollary~\ref{cr:local-opr}, the
(continuous) characteristic functions of the clopen components lie in
the domain of $\d = \ad\D$ in~$A$. A little more thought shows that
these continuous characteristic functions actually lie in~$\A$.
Therefore, \emph{provided we can ensure the existence of the metric
defined by $\D$ on each clopen piece of the decomposition of $X$}, we
may proceed as follows:
\begin{enumerate}
\item
On each clopen subset of $X$ with continuous characteristic function,
run the existing argument (with modifications for nonconstant rank
of~$S$) to deduce that this clopen piece is a topological manifold;
\item
observe that topological manifolds are locally connected (indeed,
locally path connected);
\item
deduce that all connected components of $X$ are clopen;
\item
conclude that there are only finitely many components, and that $X$ is
a disjoint union of finitely many connected compact topological
manifolds.
\end{enumerate}

Since we should be able to deal with disconnected manifolds as well as
connected ones, we propose the following alternative to
Condition~\ref{cn:irred}.

\begin{cond}[Connectivity]
\label{cn:conn}
There is an orthogonal family of projectors $p_j \in \A$ such that
$1_\A = \sum_j p_j$ and
$$
(a \in \A \word{with} [\D,a] = 0)  \iff a = \tsum_j \la_j p_j
\quad\text{for some } \{\la_j\} \subset \C.
$$
\end{cond}

The separability of $A$ guarantees that these sums are finite, so
convergence questions are moot. We may assume the family $\{p_j\}$ to
be maximal. The original spectral triple breaks up as a finite direct
sum of spectral triples:
$$
(\A,\H,\D) = \bigoplus\nolimits_j (p_j\A p_j, p_j\H, p_j\D p_j),
$$
each of which satisfies the remaining conditions. By
Proposition~\ref{pr:metric-exists} and Condition~\ref{cn:conn}, each
$p_j \A p_j$ has a connected spectrum carrying a metric defined by
$p_j \D p_j$. The procedure above allows us to deduce that
$X = \spec(\A)$ is a disjoint union of finitely many compact 
connected manifolds~$X_j$.

We may then introduce Morita equivalence, as in
Condition~\ref{cn:pmorita}, to control the irreducibility of the
vector bundle $S \to X$ as a $\CDA$-module. Write
$\A_j := p_j\A p_j$, $\H_{\infty,j} := p_j\H_\infty$ and
$\D_j := p_j\D p_j$, and note that $\A$ acts block-diagonally, that 
is, $\H_{\infty,j} \cdot \A = \H_{\infty,j} \cdot \A_j$. Then, since
$\End_\A(\H_\infty)$ consists of all adjointable linear operators on
$\H_\infty$ commuting with the right action of~$\A$, we see that
$$
\End_\A(\H_\infty)
= \End_\A\Bigl( \bigoplus\nolimits_j \H_{\infty,j} \Bigr)
= \bigoplus\nolimits_{i,j} \Hom_\A(\H_{\infty,i}, \H_{\infty,j})
\supseteq \bigoplus\nolimits_j \End_{\A_j}(\H_{\infty,j}).
$$
(If preferred, one may restate these relations using the norm
completions.) 

The Morita equivalence condition amounts to $\CDA$ being generated
by rank-one endomorphisms $\Theta_{\xi,\eta}$, for
$\xi,\eta \in \H_\infty$; and $\Theta_{\xi,\eta} = 0$ if
$\xi \in \H_{\infty,j}$ and $\eta \in \H_{\infty,k}$ with $j \neq k$.
Hence the algebra generated by the rank-one operators is block
diagonal with respect to the direct sum decomposition, and therefore
$$
\CDA = \bigoplus\nolimits_j \End_{\A_j}(\H_{\infty,j})
= \bigoplus\nolimits_j \mathcal{C}_{\D_j}(\A_j).
$$
Thus the Morita equivalence condition for $(\A,\H,\D)$ implies the
same condition for each $(\A_j,\H_j,\D_j)$.

Let $T \in \B(\H_j)$ (strongly) commute with $\D_j$ and $\A_j$, for
some~$j$. Then $T$ preserves the domain of $\D_j$ (since it lies in
the commutant of the von Neumann algebra generated by $\A_j$ and the
spectral projectors of $\D_j$), and so maps $\H_{\infty,j}$ to itself.
Since $T$ commutes with the functions in $\A_j$, it is local and so is
an endomorphism of the bundle $S\bigr|_{X_j}$, and from Morita
equivalence these endomorphisms are precisely the elements of the norm
closure of~$p_j \CDA p_j$. Hence $T$ is a central element of this norm
closure, and since $p_j \CDA p_j$ is an irreducible representation of
a complex Clifford algebra, $T$ is a function in $C(X_j)$. Since $T$
preserves $\H_\infty$, it is smooth and thus $T \in \A_j$; and
$[\D_j,T] = 0$ now implies that $T$ is constant on~$X_j$.

Hence over each connected component of $X$ the irreducibility
condition holds automatically when we assume both the connectivity and
Morita equivalence conditions.

\end{document}